\documentclass{amsart}

\usepackage{amssymb}
\usepackage{amsmath}
\usepackage{amsthm}
\usepackage{amsfonts}

\usepackage{a4wide}
\usepackage{graphicx}
\usepackage[usenames,dvipsnames]{xcolor}

\newtheorem{theorem}{Theorem}
\newtheorem{proposition}[theorem]{Proposition}

\newtheorem{conjecture}[theorem]{Conjecture}

\theoremstyle{definition}
\newtheorem{definition}[theorem]{Definition}
\newtheorem{remark}[theorem]{Remark}
\newtheorem{example}[theorem]{Example}
\newtheorem{hypothesis}[theorem]{Hypothesis}
\newtheorem{counterexample}[theorem]{Counterexample}

\numberwithin{equation}{section}

\begin{document}


\title
[Extremal root paths of Schur \(\sigma\)-groups]
{Extremal root paths of Schur \(\sigma\)-groups \\ and first \(3\)-class field towers with four stages}

\author{Daniel C. Mayer}
\address{Naglergasse 53 \\ 8010 Graz \\ Austria}
\email{algebraic.number.theory@algebra.at}
\urladdr{http://www.algebra.at}


\thanks{Research
supported by the Austrian Science Fund (FWF): P 26008-N25}


\date{April 10, 2020}


\begin{abstract}
An extremal property of finite Schur \(\sigma\)-groups \(G\)
is described in terms of their path to the root
in the descendant tree of their abelianization \(G/G^\prime\).
The phenomenon is illustrated and verified by
all known examples of Galois groups \(G=\mathrm{Gal}(\mathrm{F}_3^\infty(K)/K)\)
of \(3\)-class field towers
\(K=\mathrm{F}_3^0(K)<\mathrm{F}_3^1(K)<\mathrm{F}_3^2(K)\le\ldots\le\mathrm{F}_3^\infty(K)\)
of imaginary quadratic number fields \(K=\mathbb{Q}(\sqrt{d})\), \(d<0\),
with elementary \(3\)-class group \(\mathrm{Cl}_3(K)\) of rank two.
Such Galois groups must be Schur \(\sigma\)-groups
and the \textbf{existence of towers with at least four stages} is justified
by showing the non-existence of suitable Schur \(\sigma\)-groups \(G\)
with derived length \(\mathrm{dl}(G)\le 3\).
By means of counter-examples,
it is emphasized that real quadratic number fields
with the same type of \(3\)-class group
reveal a totally different behavior,
usually without extremal path.
\end{abstract}

\maketitle



\section{Introduction}
\label{s:Intro}

\noindent
The statement of our Main Conjecture
\ref{cnj:ExtremalPath} 
in this paper
requires a few preparatory sections on
\(p\)-class field towers,
relation ranks,
transfer kernel types, and
root paths in descendant trees.



\subsection{Galois group \(G\) of the \(p\)-class field tower}
\label{ss:TowerGroup}

\noindent
For an algebraic number field \(K/\mathbb{Q}\),
that is a finite extension of the rational number field \(\mathbb{Q}\),
and a prime number \(p\in\mathbb{P}\),
the \(p\)-\textit{class field tower} \(\mathrm{F}_p^\infty(K)\) of \(K\)
is the maximal unramified pro-\(p\)-extension of \(K\).
The Galois group \(G=\mathrm{Gal}(\mathrm{F}_p^\infty(K)/K)\)
is a potentially infinite pro-\(p\)-group
with finite abelianization \(G/G^\prime\)
isomorphic to the \(p\)-class group \(\mathrm{Cl}_p(K)\) of \(K\).
The derived series of \(G\), which is defined recursively by
\begin{equation}
\label{eqn:DerivedSeries}
G^{(0)}:=G, \text{ and } G^{(j)}:=\lbrack G^{(j-1)},G^{(j-1)}\rbrack, \text{ for all } j\ge 1,
\end{equation}
determines the various \textit{stages} of the tower
\begin{equation}
\label{eqn:Tower}
K=\mathrm{F}_p^0(K)\le\mathrm{F}_p^1(K)\le\mathrm{F}_p^2(K)\le\ldots\le\mathrm{F}_p^\infty(K),
\end{equation}
where \(\mathrm{F}_p^j(K)=\mathrm{Fix}(G^{(j)})\), that is \(G^{(j)}=\mathrm{Gal}(\mathrm{F}_p^\infty(K)/\mathrm{F}_p^j(K))\),
for all \(j\ge 0\), and
\begin{equation}
\label{eqn:Quotient}
G/G^{(j)}=\mathrm{Gal}(\mathrm{F}_p^\infty(K)/K)/\mathrm{Gal}(\mathrm{F}_p^\infty(K)/\mathrm{F}_p^j(K))\simeq\mathrm{Gal}(\mathrm{F}_p^j(K)/K),
\end{equation}
in part.
\(G/G^\prime=G/G^{(1)}\simeq\mathrm{Gal}(\mathrm{F}_p^1(K)/K)\simeq\mathrm{Cl}_p(K)\),
according to Artin's reciprocity law
\cite{Ar1,Ar2}.



\subsection{Bounds for the relation rank of \(G\)}
\label{ss:Shafarevich}

\noindent
The \textit{relation rank}
\(d_2(G):=\dim_{\mathbb{F}_p}(\mathrm{H}^2(G,\mathbb{F}_p))\) of \(G\)
frequently admits a decision concerning
the \textit{length} \(\ell_p(K)\) of the \(p\)-class field tower of \(K\),
based on the \textit{Shafarevich cohomology criterion}
\cite{Sh}
\begin{equation}
\label{eqn:Shafarevich}
d_1(G)\le d_2(G)\le d_1(G)+r+\theta
\end{equation}
in dependence on
the generator rank \(d_1(G):=\dim_{\mathbb{F}_p}(\mathrm{H}^1(G,\mathbb{F}_p))\) of \(G\),
equal to the \(p\)-class rank \(\mathrm{rk}_p(K)\),
the signature \((r_1,r_2)\),
and the torsion free Dirichlet unit rank
\(r=r_1+r_2-1\) of the number field \(K\).
See also
\cite[Thm. 5.1, p. 28]{Ma2015d},
\cite[\S\S\ 1.2--1.3, pp. 75--76]{Ma2016d},
\cite{Ma2016e}.
Here, \(\theta\) denotes the flag
\begin{equation}
\label{eqn:Flag}
\theta=
\begin{cases}
1, & \text{ if } K \text{ contains the } p\text{-th roots of unity,} \\
0, & \text{ otherwise.}
\end{cases}
\end{equation}



\noindent
In particular,
an \textit{imaginary quadratic} number field \(K\)
has signature \((r_1,r_2)=(0,1)\)
and Dirichlet unit rank \(r=0\).
If \(p\) is an \textit{odd} prime,
and \(K\) has a non-cyclic \(p\)-class group
with rank \(\mathrm{rk}_p(K)\ge 2\),
then \(\theta=0\), and we have the following theorem,
according to Koch and Venkov
\cite{KoVe}
or
\cite{Ag}.


{\color{blue}
\begin{theorem}
\label{thm:SchurSigma}
The Galois group \(G\) of the \(p\)-class field tower \(\mathrm{F}_p^\infty(K)\)
of an imaginary quadratic number field \(K\)
with non-cyclic \(p\)-class group \(\mathrm{Cl}_p(K)\)
for an odd prime number \(p\)
must be a so-called \textbf{Schur \(\sigma\)-group},
that is a pro-\(p\)-group \(G\) with balanced presentation \(d_2(G)=d_1(G)\),
which possesses an automorphism \(\sigma\in\mathrm{Aut}(G)\)
acting as inversion \(x^\sigma=x^{-1}\) on the abelianization \(G/G^\prime\).
\end{theorem}
}

\begin{proof}
This is an immediate consequence of the Formulas
\eqref{eqn:Shafarevich}
and 
\eqref{eqn:Flag}
with \(r=\theta=0\),
and the fact that
the non-trivial automorphism \(\tau\in\mathrm{Gal}(K/\mathbb{Q})\)
acts as inversion on \(\mathrm{Cl}_p(K)\simeq G/G^\prime\)
(since \(1+\tau=\mathrm{Norm}_{K/\mathbb{Q}}\) and \(\mathrm{Cl}_p(\mathbb{Q})=1\)),
and has an extension to \(\sigma\in\mathrm{Aut}(G)\)
with \(\sigma\vert_{G/G^\prime}=\tau\).
\end{proof}



\subsection{Transfer kernel type of \(G\)}
\label{ss:Capitulation}

\noindent
In the present paper,
our focus will be on the \textit{smallest odd} prime number \(p=3\)
and \textit{imaginary quadratic} base fields 
\(K=\mathbb{Q}(\sqrt{d})\), \(d<0\),
with elementary \(3\)-class group \(\mathrm{Cl}_3(K)\simeq C_3\times C_3\) of rank two.
The \(3\)-class field tower of such fields
has been investigated for the first time by Arnold Scholz and Olga Taussky in \(1934\)
\cite{SoTa}.
These authors coined the concept of \(3\)-\textit{capitulation type} \(\varkappa(K)\) of \(K\),
defined by the transfer homomorphisms \(T_i:\,\mathrm{Cl}_3(K)\to\mathrm{Cl}_3(E_i)\) of \(3\)-classes
from \(K\) to its four unramified cyclic cubic extensions \(E_i\), \(1\le i\le 4\),
within the Hilbert \(3\)-class field \(\mathrm{F}_3^1(K)\):
\begin{equation}
\label{eqn:Capitulation}
\varkappa(K):=(\varkappa(1),\ldots,\varkappa(4)), \text{ where } \ker(T_i)=N_{\varkappa(i)}.
\end{equation}
Here, \(N_i\) denotes the norm class group \(\mathrm{Norm}_{E_i/K}(\mathrm{Cl}_3(E_i))\), for each \(1\le i\le 4\).
According to the Theorem on Herbrand's quotient
\(\frac{\#\mathrm{H}^1(A,U_i)}{\#\hat{\mathrm{H}}^0(A,U_i)}=\lbrack E_i:K\rbrack=3\)
of the unit group \(U_i:=U(E_i)\)
as a Galois module over the automorphism group \(A:=\mathrm{Gal}(E_i/K)\simeq C_3\),
all the kernels \(\ker(T_i)\)
are cyclic of order \(3\), if \(K\) is imaginary quadratic
(and thus \(\hat{\mathrm{H}}^0(A,U_i)\simeq U(K)/\mathrm{Norm}_{E_i/K}(U_i)=1\)).



Scholz and Taussky arranged the essential \(19\) types of \(3\)-capitulation
(we also speak about the \textit{transfer kernel type}, TKT)
in several sections denoted by upper case letters
\cite[pp. 34--37]{SoTa}.
The designation of individual types by numbers from \(1\) to \(19\)
is taken from
\cite[p. 80]{Ma1991}:


\begin{definition}
\label{dfn:TKT}
Representatives for the \(3\)-capitulation types are arranged in the following sections:
\begin{itemize}
\item
Section A: type A.1 \(\varkappa=(\mathbf{1}111)\),
\item
Section B: types B.2 \(\varkappa=(\mathbf{12}11)\), B.3 \(\varkappa=(\mathbf{1}112)\),
\item
Section C: types C.15 \(\varkappa=(\mathbf{1234})\), C.17 \(\varkappa=(\mathbf{1}342)\), C.18 \(\varkappa=(2341)\),
\item
Section D: types D.5 \(\varkappa=(\mathbf{12}12)\), D.10 \(\varkappa=(\mathbf{1}123)\),
\item
Section E: types E.6 \(\varkappa=(\mathbf{1}122)\), E.8 \(\varkappa=(\mathbf{123}1)\), E.9 \(\varkappa=(\mathbf{1}1\mathbf{3}2)\), E.14 \(\varkappa=(2311)\),
\item
Section F: types F.7 \(\varkappa=(2112)\), F.11 \(\varkappa=(\mathbf{1}321)\), F.12 \(\varkappa=(21\mathbf{3}1)\), F.13 \(\varkappa=(2113)\),
\item
Section G: types G.16 \(\varkappa=(\mathbf{12}43)\), G.19 \(\varkappa=(2143)\),
\item
Section H: type H.4 \(\varkappa=(2111)\).
\end{itemize}
Fixed point capitulation \(\ker(T_i)=N_{i}\), i.e. \(\varkappa(i)=i\), is always marked by using \textbf{boldface} font.
\end{definition}

\noindent
Every other \(3\)-capitulation type is equivalent to some of the representatives in Dfn.
\ref{dfn:TKT}
\cite{Ma2010}.


{\color{blue}
\begin{proposition}
\label{prp:Forbidden}
The types in sections \(\mathrm{B}\) and \(\mathrm{C}\) are generally forbidden
for any number field \(K\).
Type \(\mathrm{A}.1\) enforces the extra special \(3\)-group \(\langle 27,4\rangle\)
of order \(27\) and exponent \(9\)
as Galois group \(G\) and is forbidden for any quadratic field (imaginary and real).
The types in all the other sections give rise to Galois groups \(G\)
which share the extra special \(3\)-group \(\langle 27,3\rangle\)
of order \(27\) and exponent \(3\)
as their common class-\(2\) quotient \(G/\gamma_3(G)\).
(See
\cite{SoTa}.
Here, the groups are denoted by their identifiers in the SmallGroups Library
\cite{BEO2},
and \(\gamma_3(G)\) is the third term in Formula
\eqref{eqn:MembersLCS}.)
\end{proposition}
}



\subsection{Root path of \(G\)}
\label{ss:RootPaths}

\noindent
Let \(p\) be a prime number.
A finite non-abelian \(p\)-group \(G\)
gives rise to a \textit{root path}
in the descendant tree of its abelianization \(G/G^\prime\).
The vertices \(V\) of this directed tree are
isomorphism classes of finite non-abelian \(p\)-groups
sharing a common abelianization.
Two vertices \(V\) and \(W\) are connected by a directed edge \(V\rightarrow W\),
if \(W\) is isomorphic to the quotient of \(V\) by the last non-trivial member \(\gamma_c(V)\)
of its \textit{lower central series} (or \textit{descending} central series),
\begin{equation}
\label{eqn:MembersLCS}
\gamma_1(V)=V>\gamma_2(V)=V^\prime>\gamma_3(V)>\ldots>\gamma_c(V)>\gamma_{c+1}(V)=1,\ c=\mathrm{cl}(V),
\end{equation}
which is defined recursively by
\begin{equation}
\label{eqn:RecursionLCS}
\gamma_1(V):=V, \text{ and } \gamma_j(V):=\lbrack\gamma_{j-1}(V),V\rbrack, \text{ for all } j\ge 2.
\end{equation}
and becomes trivial for \(j>c\) bigger than the \textit{nilpotency class} of \(V\).
In this case, \(W=\pi(V)\) is called \textit{the parent} of \(V\)
and \(V\) is called \textit{an immediate descendant} of \(W\).
The construction of parents can be iterated, and the \textit{root path} of \(V\)
(in the sequel without the last edge)
is given by
\begin{equation}
\label{eqn:RootPath}
V\to\pi(V)\simeq V/\gamma_c(V)\to\pi^2(V)\simeq V/\gamma_{c-1}(V)\to\ldots\to\pi^{c-1}(V)\simeq V/\gamma_2(V)=V/V^\prime.
\end{equation}
The \textit{nuclear rank} \(\nu(W)\) of the parent
determines the possible step sizes \(1\le s\le\nu(W)\)
of the edge \(V\rightarrow W\).
The increment of the logarithmic order, nilpotency class and coclass is given by 
\(\mathrm{lo}(V)-\mathrm{lo}(W)=s\),
\(\mathrm{cl}(V)-\mathrm{cl}(W)=1\) and
\(\mathrm{cc}(V)-\mathrm{cc}(W)=s-1\).
Therefore, in a coclass tree (with constant coclass),
the step size of all edges is restricted to the minimal value \(s=1\).



We are now in the position to state the \textbf{Main Conjecture}
\ref{cnj:ExtremalPath}
of this paper.

{\color{red}
\begin{conjecture}
\label{cnj:ExtremalPath}
\textbf{(Extremal root path of \(G\))}
In the descendant tree \(\mathcal{T}(C_3\times C_3)\)
of the abelian \(3\)-group \(C_3\times C_3\),
the root path \((\pi^{i-1}(G)\to\pi^i(G))_{1\le i\le c-2}\)
of the Galois group \(G=\mathrm{Gal}(\mathrm{F}_3^\infty(K)/K)\)
of the \(3\)-class field tower \(\mathrm{F}_3^\infty(K)\)
of an imaginary quadratic field \(K\)
with elementary \(3\)-class group \(\mathrm{Cl}_3(K)\simeq C_3\times C_3\) of rank two
consists of edges with \textbf{maximal step size (=nuclear rank)}, i.e.
\begin{equation}
\label{eqn:ExtremalPath}
\pi^{i-1}(G)\to\pi^i(G) \text{ is an edge of step size } s_i=\nu(\pi^i(G)), \text{ for each } 1\le i\le c-2,
\end{equation}
where \(3\le c=\mathrm{cl}(G)<\infty\) denotes the nilpotency class of \(G\).
\end{conjecture}
}



\subsection{Layout of this paper}
\label{ss:Overview}

\noindent
The Main Conjecture
\ref{cnj:ExtremalPath}
has never been violated by any known situation involving Schur \(\sigma\)-groups.
It is our desire to underpin the conjecture
with infinite series of parametrized Schur \(\sigma\)-groups,
having proven extremal path property.

At the beginning, in \S\
\ref{ss:Metabelian},
we treat the finitely many metabelian Schur \(\sigma\)-groups.
They have TKTs in Section D.

An infinitude of Schur \(\sigma\)-groups,
arising as quotients of an infinite limit group,
was discovered by Bartholdi and Bush
\cite{BaBu}.
They are investigated in the present context in \S\
\ref{ss:BartholdiBush}.
The members have TKT H.4 and unbounded derived length.

Other infinite sequences of Schur \(\sigma\)-groups,
arising as quotients of infinite limit groups,
were constructed by Newman and ourselves in
\cite{Ma2018b}.
Their root paths are checked for the extremal property in \S\
\ref{ss:ThreeStage}.
All their TKTs are contained in Section E
and they share the common derived length \(3\).
Counts of Schur \(\sigma\)-groups up to order \(3^{14}\)
are given in \S\
\ref{ss:NonMetabelian}.

Root paths with extreme complexity are presented in \S\
\ref{ss:SectionF}.
They occur for Schur \(\sigma\)-groups with TKTs in Section F.
In immediate prosecution,
the striking novelty of the first proven four-stage towers
of \(3\)-class fields over imaginary quadratic fields
is heralded as the \textbf{Main Theorem}
\ref{thm:FourStageTowerF13}
and the beginning of a new era of research on maximal unramified pro-\(p\)-extensions in \S\
\ref{ss:FourStageTower}.

In \S\
\ref{ss:SectionG},
we supplement the exposition with further root paths,
not covered by the preceding developments.
Counter-examples with real quadratic base fields are provided in \S\
\ref{s:CounterExamples}.

A summary is given in \S\
\ref{s:Conclusion}
and personal historical remarks illuminate
the arduous long and winding road
towards our present elevated perspective of 
finite \(p\)-class field towers
in \S\S\
\ref{s:HistorySectionF}
and
\ref{s:PersonalHistory}.



\section{Proving path extremality of \(G\) for infinite series of Schur \(\sigma\)-groups}
\label{s:Proofs}

\noindent
Throughout this paper
we identify isomorphism classes of finite \(p\)-groups with the aid of
their representative in the SmallGroups Database
\cite{BEO2,BEO1},
which is implemented in the computer algebra systems GAP
\cite{GAP}
and Magma
\cite{MAGMA,BCP,BCFS}.
Identifiers (Id) have the format \(\langle\text{order},\text{counter}\rangle\).

\newpage


\subsection{Sporadic metabelian Schur \(\sigma\)-groups \(G\)}
\label{ss:Metabelian}
(\(\mathrm{dl}(G)=2\))

{\color{blue}
\begin{proposition}
\label{prp:SectionD}
Among the finite \(3\)-groups \(V\) with abelian quotient invariants \(V/V^\prime\simeq (3,3)\)
there exist \textbf{precisely two metabelian} Schur \(\sigma\)-groups.
Both are of order \(\mathrm{ord}=3^5=243\), nilpotency class \(\mathrm{cl}=3\) and coclass \(\mathrm{cc}=2\).
The first is \(\langle 243,5\rangle\) with TKT \(\mathrm{D}.10\),
and the second is \(\langle 243,7\rangle\) with TKT \(\mathrm{D}.5\).
They are the unique terminal immediate descendants with step size \(s=2\)
of the extra special group \(\langle 27,3\rangle\),
and thus isolated orphans without parents in the coclass forest \(\mathcal{F}(2)\).
\end{proposition}
}

\begin{proof}
A search for metabelian Schur groups with balanced presentation \(d_2=d_1=2\)
additionally yields the extra special group \(\langle 27,4\rangle\) with TKT A.1,
but this is not a \(\sigma\)-group with generator inverting automorphism.
\end{proof}

\begin{remark}
\label{rmk:SectionD}
According to Nebelung
\cite{Ne1,Ne2},
there do not exist other
finite metabelian \(3\)-groups \(V\) with abelian quotient invariants \(V/V^\prime\simeq (3,3)\)
and TKT in Section A or Section D.
Since an epimorphism onto a Schur group must be an isomorphism,
there cannot exist non-metabelian \(3\)-groups \(V\) with abelian quotient invariants \(V/V^\prime\simeq (3,3)\)
and TKT in one of these two Sections.
Sections A, B, C and D are the only sections with a finite number of associated Schur \(\sigma\)-groups.
\end{remark}


{\color{blue}
\begin{theorem}
\label{thm:SectionD}
The \(3\)-class field tower of an imaginary quadratic field \(K\)
with \(3\)-capitulation type in section \(\mathrm{D}\)
is metabelian with length \(\ell_3(K)=2\).
The Galois group \(G=\mathrm{Gal}(\mathrm{F}_3^\infty(K)/K)\) satisfies the extremal path property with
\(s_1=2\).
\end{theorem}
}

\begin{proof}
The length was proven by Scholz and Taussky in
\cite{SoTa}.
The extremal property is a consequence of Table
\ref{tbl:D10}
and a similar table for type D.5 with \(\langle 243,5\rangle\) replaced by \(\langle 243,7\rangle\).
\end{proof}


\renewcommand{\arraystretch}{1.1}

\begin{table}[ht]
\caption{Root path of \(G\) for transfer kernel type D.10}
\label{tbl:D10}
\begin{center}
\begin{tabular}{|c|c||c|c|c|}
\hline
 Ancestor   & Id                       & \((\nu,\mu)\)                   & \((N_s/C_s)_{1\le s\le\nu}\)        & TKT  \\
\hline
 \(\pi(G)\) & \(\langle 27,3\rangle\)  & \(({\color{red}\mathbf{2,4}})\) & \((4/1,{\color{red}\mathbf{7/5}})\) & a.1  \\
 \(G\)      & \(\langle 243,5\rangle\) & \(({\color{red}\mathbf{0,2}})\) &                                     & D.10 \\
\hline
\end{tabular}
\end{center}
\end{table}



\noindent
Table
\ref{tbl:CQ3x3}
shows the statistical distribution of the simplest \(3\)-capitulation types
in three ranges of discriminants \(d\).
The range \(-10^6<d<0\) was investigated by ourselves in
\cite{Ma2012,Ma2014},
in
\cite{Ma2016b}
we extended to \(-10^7<d<0\),
and the range \(-10^8<d<0\) is due to Boston, Bush and Hajir
\cite{BBH}.
Without doubt,
the metabelian Schur \(\sigma\)-group \(\langle 243,5\rangle\)
with TKT D.10, resp. all four cases (see Fig.
\ref{fig:SporCc2}),
enjoy the most dense population,
covering nearly one third, resp. two thirds, of all cases.

\noindent
Here and in the sequel,
we additionally need the second component of the \textit{Artin pattern}
\(\mathrm{AP}(G)=(\varkappa(G),\tau(G))\).
The \textit{transfer target type} (TTT) \(\tau(G)\)
consists of the logarithmic abelian type invariants of the four \(3\)-class groups
\(\mathrm{Cl}_3(E_i)\), \(1\le i\le 4\),
corresponding to the TKT \(\varkappa(G)\),
and \(\varepsilon\) denotes the number of components of \(\tau(G)\)
isomorphic to \(1^3\simeq C_3\times C_3\times C_3\).
(Similarly, \(21\simeq C_9\times C_3\).)


\renewcommand{\arraystretch}{1.1}

\begin{table}[ht]
\caption{Smallest \(G/G^{(2)}\) for imaginary quadratic fields of type \((3,3)\)}
\label{tbl:CQ3x3}
\begin{center}
\begin{tabular}{|c|r||l|l|l|l||l|}
\hline
 Discriminant     & Total\#      & \(\varepsilon=1\) & \(\varepsilon=2\)    & \(\varepsilon=3\) & \(\varepsilon=0\) & \(\Sigma\%\)  \\
\hline
 \(-10^6<d<0\)    &   \(2\,020\) & \(33.0\%\)        & \(13.3\%\)           & \(14.7\%\)        & \(4.7\%\)         & \(65.7\%\)    \\
 \(-10^7<d<0\)    &  \(24\,476\) & \(31.14\%\)       & \(14.81\%\)          & \(14.79\%\)       & \(4.163\%\)       & \(64.903\%\)  \\
 \(-10^8<d<0\)    & \(276\,375\) & \(30.159\%\)      & \(14.979\%\)         & \(14.823\%\)      & \(3.7724\%\)      & \(63.7334\%\) \\
\hline
 \(\tau(G)\)      &              & \((1^3,(21)^3)\)  & \(((1^3)^2,(21)^2)\) & \(((1^3)^3,21)\)  & \(((21)^4)\)      &               \\
 \(\varkappa(G)\) &              & \((2241)\)        & \((4224)\)           & \((4443)\)        & \((2143)\)        &               \\
 TKT              &              & D.10              & D.5                  & H.4               & G.19              &               \\
 \(G/G^{(2)}\simeq\) & & \(\langle 243,5\rangle\) & \(\langle 243,7\rangle\) & \(\langle 729,45\rangle\) & \(\langle 729,57\rangle\) & \\
\hline
\end{tabular}
\end{center}
\end{table}

\newpage


{\tiny

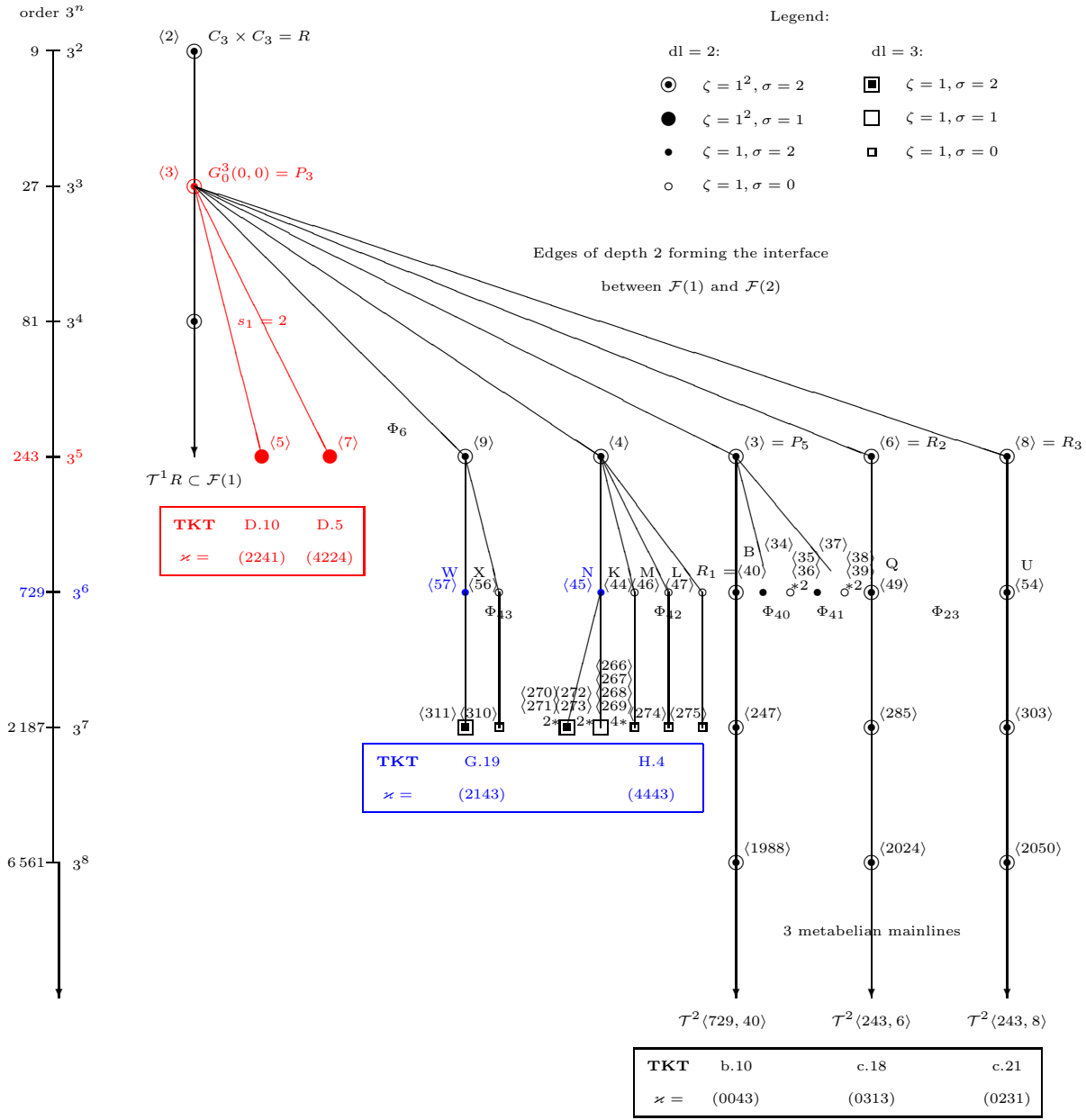
\begin{figure}[hb]
\caption{Extremal paths to Schur \(\sigma\)-groups \(\langle 5\rangle\), \(\langle 7\rangle\), log ord \(5\), in coclass forest \(\mathcal{F}(2)\)}
\label{fig:SporCc2}


{\tiny

\begin{center}

\setlength{\unitlength}{1.0cm}
\begin{picture}(16,17)(0,-14)

\put(0,2.5){\makebox(0,0)[cb]{order \(3^n\)}}

\put(0,2){\line(0,-1){12}}
\multiput(-0.1,2)(0,-2){7}{\line(1,0){0.2}}

\put(-0.2,2){\makebox(0,0)[rc]{\(9\)}}
\put(0.2,2){\makebox(0,0)[lc]{\(3^2\)}}
\put(-0.2,0){\makebox(0,0)[rc]{\(27\)}}
\put(0.2,0){\makebox(0,0)[lc]{\(3^3\)}}
\put(-0.2,-2){\makebox(0,0)[rc]{\(81\)}}
\put(0.2,-2){\makebox(0,0)[lc]{\(3^4\)}}
{\color{red}
\put(-0.2,-4){\makebox(0,0)[rc]{\(243\)}}
\put(0.2,-4){\makebox(0,0)[lc]{\(3^5\)}}
}
{\color{blue}
\put(-0.2,-6){\makebox(0,0)[rc]{\(729\)}}
\put(0.2,-6){\makebox(0,0)[lc]{\(3^6\)}}
}
\put(-0.2,-8){\makebox(0,0)[rc]{\(2\,187\)}}
\put(0.2,-8){\makebox(0,0)[lc]{\(3^7\)}}
\put(-0.2,-10){\makebox(0,0)[rc]{\(6\,561\)}}
\put(0.2,-10){\makebox(0,0)[lc]{\(3^8\)}}

\put(0,-10){\vector(0,-1){2}}

\put(10.5,2.5){\makebox(0,0)[lc]{Legend:}}

\put(9,2){\makebox(0,0)[lc]{\(\mathrm{dl}=2\):}}
\put(9,1.5){\circle{0.2}}
\put(9,1.5){\circle*{0.1}}
\put(9,1){\circle*{0.2}}
\put(9,0.5){\circle*{0.1}}
\put(9,0){\circle{0.1}}
\put(9.5,1.5){\makebox(0,0)[lc]{\(\zeta=1^2,\sigma=2\)}}
\put(9.5,1){\makebox(0,0)[lc]{\(\zeta=1^2,\sigma=1\)}}
\put(9.5,0.5){\makebox(0,0)[lc]{\(\zeta=1,\sigma=2\)}}
\put(9.5,0){\makebox(0,0)[lc]{\(\zeta=1,\sigma=0\)}}

\put(12,2){\makebox(0,0)[lc]{\(\mathrm{dl}=3\):}}
\put(11.9,1.4){\framebox(0.2,0.2){\vrule height 3pt width 3pt}}
\put(11.9,0.9){\framebox(0.2,0.2){}}
\put(11.95,0.45){\framebox(0.1,0.1){}}
\put(12.5,1.5){\makebox(0,0)[lc]{\(\zeta=1,\sigma=2\)}}
\put(12.5,1){\makebox(0,0)[lc]{\(\zeta=1,\sigma=1\)}}
\put(12.5,0.5){\makebox(0,0)[lc]{\(\zeta=1,\sigma=0\)}}

\put(2.2,2.2){\makebox(0,0)[lc]{\(C_3\times C_3=R\)}}
\put(1.8,2.2){\makebox(0,0)[rc]{\(\langle 2\rangle\)}}
\put(2,2){\circle{0.2}}
\put(2,2){\circle*{0.1}}
\put(2,2){\line(0,-1){2}}
{\color{red}
\put(2,0){\circle{0.2}}
\put(2,0){\circle*{0.1}}
\put(2.2,0.2){\makebox(0,0)[lc]{\(G^3_0(0,0)=P_3\)}}
\put(1.8,0.2){\makebox(0,0)[rc]{\(\langle 3\rangle\)}}
}

\multiput(2,-2)(0,-2){1}{\circle{0.2}}
\multiput(2,-2)(0,-2){1}{\circle*{0.1}}

\put(2,0){\vector(0,-1){4}}
\put(2,-4.2){\makebox(0,0)[ct]{\(\mathcal{T}^1R\subset\mathcal{F}(1)\)}}

{\color{red}
\put(2,0){\line(1,-4){1}}
\put(2,0){\line(1,-2){2}}
\put(3,-2){\makebox(0,0)[cc]{\(s_1=2\)}}
}
\put(2,0){\line(1,-1){4}}
\put(2,0){\line(3,-2){6}}
\put(2,0){\line(2,-1){8}}
\put(2,0){\line(5,-2){10}}
\put(2,0){\line(3,-1){12}}
\put(7,-1){\makebox(0,0)[lc]{Edges of depth \(2\) forming the interface}}
\put(8,-1.5){\makebox(0,0)[lc]{between \(\mathcal{F}(1)\) and \(\mathcal{F}(2)\)}}

\put(5,-3.6){\makebox(0,0)[cc]{\(\Phi_6\)}}
{\color{red}
\multiput(3,-4)(1,0){2}{\circle*{0.2}}
\put(3.1,-3.9){\makebox(0,0)[lb]{\(\langle 5\rangle\)}}
\put(4.1,-3.9){\makebox(0,0)[lb]{\(\langle 7\rangle\)}}
}
\multiput(6,-4)(2,0){5}{\circle{0.2}}
\multiput(6,-4)(2,0){5}{\circle*{0.1}}
\put(6.1,-3.9){\makebox(0,0)[lb]{\(\langle 9\rangle\)}}
\put(8.1,-3.9){\makebox(0,0)[lb]{\(\langle 4\rangle\)}}
\put(10.1,-3.9){\makebox(0,0)[lb]{\(\langle 3\rangle=P_5\)}}
\put(12.1,-3.9){\makebox(0,0)[lb]{\(\langle 6\rangle=R_2\)}}
\put(14.1,-3.9){\makebox(0,0)[lb]{\(\langle 8\rangle=R_3\)}}

{\color{red}
\put(2,-5){\makebox(0,0)[cc]{\textbf{TKT}}}
\put(3,-5){\makebox(0,0)[cc]{D.10}}
\put(4,-5){\makebox(0,0)[cc]{D.5}}
\put(2,-5.5){\makebox(0,0)[cc]{\(\varkappa=\)}}
\put(3,-5.5){\makebox(0,0)[cc]{\((2241)\)}}
\put(4,-5.5){\makebox(0,0)[cc]{\((4224)\)}}
\put(1.5,-5.75){\framebox(3,1){}}
}

\put(6.5,-6.3){\makebox(0,0)[cc]{\(\Phi_{43}\)}}
\put(6,-4){\line(0,-1){2}}
{\color{blue}
\put(5.9,-5.9){\makebox(0,0)[rc]{\(\langle 57\rangle\)}}
\put(5.9,-5.7){\makebox(0,0)[rc]{W}}
}
\put(6.5,-5.9){\makebox(0,0)[rc]{\(\langle 56\rangle\)}}
\put(6.3,-5.7){\makebox(0,0)[rc]{X}}
\put(5.9,-7.8){\makebox(0,0)[rc]{\(\langle 311\rangle\)}}
\put(6.5,-7.8){\makebox(0,0)[rc]{\(\langle 310\rangle\)}}
\put(6,-4){\line(1,-4){0.5}}
{\color{blue}
\put(6,-6){\circle*{0.1}}
}
\put(6.5,-6){\circle{0.1}}
\multiput(6,-6)(0.5,0){2}{\line(0,-1){2}}
\put(5.9,-8.1){\framebox(0.2,0.2){\vrule height 3pt width 3pt}}
\put(6.45,-8.05){\framebox(0.1,0.1){}}

\put(9,-6.3){\makebox(0,0)[cc]{\(\Phi_{42}\)}}
\put(8,-4){\line(0,-1){2}}
{\color{blue}
\put(7.9,-5.9){\makebox(0,0)[rc]{\(\langle 45\rangle\)}}
\put(7.9,-5.7){\makebox(0,0)[rc]{N}}
}
\put(8.5,-5.9){\makebox(0,0)[rc]{\(\langle 44\rangle\)}}
\put(8.3,-5.7){\makebox(0,0)[rc]{K}}
\put(8.9,-5.9){\makebox(0,0)[rc]{\(\langle 46\rangle\)}}
\put(8.8,-5.7){\makebox(0,0)[rc]{M}}
\put(9.4,-5.9){\makebox(0,0)[rc]{\(\langle 47\rangle\)}}
\put(9.2,-5.7){\makebox(0,0)[rc]{L}}
\put(7.4,-7.5){\makebox(0,0)[rc]{\(\langle 270\rangle\)}}
\put(7.4,-7.7){\makebox(0,0)[rc]{\(\langle 271\rangle\)}}
\put(7.9,-7.5){\makebox(0,0)[rc]{\(\langle 272\rangle\)}}
\put(7.9,-7.7){\makebox(0,0)[rc]{\(\langle 273\rangle\)}}
\put(8.5,-7.1){\makebox(0,0)[rc]{\(\langle 266\rangle\)}}
\put(8.5,-7.3){\makebox(0,0)[rc]{\(\langle 267\rangle\)}}
\put(8.5,-7.5){\makebox(0,0)[rc]{\(\langle 268\rangle\)}}
\put(8.5,-7.7){\makebox(0,0)[rc]{\(\langle 269\rangle\)}}
\put(9.0,-7.8){\makebox(0,0)[rc]{\(\langle 274\rangle\)}}
\put(9.6,-7.8){\makebox(0,0)[rc]{\(\langle 275\rangle\)}}
\put(8,-4){\line(1,-4){0.5}}
\put(8,-4){\line(1,-2){1}}
\put(8,-4){\line(3,-4){1.5}}
{\color{blue}
\put(8,-6){\circle*{0.1}}
}
\multiput(8.5,-6)(0.5,0){3}{\circle{0.1}}
\put(8,-6){\line(-1,-4){0.5}}
\multiput(8,-6)(0.5,0){4}{\line(0,-1){2}}
\put(7.4,-8.1){\framebox(0.2,0.2){\vrule height 3pt width 3pt}}
\put(7.9,-8.1){\framebox(0.2,0.2){}}
\multiput(8.45,-8.05)(0.5,0){3}{\framebox(0.1,0.1){}}
\multiput(7.4,-7.9)(0.5,0){2}{\makebox(0,0)[rc]{\(2*\)}}
\put(8.4,-7.9){\makebox(0,0)[rc]{\(4*\)}}

{\color{blue}
\put(5,-8.5){\makebox(0,0)[cc]{\textbf{TKT}}}
\put(6.25,-8.5){\makebox(0,0)[cc]{G.19}}
\put(8.75,-8.5){\makebox(0,0)[cc]{H.4}}
\put(5,-9){\makebox(0,0)[cc]{\(\varkappa=\)}}
\put(6.25,-9){\makebox(0,0)[cc]{\((2143)\)}}
\put(8.75,-9){\makebox(0,0)[cc]{\((4443)\)}}
\put(4.5,-9.25){\framebox(5,1){}}
}


\put(10.6,-6.3){\makebox(0,0)[cc]{\(\Phi_{40}\)}}
\put(11.4,-6.3){\makebox(0,0)[cc]{\(\Phi_{41}\)}}
\multiput(10,-4)(0,-2){3}{\line(0,-1){2}}
\multiput(10,-6)(0,-2){3}{\circle{0.2}}
\multiput(10,-6)(0,-2){3}{\circle*{0.1}}
\put(10.0,-5.7){\makebox(0,0)[rc]{\(R_1=\)}}
\put(10.0,-5.7){\makebox(0,0)[lc]{\(\langle 40\rangle\)}}
\put(10.1,-5.4){\makebox(0,0)[lc]{B}}
\put(10.4,-5.3){\makebox(0,0)[lc]{\(\langle 34\rangle\)}}
\put(10.8,-5.5){\makebox(0,0)[lc]{\(\langle 35\rangle\)}}
\put(10.8,-5.7){\makebox(0,0)[lc]{\(\langle 36\rangle\)}}
\put(11.2,-5.3){\makebox(0,0)[lc]{\(\langle 37\rangle\)}}
\put(11.6,-5.5){\makebox(0,0)[lc]{\(\langle 38\rangle\)}}
\put(11.6,-5.7){\makebox(0,0)[lc]{\(\langle 39\rangle\)}}
\put(10.1,-7.8){\makebox(0,0)[lc]{\(\langle 247\rangle\)}}
\put(10.1,-9.8){\makebox(0,0)[lc]{\(\langle 1988\rangle\)}}
\put(10,-4){\line(1,-4){0.4}}
\put(10,-4){\line(5,-6){1.4}}
\multiput(10.4,-6)(0.8,0){2}{\circle*{0.1}}
\multiput(10.8,-6)(0.8,0){2}{\circle{0.1}}
\multiput(10.8,-5.9)(0.8,0){2}{\makebox(0,0)[lc]{\(*2\)}}
\put(10,-10){\vector(0,-1){2}}
\put(9.8,-12.2){\makebox(0,0)[ct]{\(\mathcal{T}^2\langle 729,40\rangle\)}}


\multiput(12,-4)(0,-2){3}{\line(0,-1){2}}
\multiput(12,-6)(0,-2){3}{\circle{0.2}}
\multiput(12,-6)(0,-2){3}{\circle*{0.1}}
\put(12.1,-5.9){\makebox(0,0)[lc]{\(\langle 49\rangle\)}}
\put(12.2,-5.6){\makebox(0,0)[lc]{Q}}
\put(12.1,-7.8){\makebox(0,0)[lc]{\(\langle 285\rangle\)}}
\put(12.1,-9.8){\makebox(0,0)[lc]{\(\langle 2024\rangle\)}}
\put(12,-10){\vector(0,-1){2}}
\put(12,-12.2){\makebox(0,0)[ct]{\(\mathcal{T}^2\langle 243,6\rangle\)}}

\put(13.1,-6.3){\makebox(0,0)[cc]{\(\Phi_{23}\)}}
\put(12,-11){\makebox(0,0)[cc]{\(3\) metabelian mainlines}}

\multiput(14,-4)(0,-2){3}{\line(0,-1){2}}
\multiput(14,-6)(0,-2){3}{\circle{0.2}}
\multiput(14,-6)(0,-2){3}{\circle*{0.1}}
\put(14.1,-5.9){\makebox(0,0)[lc]{\(\langle 54\rangle\)}}
\put(14.2,-5.6){\makebox(0,0)[lc]{U}}
\put(14.1,-7.8){\makebox(0,0)[lc]{\(\langle 303\rangle\)}}
\put(14.1,-9.8){\makebox(0,0)[lc]{\(\langle 2050\rangle\)}}
\put(14,-10){\vector(0,-1){2}}
\put(14,-12.2){\makebox(0,0)[ct]{\(\mathcal{T}^2\langle 243,8\rangle\)}}

\put(9,-13){\makebox(0,0)[cc]{\textbf{TKT}}}
\put(10,-13){\makebox(0,0)[cc]{b.10}}
\put(12,-13){\makebox(0,0)[cc]{c.18}}
\put(14,-13){\makebox(0,0)[cc]{c.21}}
\put(9,-13.5){\makebox(0,0)[cc]{\(\varkappa=\)}}
\put(10,-13.5){\makebox(0,0)[cc]{\((0043)\)}}
\put(12,-13.5){\makebox(0,0)[cc]{\((0313)\)}}
\put(14,-13.5){\makebox(0,0)[cc]{\((0231)\)}}
\put(8.5,-13.75){\framebox(6,1){}}

\end{picture}

\end{center}

}

\end{figure}

}

\newpage


\subsection{Infinite series of Schur \(\sigma\)-groups \(G\) with unbounded derived length}
\label{ss:BartholdiBush}
(\(\mathrm{dl}(G)\ge 3\))

{\color{blue}
\begin{proposition}
\label{prp:SporadicSectionH}
Among the finite \(3\)-groups \(V\) with abelian quotient invariants \(V/V^\prime\simeq (3,3)\),
transfer kernel type \(\mathrm{H}.4\)
and transfer target type \(\tau(V)=((1^3)^3,21)\),
also uniquely characterized by the invariant \(\varepsilon=3\),
there exists \textbf{an infinite sequence} of Schur \(\sigma\)-groups \((S_n)_{n\ge 0}\)
with unbounded derived length \(\mathrm{dl}(S_n)\ge 3\).
(The law for the dependence of \(\mathrm{dl}(S_n)\) on \(n\) is given in
\cite{BaBu}.)
Their logarithmic order, nilpotency class and coclass are given by the following laws:
\begin{equation}
\label{eqn:SporadicSectionH}
\mathrm{lo}(S_n)=8+3n, \quad \mathrm{cl}(S_n)=5+2n, \quad \mathrm{cc}(S_n)=3+n, \text{ for all } n\ge 0.
\end{equation}
\end{proposition}
}

\begin{proof}
The existence and the unbounded derived length of the infinite sequence was proved by Bartholdi and Bush
\cite{BaBu}.
The deterministic laws for invariants were deduced by ourselves
\cite{Ma2015b}
from the structure of the purged descendant tree of \(R:=\langle 243,4\rangle\),
restricted to \(\sigma\)-groups with generator inverting automorphism.
The tree is not a coclass tree, but it has an infinite main trunk
with strictly alternating step sizes \(s=1\) and \(s=2\),
and thus contains periodic bifurcations to higher coclass.
Since \(S_n=R(-\#1;a_j-\#2;b_j)_{1\le j\le n+1}\) with certain \(1\le a_j\le 4\), \(1\le b_j\le 2\),
for each \(n\ge 0\),
it follows that
\(\mathrm{lo}(S_n)=\mathrm{lo}(R)+3(n+1)=5+3n+3=8+3n\),
\(\mathrm{cl}(S_n)=\mathrm{cl}(R)+2(n+1)=3+2n+2=5+2n\), and
\(\mathrm{cc}(S_n)=\mathrm{cc}(R)+(n+1)=2+n+1=3+n\).
\end{proof}


{\color{blue}
\begin{theorem}
\label{thm:SporadicSectionH}
The \(3\)-class field tower of an imaginary quadratic field
with TKT \(\mathrm{H}.4\) and TTT \(\tau=((1^3)^3,21)\)
is non-metabelian with unbounded length \(\ell_3(K)\ge 3\).
In the simplest case, the Galois group \(G\) satisfies the extremal property with
\(s_1=2\), \(s_2=1\) and \(s_3=2\).
Generally, the extremal path property of all cases
is satisfied with strictly alternating step sizes \(s=2\) and \(s=1\).
\end{theorem}
}

\begin{proof}
The unbounded length was proved by Bartholdi and Bush in
\cite{BaBu}.
The extremal property for all cases was proved by ourselves in
\cite{Ma2015b,Ma2017}.
For the simplest case, see Table
\ref{tbl:H4}.
\end{proof}


\renewcommand{\arraystretch}{1.1}

\begin{table}[ht]
\caption{Root path of \(G=S_0\) for the simplest case of transfer kernel type H.4}
\label{tbl:H4}
\begin{center}
\begin{tabular}{|c|c||c|c|c|}
\hline
 Ancestor     & Id                          & \((\nu,\mu)\)                   & \((N_s/C_s)_{1\le s\le\nu}\)        & TKT \\
\hline
 \(\pi^3(G)\) & \(\langle 27,3\rangle\)     & \(({\color{red}\mathbf{2,4}})\) & \((4/1,{\color{red}\mathbf{7/5}})\) & a.1 \\
 \(\pi^2(G)\) & \(\langle 243,4\rangle\)    & \(({\color{red}\mathbf{1,3}})\) & \(({\color{red}\mathbf{4/4}})\)     & H.4 \\
 \(\pi(G)\)   & \(\langle 729,45\rangle\)   & \(({\color{red}\mathbf{2,4}})\) & \((4/0,{\color{red}\mathbf{2/1}})\) & H.4 \\
 \(G\)        & \(\langle 6561,606\rangle\) & \(({\color{red}\mathbf{0,2}})\) &                                     & H.4 \\
\hline
\end{tabular}
\end{center}
\end{table}


\begin{example}
\label{exm:H4}
The root path of the smallest Schur \(\sigma\)-group
with TKT H.4 and TTT \(\tau=((1^3)^3,21)\)
is illustrated with red color in Figure
\ref{fig:H4S0}.
The group is denoted by \(S_0=\langle 6561,606\rangle\).
In
\cite{Ma2015b,Ma2017}
it was proved to be the Galois group \(G=\mathrm{Gal}(\mathrm{F}_3^\infty(K)/K)\)
of the three-stage \(3\)-class field tower
of the quadratic fields \(\mathbb{Q}(\sqrt{d})\) with discriminants
\(d\in\lbrace -3896, -25447, -27355\rbrace\).

The root path of the next Schur \(\sigma\)-group
\(S_1=\langle 6561,605\rangle-\#1;2-\#2;2\)
with order \(3^{11}=177147\)
is drawn with red color in Figure
\ref{fig:H4S1}.
In
\cite{Ma2015b,Ma2017}
it was proved that
the Galois group \(G=\mathrm{Gal}(\mathrm{F}_3^\infty(K)/K)\)
for the quadratic fields \(\mathbb{Q}(\sqrt{d})\) with discriminants
\(d\in\lbrace -6583, -23428, -27991\rbrace\)
is certainly not isomorphic to \(S_0\).
It is currently beyond the reach of actual computations to decide
whether \(G\simeq S_i\) with \(1\le i\le 2\) and \(\mathrm{dl}(G)=3\)
or \(G\simeq S_i\) with \(i\ge 3\) and \(\mathrm{dl}(G)\ge 4\).

The root path of the Schur \(\sigma\)-groups
\(S_2\), resp. \(S_3\),
is shown in Figure
\ref{fig:H4S2},
resp.
\ref{fig:H4S3}.

\noindent
In all figures, the tree is \textit{purged}
in the sense that it is restricted to \(\sigma\)-groups \(V\)
with generator inverting automorphism \(\sigma\in\mathrm{Aut}(V)\).
\end{example}

\newpage


{\tiny

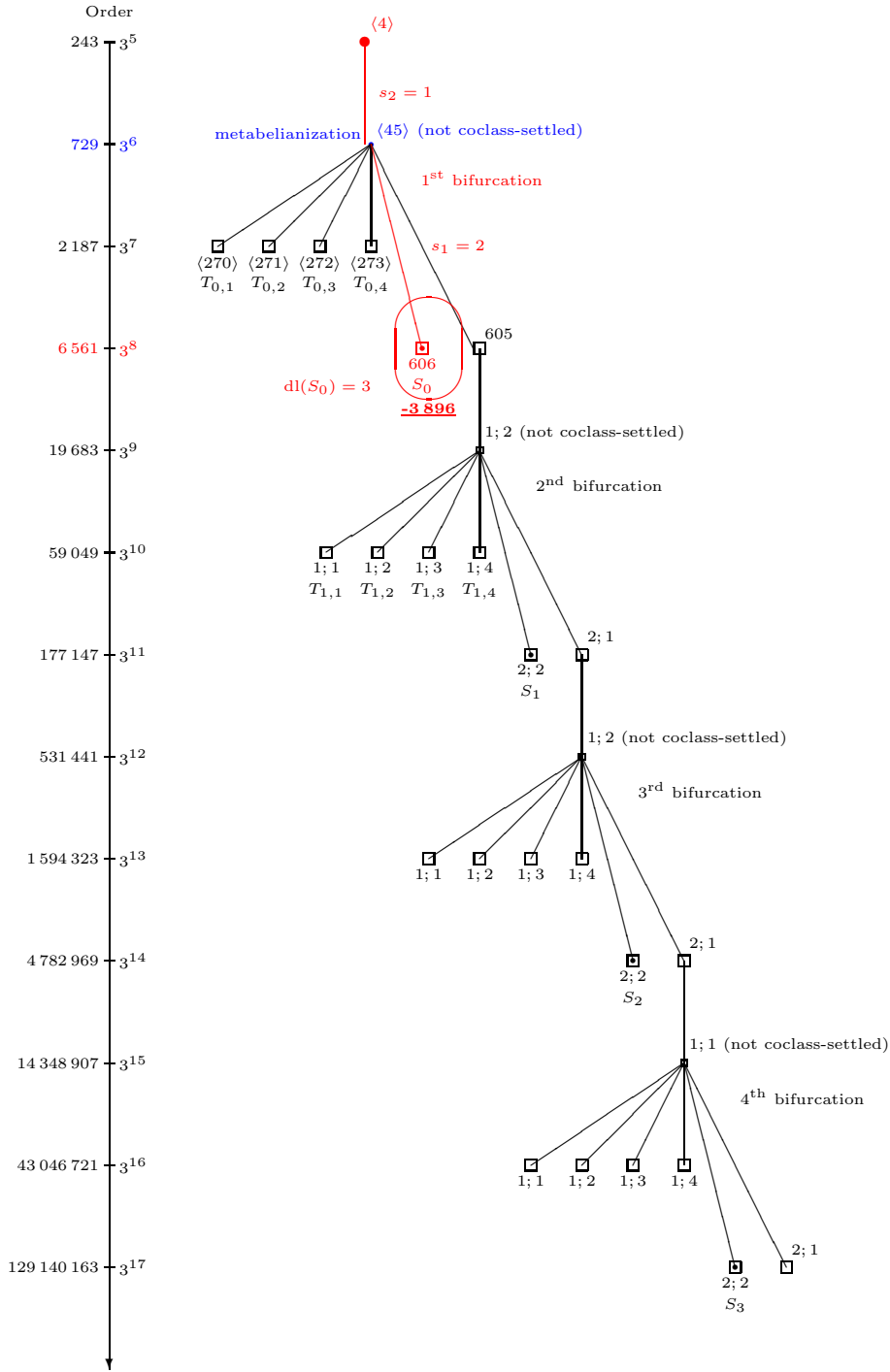
\begin{figure}[hb]
\caption{Extremal path to Schur \(\sigma\)-group \(S_0\), log ord \(8\), on purged tree \(\mathcal{T}_\ast(\langle 243,4\rangle)\)}
\label{fig:H4S0}


\setlength{\unitlength}{0.7cm}
\begin{picture}(18,26.5)(-6,-25.5)

\put(-5,0.5){\makebox(0,0)[cb]{Order}}
\put(-5,0){\line(0,-1){24}}
\multiput(-5.1,0)(0,-2){13}{\line(1,0){0.2}}
\put(-5.2,0){\makebox(0,0)[rc]{\(243\)}}
\put(-4.8,0){\makebox(0,0)[lc]{\(3^5\)}}
{\color{blue}
\put(-5.2,-2){\makebox(0,0)[rc]{\(729\)}}
\put(-4.8,-2){\makebox(0,0)[lc]{\(3^6\)}}
}
\put(-5.2,-4){\makebox(0,0)[rc]{\(2\,187\)}}
\put(-4.8,-4){\makebox(0,0)[lc]{\(3^7\)}}
{\color{red}
\put(-5.2,-6){\makebox(0,0)[rc]{\(6\,561\)}}
\put(-4.8,-6){\makebox(0,0)[lc]{\(3^8\)}}
}
\put(-5.2,-8){\makebox(0,0)[rc]{\(19\,683\)}}
\put(-4.8,-8){\makebox(0,0)[lc]{\(3^9\)}}
\put(-5.2,-10){\makebox(0,0)[rc]{\(59\,049\)}}
\put(-4.8,-10){\makebox(0,0)[lc]{\(3^{10}\)}}
\put(-5.2,-12){\makebox(0,0)[rc]{\(177\,147\)}}
\put(-4.8,-12){\makebox(0,0)[lc]{\(3^{11}\)}}
\put(-5.2,-14){\makebox(0,0)[rc]{\(531\,441\)}}
\put(-4.8,-14){\makebox(0,0)[lc]{\(3^{12}\)}}
\put(-5.2,-16){\makebox(0,0)[rc]{\(1\,594\,323\)}}
\put(-4.8,-16){\makebox(0,0)[lc]{\(3^{13}\)}}
\put(-5.2,-18){\makebox(0,0)[rc]{\(4\,782\,969\)}}
\put(-4.8,-18){\makebox(0,0)[lc]{\(3^{14}\)}}
\put(-5.2,-20){\makebox(0,0)[rc]{\(14\,348\,907\)}}
\put(-4.8,-20){\makebox(0,0)[lc]{\(3^{15}\)}}
\put(-5.2,-22){\makebox(0,0)[rc]{\(43\,046\,721\)}}
\put(-4.8,-22){\makebox(0,0)[lc]{\(3^{16}\)}}
\put(-5.2,-24){\makebox(0,0)[rc]{\(129\,140\,163\)}}
\put(-4.8,-24){\makebox(0,0)[lc]{\(3^{17}\)}}
\put(-5,-24){\vector(0,-1){2}}

{\color{red}
\put(0,0){\circle*{0.2}}
\put(0,0){\line(0,-1){2}}
\put(0.8,-1){\makebox(0,0)[cc]{\(s_2=1\)}}
\put(0.1,0.2){\makebox(0,0)[lb]{\(\langle 4\rangle\)}}
\put(1.1,-2.8){\makebox(0,0)[lb]{\(1^{\text{st}}\) bifurcation}}
}
{\color{blue}
\put(0,-2){\circle*{0.1}}
\put(0.1,-1.9){\makebox(0,0)[lb]{\(\langle 45\rangle\) (not coclass-settled)}}
\put(-0.2,-1.9){\makebox(0,0)[rb]{metabelianization}}
}
\multiput(0,-2)(0,-4){1}{\line(-3,-2){3}}
\multiput(0,-2)(0,-4){1}{\line(-1,-1){2}}
\multiput(0,-2)(0,-4){1}{\line(-1,-2){1}}
\multiput(0,-2)(0,-2){1}{\line(0,-1){2}}
\multiput(-3.1,-4.1)(1,0){4}{\framebox(0.2,0.2){}}
\put(-3,-4.2){\makebox(0,0)[ct]{\(\langle 270\rangle\)}}
\put(-2,-4.2){\makebox(0,0)[ct]{\(\langle 271\rangle\)}}
\put(-1,-4.2){\makebox(0,0)[ct]{\(\langle 272\rangle\)}}
\put(0,-4.2){\makebox(0,0)[ct]{\(\langle 273\rangle\)}}
\put(-3,-4.6){\makebox(0,0)[ct]{\(T_{0,1}\)}}
\put(-2,-4.6){\makebox(0,0)[ct]{\(T_{0,2}\)}}
\put(-1,-4.6){\makebox(0,0)[ct]{\(T_{0,3}\)}}
\put(0,-4.6){\makebox(0,0)[ct]{\(T_{0,4}\)}}
\put(0,-2){\line(1,-2){2}}
{\color{red}
\put(0,-2){\line(1,-4){1}}
\put(1.7,-4){\makebox(0,0)[cc]{\(s_1=2\)}}
\multiput(0.9,-6.1)(1,0){1}{\framebox(0.2,0.2){}}
\multiput(1,-6)(2,-6){1}{\circle*{0.1}}
\put(1,-6.2){\makebox(0,0)[ct]{\(606\)}}
\put(1,-6.6){\makebox(0,0)[ct]{\(S_0\)}}
\put(0,-6.6){\makebox(0,0)[rt]{\(\mathrm{dl}(S_0)=3\)}}
}




{\color{red}
\put(1,-6){\oval(1.3,2.0)}
\put(1,-7.2){\makebox(0,0)[cc]{\underbar{\textbf{-3\,896}}}}
}



\multiput(2,-6)(0,-2){2}{\line(0,-1){2}}
\put(1.9,-6.1){\framebox(0.2,0.2){}}
\put(1.95,-8.05){\framebox(0.1,0.1){}}
\put(2.1,-5.8){\makebox(0,0)[lb]{\(605\)}}
\put(2.1,-7.8){\makebox(0,0)[lb]{\(1;2\) (not coclass-settled)}}
\put(3.1,-8.8){\makebox(0,0)[lb]{\(2^{\text{nd}}\) bifurcation}}
\multiput(2,-8)(0,-4){1}{\line(-3,-2){3}}
\multiput(2,-8)(0,-4){1}{\line(-1,-1){2}}
\multiput(2,-8)(0,-4){1}{\line(-1,-2){1}}
\multiput(-1.1,-10.1)(1,0){4}{\framebox(0.2,0.2){}}
\put(-1,-10.2){\makebox(0,0)[ct]{\(1;1\)}}
\put(0,-10.2){\makebox(0,0)[ct]{\(1;2\)}}
\put(1,-10.2){\makebox(0,0)[ct]{\(1;3\)}}
\put(2,-10.2){\makebox(0,0)[ct]{\(1;4\)}}
\put(-1,-10.6){\makebox(0,0)[ct]{\(T_{1,1}\)}}
\put(0,-10.6){\makebox(0,0)[ct]{\(T_{1,2}\)}}
\put(1,-10.6){\makebox(0,0)[ct]{\(T_{1,3}\)}}
\put(2,-10.6){\makebox(0,0)[ct]{\(T_{1,4}\)}}
\put(2,-8){\line(1,-2){2}}
\put(2,-8){\line(1,-4){1}}
\multiput(2.9,-12.1)(1,0){1}{\framebox(0.2,0.2){}}
\multiput(3,-12)(2,-6){1}{\circle*{0.1}}
\put(3,-12.2){\makebox(0,0)[ct]{\(2;2\)}}
\put(3,-12.6){\makebox(0,0)[ct]{\(S_1\)}}

\multiput(4,-12)(0,-2){2}{\line(0,-1){2}}
\put(3.9,-12.1){\framebox(0.2,0.2){}}
\put(3.95,-14.05){\framebox(0.1,0.1){}}
\put(4.1,-11.8){\makebox(0,0)[lb]{\(2;1\)}}
\put(4.1,-13.8){\makebox(0,0)[lb]{\(1;2\) (not coclass-settled)}}
\put(5.1,-14.8){\makebox(0,0)[lb]{\(3^{\text{rd}}\) bifurcation}}
\multiput(4,-14)(0,-4){1}{\line(-3,-2){3}}
\multiput(4,-14)(0,-4){1}{\line(-1,-1){2}}
\multiput(4,-14)(0,-4){1}{\line(-1,-2){1}}
\multiput(0.9,-16.1)(1,0){4}{\framebox(0.2,0.2){}}
\put(1,-16.2){\makebox(0,0)[ct]{\(1;1\)}}
\put(2,-16.2){\makebox(0,0)[ct]{\(1;2\)}}
\put(3,-16.2){\makebox(0,0)[ct]{\(1;3\)}}
\put(4,-16.2){\makebox(0,0)[ct]{\(1;4\)}}
\put(4,-14){\line(1,-2){2}}
\put(4,-14){\line(1,-4){1}}
\multiput(4.9,-18.1)(1,0){1}{\framebox(0.2,0.2){}}
\multiput(5,-18)(2,-6){1}{\circle*{0.1}}
\put(5,-18.2){\makebox(0,0)[ct]{\(2;2\)}}
\put(5,-18.6){\makebox(0,0)[ct]{\(S_2\)}}

\multiput(6,-18)(0,-2){2}{\line(0,-1){2}}
\put(5.9,-18.1){\framebox(0.2,0.2){}}
\put(5.95,-20.05){\framebox(0.1,0.1){}}
\put(6.1,-17.8){\makebox(0,0)[lb]{\(2;1\)}}
\put(6.1,-19.8){\makebox(0,0)[lb]{\(1;1\) (not coclass-settled)}}
\put(7.1,-20.8){\makebox(0,0)[lb]{\(4^{\text{th}}\) bifurcation}}
\multiput(6,-20)(0,-4){1}{\line(-3,-2){3}}
\multiput(6,-20)(0,-4){1}{\line(-1,-1){2}}
\multiput(6,-20)(0,-4){1}{\line(-1,-2){1}}
\multiput(2.9,-22.1)(1,0){4}{\framebox(0.2,0.2){}}
\put(3,-22.2){\makebox(0,0)[ct]{\(1;1\)}}
\put(4,-22.2){\makebox(0,0)[ct]{\(1;2\)}}
\put(5,-22.2){\makebox(0,0)[ct]{\(1;3\)}}
\put(6,-22.2){\makebox(0,0)[ct]{\(1;4\)}}
\put(6,-20){\line(1,-2){2}}
\put(6,-20){\line(1,-4){1}}
\multiput(6.9,-24.1)(1,0){1}{\framebox(0.2,0.2){}}
\multiput(7,-24)(2,-6){1}{\circle*{0.1}}
\put(7,-24.2){\makebox(0,0)[ct]{\(2;2\)}}
\put(7,-24.6){\makebox(0,0)[ct]{\(S_3\)}}

\multiput(7.9,-24.1)(0,-2){1}{\framebox(0.2,0.2){}}
\put(8.1,-23.8){\makebox(0,0)[lb]{\(2;1\)}}

\end{picture}

\end{figure}

}

\newpage


{\tiny

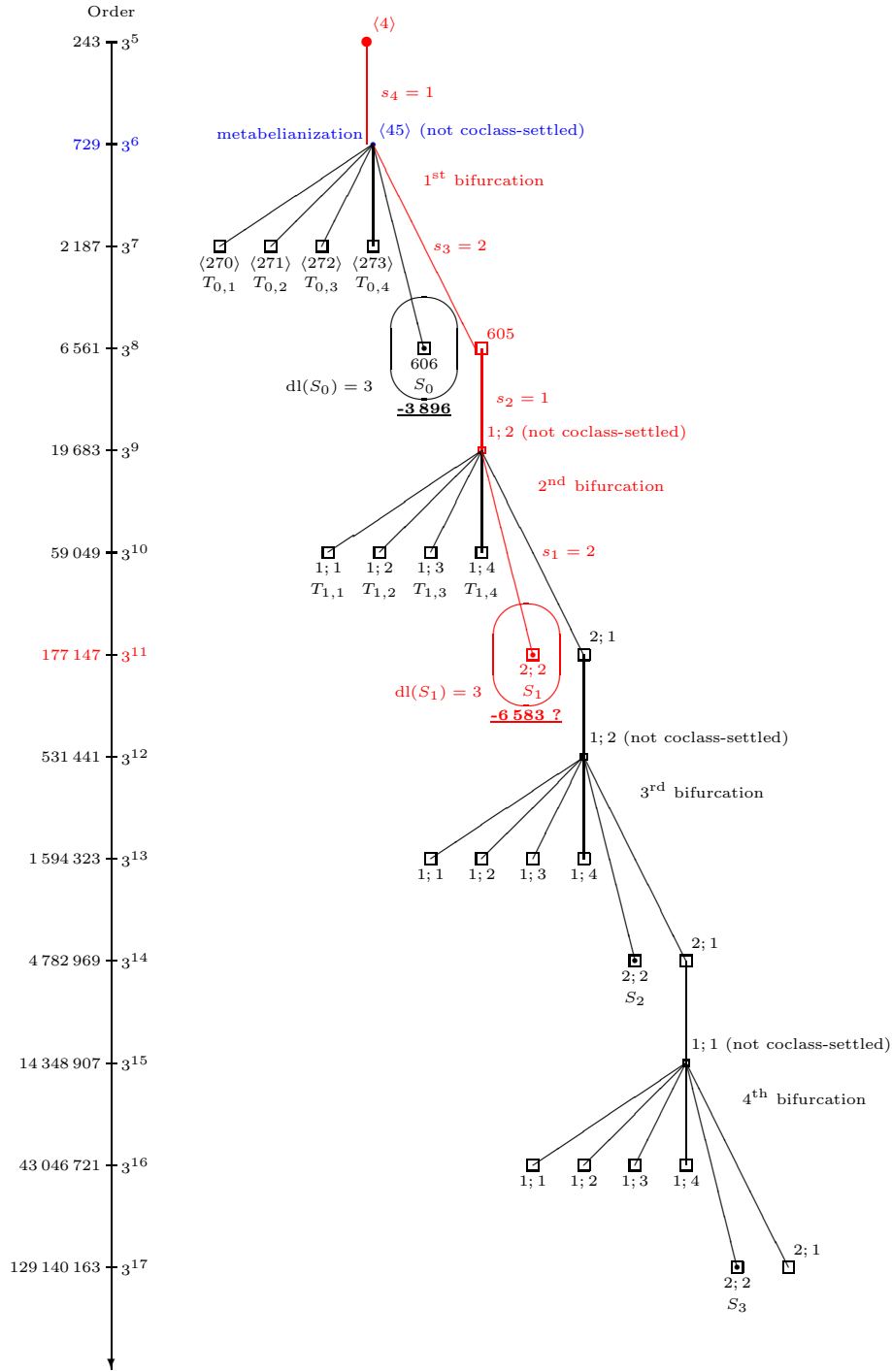
\begin{figure}[ht]
\caption{Extremal path to Schur \(\sigma\)-group \(S_1\), log ord \(11\), on purged tree \(\mathcal{T}_\ast(\langle 243,4\rangle)\)}
\label{fig:H4S1}


\setlength{\unitlength}{0.7cm}
\begin{picture}(18,26.5)(-6,-25.5)

\put(-5,0.5){\makebox(0,0)[cb]{Order}}
\put(-5,0){\line(0,-1){24}}
\multiput(-5.1,0)(0,-2){13}{\line(1,0){0.2}}
\put(-5.2,0){\makebox(0,0)[rc]{\(243\)}}
\put(-4.8,0){\makebox(0,0)[lc]{\(3^5\)}}
{\color{blue}
\put(-5.2,-2){\makebox(0,0)[rc]{\(729\)}}
\put(-4.8,-2){\makebox(0,0)[lc]{\(3^6\)}}
}
\put(-5.2,-4){\makebox(0,0)[rc]{\(2\,187\)}}
\put(-4.8,-4){\makebox(0,0)[lc]{\(3^7\)}}
\put(-5.2,-6){\makebox(0,0)[rc]{\(6\,561\)}}
\put(-4.8,-6){\makebox(0,0)[lc]{\(3^8\)}}
\put(-5.2,-8){\makebox(0,0)[rc]{\(19\,683\)}}
\put(-4.8,-8){\makebox(0,0)[lc]{\(3^9\)}}
\put(-5.2,-10){\makebox(0,0)[rc]{\(59\,049\)}}
\put(-4.8,-10){\makebox(0,0)[lc]{\(3^{10}\)}}
{\color{red}
\put(-5.2,-12){\makebox(0,0)[rc]{\(177\,147\)}}
\put(-4.8,-12){\makebox(0,0)[lc]{\(3^{11}\)}}
}
\put(-5.2,-14){\makebox(0,0)[rc]{\(531\,441\)}}
\put(-4.8,-14){\makebox(0,0)[lc]{\(3^{12}\)}}
\put(-5.2,-16){\makebox(0,0)[rc]{\(1\,594\,323\)}}
\put(-4.8,-16){\makebox(0,0)[lc]{\(3^{13}\)}}
\put(-5.2,-18){\makebox(0,0)[rc]{\(4\,782\,969\)}}
\put(-4.8,-18){\makebox(0,0)[lc]{\(3^{14}\)}}
\put(-5.2,-20){\makebox(0,0)[rc]{\(14\,348\,907\)}}
\put(-4.8,-20){\makebox(0,0)[lc]{\(3^{15}\)}}
\put(-5.2,-22){\makebox(0,0)[rc]{\(43\,046\,721\)}}
\put(-4.8,-22){\makebox(0,0)[lc]{\(3^{16}\)}}
\put(-5.2,-24){\makebox(0,0)[rc]{\(129\,140\,163\)}}
\put(-4.8,-24){\makebox(0,0)[lc]{\(3^{17}\)}}
\put(-5,-24){\vector(0,-1){2}}

{\color{red}
\multiput(0,0)(0,-2){1}{\circle*{0.2}}
\multiput(0,0)(0,-2){1}{\line(0,-1){2}}
\put(0.8,-1){\makebox(0,0)[cc]{\(s_4=1\)}}
\put(0.1,0.2){\makebox(0,0)[lb]{\(\langle 4\rangle\)}}
\put(1.1,-2.8){\makebox(0,0)[lb]{\(1^{\text{st}}\) bifurcation}}
}
{\color{blue}
\put(0,-2){\circle*{0.1}}
\put(0.1,-1.9){\makebox(0,0)[lb]{\(\langle 45\rangle\) (not coclass-settled)}}
\put(-0.2,-1.9){\makebox(0,0)[rb]{metabelianization}}
}
\multiput(0,-2)(0,-4){1}{\line(-3,-2){3}}
\multiput(0,-2)(0,-4){1}{\line(-1,-1){2}}
\multiput(0,-2)(0,-4){1}{\line(-1,-2){1}}
\multiput(0,-2)(0,-2){1}{\line(0,-1){2}}
\multiput(-3.1,-4.1)(1,0){4}{\framebox(0.2,0.2){}}
\put(-3,-4.2){\makebox(0,0)[ct]{\(\langle 270\rangle\)}}
\put(-2,-4.2){\makebox(0,0)[ct]{\(\langle 271\rangle\)}}
\put(-1,-4.2){\makebox(0,0)[ct]{\(\langle 272\rangle\)}}
\put(0,-4.2){\makebox(0,0)[ct]{\(\langle 273\rangle\)}}
\put(-3,-4.6){\makebox(0,0)[ct]{\(T_{0,1}\)}}
\put(-2,-4.6){\makebox(0,0)[ct]{\(T_{0,2}\)}}
\put(-1,-4.6){\makebox(0,0)[ct]{\(T_{0,3}\)}}
\put(0,-4.6){\makebox(0,0)[ct]{\(T_{0,4}\)}}
{\color{red}
\put(0,-2){\line(1,-2){2}}
\put(1.7,-4){\makebox(0,0)[cc]{\(s_3=2\)}}
}
\put(0,-2){\line(1,-4){1}}
\multiput(0.9,-6.1)(1,0){1}{\framebox(0.2,0.2){}}
\multiput(1,-6)(2,-6){1}{\circle*{0.1}}
\put(1,-6.2){\makebox(0,0)[ct]{\(606\)}}
\put(1,-6.6){\makebox(0,0)[ct]{\(S_0\)}}
\put(0,-6.6){\makebox(0,0)[rt]{\(\mathrm{dl}(S_0)=3\)}}






\put(1,-6){\oval(1.3,2.0)}
\put(1,-7.2){\makebox(0,0)[cc]{\underbar{\textbf{-3\,896}}}}
{\color{red}
\put(3,-12){\oval(1.3,2.0)}
\put(3,-13.2){\makebox(0,0)[cc]{\underbar{\textbf{-6\,583\ ?}}}}
}


{\color{red}
\multiput(2,-6)(0,-2){1}{\line(0,-1){2}}
\put(2.8,-7){\makebox(0,0)[cc]{\(s_2=1\)}}
\put(1.9,-6.1){\framebox(0.2,0.2){}}
\put(1.95,-8.05){\framebox(0.1,0.1){}}
\put(2.1,-5.8){\makebox(0,0)[lb]{\(605\)}}
\put(2.1,-7.8){\makebox(0,0)[lb]{\(1;2\) (not coclass-settled)}}
\put(3.1,-8.8){\makebox(0,0)[lb]{\(2^{\text{nd}}\) bifurcation}}
}
\multiput(2,-8)(0,-2){1}{\line(0,-1){2}}
\multiput(2,-8)(0,-4){1}{\line(-3,-2){3}}
\multiput(2,-8)(0,-4){1}{\line(-1,-1){2}}
\multiput(2,-8)(0,-4){1}{\line(-1,-2){1}}
\multiput(-1.1,-10.1)(1,0){4}{\framebox(0.2,0.2){}}
\put(-1,-10.2){\makebox(0,0)[ct]{\(1;1\)}}
\put(0,-10.2){\makebox(0,0)[ct]{\(1;2\)}}
\put(1,-10.2){\makebox(0,0)[ct]{\(1;3\)}}
\put(2,-10.2){\makebox(0,0)[ct]{\(1;4\)}}
\put(-1,-10.6){\makebox(0,0)[ct]{\(T_{1,1}\)}}
\put(0,-10.6){\makebox(0,0)[ct]{\(T_{1,2}\)}}
\put(1,-10.6){\makebox(0,0)[ct]{\(T_{1,3}\)}}
\put(2,-10.6){\makebox(0,0)[ct]{\(T_{1,4}\)}}
\put(2,-8){\line(1,-2){2}}
{\color{red}
\put(2,-8){\line(1,-4){1}}
\put(3.7,-10){\makebox(0,0)[cc]{\(s_1=2\)}}
\multiput(2.9,-12.1)(1,0){1}{\framebox(0.2,0.2){}}
\multiput(3,-12)(2,-6){1}{\circle*{0.1}}
\put(3,-12.2){\makebox(0,0)[ct]{\(2;2\)}}
\put(3,-12.6){\makebox(0,0)[ct]{\(S_1\)}}
\put(2,-12.6){\makebox(0,0)[rt]{\(\mathrm{dl}(S_1)=3\)}}
}
\multiput(4,-12)(0,-2){2}{\line(0,-1){2}}
\put(3.9,-12.1){\framebox(0.2,0.2){}}
\put(3.95,-14.05){\framebox(0.1,0.1){}}
\put(4.1,-11.8){\makebox(0,0)[lb]{\(2;1\)}}
\put(4.1,-13.8){\makebox(0,0)[lb]{\(1;2\) (not coclass-settled)}}
\put(5.1,-14.8){\makebox(0,0)[lb]{\(3^{\text{rd}}\) bifurcation}}
\multiput(4,-14)(0,-4){1}{\line(-3,-2){3}}
\multiput(4,-14)(0,-4){1}{\line(-1,-1){2}}
\multiput(4,-14)(0,-4){1}{\line(-1,-2){1}}
\multiput(0.9,-16.1)(1,0){4}{\framebox(0.2,0.2){}}
\put(1,-16.2){\makebox(0,0)[ct]{\(1;1\)}}
\put(2,-16.2){\makebox(0,0)[ct]{\(1;2\)}}
\put(3,-16.2){\makebox(0,0)[ct]{\(1;3\)}}
\put(4,-16.2){\makebox(0,0)[ct]{\(1;4\)}}
\put(4,-14){\line(1,-2){2}}
\put(4,-14){\line(1,-4){1}}
\multiput(4.9,-18.1)(1,0){1}{\framebox(0.2,0.2){}}
\multiput(5,-18)(2,-6){1}{\circle*{0.1}}
\put(5,-18.2){\makebox(0,0)[ct]{\(2;2\)}}
\put(5,-18.6){\makebox(0,0)[ct]{\(S_2\)}}

\multiput(6,-18)(0,-2){2}{\line(0,-1){2}}
\put(5.9,-18.1){\framebox(0.2,0.2){}}
\put(5.95,-20.05){\framebox(0.1,0.1){}}
\put(6.1,-17.8){\makebox(0,0)[lb]{\(2;1\)}}
\put(6.1,-19.8){\makebox(0,0)[lb]{\(1;1\) (not coclass-settled)}}
\put(7.1,-20.8){\makebox(0,0)[lb]{\(4^{\text{th}}\) bifurcation}}
\multiput(6,-20)(0,-4){1}{\line(-3,-2){3}}
\multiput(6,-20)(0,-4){1}{\line(-1,-1){2}}
\multiput(6,-20)(0,-4){1}{\line(-1,-2){1}}
\multiput(2.9,-22.1)(1,0){4}{\framebox(0.2,0.2){}}
\put(3,-22.2){\makebox(0,0)[ct]{\(1;1\)}}
\put(4,-22.2){\makebox(0,0)[ct]{\(1;2\)}}
\put(5,-22.2){\makebox(0,0)[ct]{\(1;3\)}}
\put(6,-22.2){\makebox(0,0)[ct]{\(1;4\)}}
\put(6,-20){\line(1,-2){2}}
\put(6,-20){\line(1,-4){1}}
\multiput(6.9,-24.1)(1,0){1}{\framebox(0.2,0.2){}}
\multiput(7,-24)(2,-6){1}{\circle*{0.1}}
\put(7,-24.2){\makebox(0,0)[ct]{\(2;2\)}}
\put(7,-24.6){\makebox(0,0)[ct]{\(S_3\)}}

\multiput(7.9,-24.1)(0,-2){1}{\framebox(0.2,0.2){}}
\put(8.1,-23.8){\makebox(0,0)[lb]{\(2;1\)}}

\end{picture}

\end{figure}

}

\newpage


{\tiny

\begin{figure}[ht]
\caption{Extremal path to Schur \(\sigma\)-group \(S_2\), log ord \(14\), on purged tree \(\mathcal{T}_\ast(\langle 243,4\rangle)\)}
\label{fig:H4S2}


\setlength{\unitlength}{0.7cm}
\begin{picture}(18,26.5)(-6,-25.5)

\put(-5,0.5){\makebox(0,0)[cb]{Order}}
\put(-5,0){\line(0,-1){24}}
\multiput(-5.1,0)(0,-2){13}{\line(1,0){0.2}}
\put(-5.2,0){\makebox(0,0)[rc]{\(243\)}}
\put(-4.8,0){\makebox(0,0)[lc]{\(3^5\)}}
{\color{blue}
\put(-5.2,-2){\makebox(0,0)[rc]{\(729\)}}
\put(-4.8,-2){\makebox(0,0)[lc]{\(3^6\)}}
}
\put(-5.2,-4){\makebox(0,0)[rc]{\(2\,187\)}}
\put(-4.8,-4){\makebox(0,0)[lc]{\(3^7\)}}
\put(-5.2,-6){\makebox(0,0)[rc]{\(6\,561\)}}
\put(-4.8,-6){\makebox(0,0)[lc]{\(3^8\)}}
\put(-5.2,-8){\makebox(0,0)[rc]{\(19\,683\)}}
\put(-4.8,-8){\makebox(0,0)[lc]{\(3^9\)}}
\put(-5.2,-10){\makebox(0,0)[rc]{\(59\,049\)}}
\put(-4.8,-10){\makebox(0,0)[lc]{\(3^{10}\)}}
\put(-5.2,-12){\makebox(0,0)[rc]{\(177\,147\)}}
\put(-4.8,-12){\makebox(0,0)[lc]{\(3^{11}\)}}
\put(-5.2,-14){\makebox(0,0)[rc]{\(531\,441\)}}
\put(-4.8,-14){\makebox(0,0)[lc]{\(3^{12}\)}}
\put(-5.2,-16){\makebox(0,0)[rc]{\(1\,594\,323\)}}
\put(-4.8,-16){\makebox(0,0)[lc]{\(3^{13}\)}}
{\color{red}
\put(-5.2,-18){\makebox(0,0)[rc]{\(4\,782\,969\)}}
\put(-4.8,-18){\makebox(0,0)[lc]{\(3^{14}\)}}
}
\put(-5.2,-20){\makebox(0,0)[rc]{\(14\,348\,907\)}}
\put(-4.8,-20){\makebox(0,0)[lc]{\(3^{15}\)}}
\put(-5.2,-22){\makebox(0,0)[rc]{\(43\,046\,721\)}}
\put(-4.8,-22){\makebox(0,0)[lc]{\(3^{16}\)}}
\put(-5.2,-24){\makebox(0,0)[rc]{\(129\,140\,163\)}}
\put(-4.8,-24){\makebox(0,0)[lc]{\(3^{17}\)}}
\put(-5,-24){\vector(0,-1){2}}

{\color{red}
\multiput(0,0)(0,-2){1}{\circle*{0.2}}
\multiput(0,0)(0,-2){1}{\line(0,-1){2}}
\put(0.8,-1){\makebox(0,0)[cc]{\(s_6=1\)}}
\put(0.1,0.2){\makebox(0,0)[lb]{\(\langle 4\rangle\)}}
\put(1.1,-2.8){\makebox(0,0)[lb]{\(1^{\text{st}}\) bifurcation}}
}
{\color{blue}
\put(0,-2){\circle*{0.1}}
\put(0.1,-1.9){\makebox(0,0)[lb]{\(\langle 45\rangle\) (not coclass-settled)}}
\put(-0.2,-1.9){\makebox(0,0)[rb]{metabelianization}}
}
\multiput(0,-2)(0,-4){1}{\line(-3,-2){3}}
\multiput(0,-2)(0,-4){1}{\line(-1,-1){2}}
\multiput(0,-2)(0,-4){1}{\line(-1,-2){1}}
\multiput(0,-2)(0,-2){1}{\line(0,-1){2}}
\multiput(-3.1,-4.1)(1,0){4}{\framebox(0.2,0.2){}}
\put(-3,-4.2){\makebox(0,0)[ct]{\(\langle 270\rangle\)}}
\put(-2,-4.2){\makebox(0,0)[ct]{\(\langle 271\rangle\)}}
\put(-1,-4.2){\makebox(0,0)[ct]{\(\langle 272\rangle\)}}
\put(0,-4.2){\makebox(0,0)[ct]{\(\langle 273\rangle\)}}
\put(-3,-4.6){\makebox(0,0)[ct]{\(T_{0,1}\)}}
\put(-2,-4.6){\makebox(0,0)[ct]{\(T_{0,2}\)}}
\put(-1,-4.6){\makebox(0,0)[ct]{\(T_{0,3}\)}}
\put(0,-4.6){\makebox(0,0)[ct]{\(T_{0,4}\)}}
{\color{red}
\put(0,-2){\line(1,-2){2}}
\put(1.7,-4){\makebox(0,0)[cc]{\(s_5=2\)}}
}
\put(0,-2){\line(1,-4){1}}
\multiput(0.9,-6.1)(1,0){1}{\framebox(0.2,0.2){}}
\multiput(1,-6)(2,-6){1}{\circle*{0.1}}
\put(1,-6.2){\makebox(0,0)[ct]{\(606\)}}
\put(1,-6.6){\makebox(0,0)[ct]{\(S_0\)}}
\put(0,-6.6){\makebox(0,0)[rt]{\(\mathrm{dl}(S_0)=3\)}}






\put(1,-6){\oval(1.3,2.0)}
\put(1,-7.2){\makebox(0,0)[cc]{\underbar{\textbf{-3\,896}}}}

\put(3,-12){\oval(1.3,2.0)}
\put(3,-13.2){\makebox(0,0)[cc]{\underbar{\textbf{-6\,583\ ?}}}}


{\color{red}
\multiput(2,-6)(0,-2){1}{\line(0,-1){2}}
\put(2.8,-7){\makebox(0,0)[cc]{\(s_4=1\)}}
\put(1.9,-6.1){\framebox(0.2,0.2){}}
\put(1.95,-8.05){\framebox(0.1,0.1){}}
\put(2.1,-5.8){\makebox(0,0)[lb]{\(605\)}}
\put(2.1,-7.8){\makebox(0,0)[lb]{\(1;2\) (not coclass-settled)}}
\put(3.1,-8.8){\makebox(0,0)[lb]{\(2^{\text{nd}}\) bifurcation}}
}
\multiput(2,-8)(0,-2){1}{\line(0,-1){2}}
\multiput(2,-8)(0,-4){1}{\line(-3,-2){3}}
\multiput(2,-8)(0,-4){1}{\line(-1,-1){2}}
\multiput(2,-8)(0,-4){1}{\line(-1,-2){1}}
\multiput(-1.1,-10.1)(1,0){4}{\framebox(0.2,0.2){}}
\put(-1,-10.2){\makebox(0,0)[ct]{\(1;1\)}}
\put(0,-10.2){\makebox(0,0)[ct]{\(1;2\)}}
\put(1,-10.2){\makebox(0,0)[ct]{\(1;3\)}}
\put(2,-10.2){\makebox(0,0)[ct]{\(1;4\)}}
\put(-1,-10.6){\makebox(0,0)[ct]{\(T_{1,1}\)}}
\put(0,-10.6){\makebox(0,0)[ct]{\(T_{1,2}\)}}
\put(1,-10.6){\makebox(0,0)[ct]{\(T_{1,3}\)}}
\put(2,-10.6){\makebox(0,0)[ct]{\(T_{1,4}\)}}
{\color{red}
\put(2,-8){\line(1,-2){2}}
\put(3.7,-10){\makebox(0,0)[cc]{\(s_3=2\)}}
}
\put(2,-8){\line(1,-4){1}}
\multiput(2.9,-12.1)(1,0){1}{\framebox(0.2,0.2){}}
\multiput(3,-12)(2,-6){1}{\circle*{0.1}}
\put(3,-12.2){\makebox(0,0)[ct]{\(2;2\)}}
\put(3,-12.6){\makebox(0,0)[ct]{\(S_1\)}}
\put(2,-12.6){\makebox(0,0)[rt]{\(\mathrm{dl}(S_1)=3\)}}

{\color{red}
\multiput(4,-12)(0,-2){1}{\line(0,-1){2}}
\put(4.8,-13){\makebox(0,0)[cc]{\(s_2=1\)}}
\put(3.9,-12.1){\framebox(0.2,0.2){}}
\put(3.95,-14.05){\framebox(0.1,0.1){}}
\put(4.1,-11.8){\makebox(0,0)[lb]{\(2;1\)}}
\put(4.1,-13.8){\makebox(0,0)[lb]{\(1;2\) (not coclass-settled)}}
\put(5.1,-14.8){\makebox(0,0)[lb]{\(3^{\text{rd}}\) bifurcation}}
}
\multiput(4,-14)(0,-4){1}{\line(-3,-2){3}}
\multiput(4,-14)(0,-4){1}{\line(-1,-1){2}}
\multiput(4,-14)(0,-4){1}{\line(-1,-2){1}}
\multiput(4,-14)(0,-2){1}{\line(0,-1){2}}
\multiput(0.9,-16.1)(1,0){4}{\framebox(0.2,0.2){}}
\put(1,-16.2){\makebox(0,0)[ct]{\(1;1\)}}
\put(2,-16.2){\makebox(0,0)[ct]{\(1;2\)}}
\put(3,-16.2){\makebox(0,0)[ct]{\(1;3\)}}
\put(4,-16.2){\makebox(0,0)[ct]{\(1;4\)}}
\put(4,-14){\line(1,-2){2}}
{\color{red}
\put(4,-14){\line(1,-4){1}}
\put(5.7,-16){\makebox(0,0)[cc]{\(s_1=2\)}}
\multiput(4.9,-18.1)(1,0){1}{\framebox(0.2,0.2){}}
\multiput(5,-18)(2,-6){1}{\circle*{0.1}}
\put(5,-18.2){\makebox(0,0)[ct]{\(2;2\)}}
\put(5,-18.6){\makebox(0,0)[ct]{\(S_2\)}}
\put(4,-18.6){\makebox(0,0)[rt]{\(\mathrm{dl}(S_2)=3\)}}
}
\multiput(6,-18)(0,-2){2}{\line(0,-1){2}}
\put(5.9,-18.1){\framebox(0.2,0.2){}}
\put(5.95,-20.05){\framebox(0.1,0.1){}}
\put(6.1,-17.8){\makebox(0,0)[lb]{\(2;1\)}}
\put(6.1,-19.8){\makebox(0,0)[lb]{\(1;1\) (not coclass-settled)}}
\put(7.1,-20.8){\makebox(0,0)[lb]{\(4^{\text{th}}\) bifurcation}}
\multiput(6,-20)(0,-4){1}{\line(-3,-2){3}}
\multiput(6,-20)(0,-4){1}{\line(-1,-1){2}}
\multiput(6,-20)(0,-4){1}{\line(-1,-2){1}}
\multiput(2.9,-22.1)(1,0){4}{\framebox(0.2,0.2){}}
\put(3,-22.2){\makebox(0,0)[ct]{\(1;1\)}}
\put(4,-22.2){\makebox(0,0)[ct]{\(1;2\)}}
\put(5,-22.2){\makebox(0,0)[ct]{\(1;3\)}}
\put(6,-22.2){\makebox(0,0)[ct]{\(1;4\)}}
\put(6,-20){\line(1,-2){2}}
\put(6,-20){\line(1,-4){1}}
\multiput(6.9,-24.1)(1,0){1}{\framebox(0.2,0.2){}}
\multiput(7,-24)(2,-6){1}{\circle*{0.1}}
\put(7,-24.2){\makebox(0,0)[ct]{\(2;2\)}}
\put(7,-24.6){\makebox(0,0)[ct]{\(S_3\)}}

\multiput(7.9,-24.1)(0,-2){1}{\framebox(0.2,0.2){}}
\put(8.1,-23.8){\makebox(0,0)[lb]{\(2;1\)}}

\end{picture}

\end{figure}
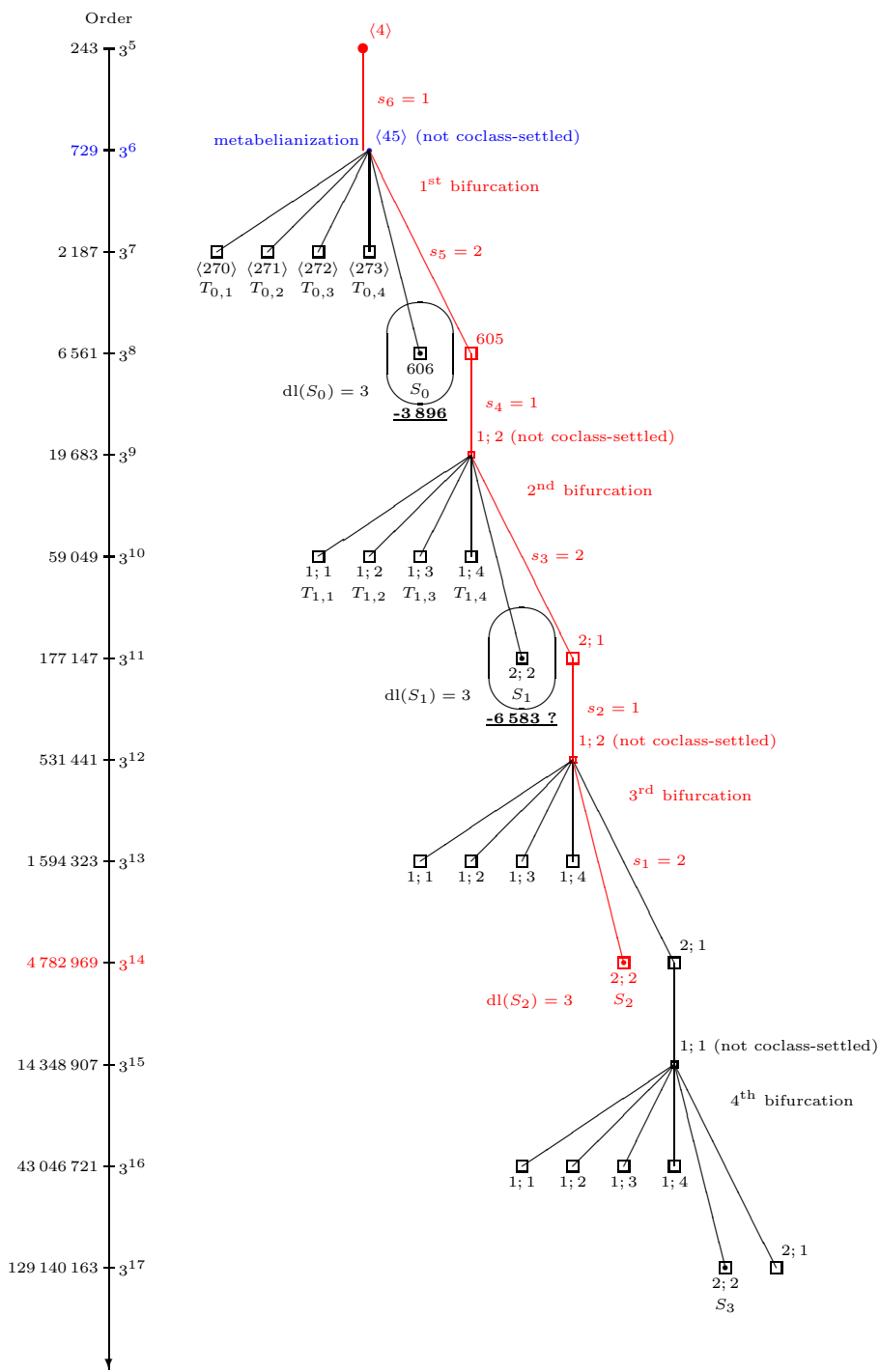

}

\newpage


{\tiny

\begin{figure}[ht]
\caption{Extremal path to Schur \(\sigma\)-group \(S_3\), log ord \(17\), on purged tree \(\mathcal{T}_\ast(\langle 243,4\rangle)\)}
\label{fig:H4S3}


\setlength{\unitlength}{0.7cm}
\begin{picture}(18,26.5)(-6,-25.5)

\put(-5,0.5){\makebox(0,0)[cb]{Order}}
\put(-5,0){\line(0,-1){24}}
\multiput(-5.1,0)(0,-2){13}{\line(1,0){0.2}}
\put(-5.2,0){\makebox(0,0)[rc]{\(243\)}}
\put(-4.8,0){\makebox(0,0)[lc]{\(3^5\)}}
{\color{blue}
\put(-5.2,-2){\makebox(0,0)[rc]{\(729\)}}
\put(-4.8,-2){\makebox(0,0)[lc]{\(3^6\)}}
}
\put(-5.2,-4){\makebox(0,0)[rc]{\(2\,187\)}}
\put(-4.8,-4){\makebox(0,0)[lc]{\(3^7\)}}
\put(-5.2,-6){\makebox(0,0)[rc]{\(6\,561\)}}
\put(-4.8,-6){\makebox(0,0)[lc]{\(3^8\)}}
\put(-5.2,-8){\makebox(0,0)[rc]{\(19\,683\)}}
\put(-4.8,-8){\makebox(0,0)[lc]{\(3^9\)}}
\put(-5.2,-10){\makebox(0,0)[rc]{\(59\,049\)}}
\put(-4.8,-10){\makebox(0,0)[lc]{\(3^{10}\)}}
\put(-5.2,-12){\makebox(0,0)[rc]{\(177\,147\)}}
\put(-4.8,-12){\makebox(0,0)[lc]{\(3^{11}\)}}
\put(-5.2,-14){\makebox(0,0)[rc]{\(531\,441\)}}
\put(-4.8,-14){\makebox(0,0)[lc]{\(3^{12}\)}}
\put(-5.2,-16){\makebox(0,0)[rc]{\(1\,594\,323\)}}
\put(-4.8,-16){\makebox(0,0)[lc]{\(3^{13}\)}}
\put(-5.2,-18){\makebox(0,0)[rc]{\(4\,782\,969\)}}
\put(-4.8,-18){\makebox(0,0)[lc]{\(3^{14}\)}}
\put(-5.2,-20){\makebox(0,0)[rc]{\(14\,348\,907\)}}
\put(-4.8,-20){\makebox(0,0)[lc]{\(3^{15}\)}}
\put(-5.2,-22){\makebox(0,0)[rc]{\(43\,046\,721\)}}
\put(-4.8,-22){\makebox(0,0)[lc]{\(3^{16}\)}}
{\color{red}
\put(-5.2,-24){\makebox(0,0)[rc]{\(129\,140\,163\)}}
\put(-4.8,-24){\makebox(0,0)[lc]{\(3^{17}\)}}
}
\put(-5,-24){\vector(0,-1){2}}

{\color{red}
\multiput(0,0)(0,-2){1}{\circle*{0.2}}
\multiput(0,0)(0,-2){1}{\line(0,-1){2}}
\put(0.8,-1){\makebox(0,0)[cc]{\(s_8=1\)}}
\put(0.1,0.2){\makebox(0,0)[lb]{\(\langle 4\rangle\)}}
\put(1.1,-2.8){\makebox(0,0)[lb]{\(1^{\text{st}}\) bifurcation}}
}
{\color{blue}
\put(0,-2){\circle*{0.1}}
\put(0.1,-1.9){\makebox(0,0)[lb]{\(\langle 45\rangle\) (not coclass-settled)}}
\put(-0.2,-1.9){\makebox(0,0)[rb]{metabelianization}}
}
\multiput(0,-2)(0,-4){1}{\line(-3,-2){3}}
\multiput(0,-2)(0,-4){1}{\line(-1,-1){2}}
\multiput(0,-2)(0,-4){1}{\line(-1,-2){1}}
\multiput(0,-2)(0,-2){1}{\line(0,-1){2}}
\multiput(-3.1,-4.1)(1,0){4}{\framebox(0.2,0.2){}}
\put(-3,-4.2){\makebox(0,0)[ct]{\(\langle 270\rangle\)}}
\put(-2,-4.2){\makebox(0,0)[ct]{\(\langle 271\rangle\)}}
\put(-1,-4.2){\makebox(0,0)[ct]{\(\langle 272\rangle\)}}
\put(0,-4.2){\makebox(0,0)[ct]{\(\langle 273\rangle\)}}
\put(-3,-4.6){\makebox(0,0)[ct]{\(T_{0,1}\)}}
\put(-2,-4.6){\makebox(0,0)[ct]{\(T_{0,2}\)}}
\put(-1,-4.6){\makebox(0,0)[ct]{\(T_{0,3}\)}}
\put(0,-4.6){\makebox(0,0)[ct]{\(T_{0,4}\)}}
{\color{red}
\put(0,-2){\line(1,-2){2}}
\put(1.7,-4){\makebox(0,0)[cc]{\(s_7=2\)}}
}
\put(0,-2){\line(1,-4){1}}
\multiput(0.9,-6.1)(1,0){1}{\framebox(0.2,0.2){}}
\multiput(1,-6)(2,-6){1}{\circle*{0.1}}
\put(1,-6.2){\makebox(0,0)[ct]{\(606\)}}
\put(1,-6.6){\makebox(0,0)[ct]{\(S_0\)}}
\put(0,-6.6){\makebox(0,0)[rt]{\(\mathrm{dl}(S_0)=3\)}}






\put(1,-6){\oval(1.3,2.0)}
\put(1,-7.2){\makebox(0,0)[cc]{\underbar{\textbf{-3\,896}}}}

\put(3,-12){\oval(1.3,2.0)}
\put(3,-13.2){\makebox(0,0)[cc]{\underbar{\textbf{-6\,583\ ?}}}}


{\color{red}
\multiput(2,-6)(0,-2){1}{\line(0,-1){2}}
\put(2.8,-7){\makebox(0,0)[cc]{\(s_6=1\)}}
\put(1.9,-6.1){\framebox(0.2,0.2){}}
\put(1.95,-8.05){\framebox(0.1,0.1){}}
\put(2.1,-5.8){\makebox(0,0)[lb]{\(605\)}}
\put(2.1,-7.8){\makebox(0,0)[lb]{\(1;2\) (not coclass-settled)}}
\put(3.1,-8.8){\makebox(0,0)[lb]{\(2^{\text{nd}}\) bifurcation}}
}
\multiput(2,-8)(0,-2){1}{\line(0,-1){2}}
\multiput(2,-8)(0,-4){1}{\line(-3,-2){3}}
\multiput(2,-8)(0,-4){1}{\line(-1,-1){2}}
\multiput(2,-8)(0,-4){1}{\line(-1,-2){1}}
\multiput(-1.1,-10.1)(1,0){4}{\framebox(0.2,0.2){}}
\put(-1,-10.2){\makebox(0,0)[ct]{\(1;1\)}}
\put(0,-10.2){\makebox(0,0)[ct]{\(1;2\)}}
\put(1,-10.2){\makebox(0,0)[ct]{\(1;3\)}}
\put(2,-10.2){\makebox(0,0)[ct]{\(1;4\)}}
\put(-1,-10.6){\makebox(0,0)[ct]{\(T_{1,1}\)}}
\put(0,-10.6){\makebox(0,0)[ct]{\(T_{1,2}\)}}
\put(1,-10.6){\makebox(0,0)[ct]{\(T_{1,3}\)}}
\put(2,-10.6){\makebox(0,0)[ct]{\(T_{1,4}\)}}
{\color{red}
\put(2,-8){\line(1,-2){2}}
\put(3.7,-10){\makebox(0,0)[cc]{\(s_5=2\)}}
}
\put(2,-8){\line(1,-4){1}}
\multiput(2.9,-12.1)(1,0){1}{\framebox(0.2,0.2){}}
\multiput(3,-12)(2,-6){1}{\circle*{0.1}}
\put(3,-12.2){\makebox(0,0)[ct]{\(2;2\)}}
\put(3,-12.6){\makebox(0,0)[ct]{\(S_1\)}}
\put(2,-12.6){\makebox(0,0)[rt]{\(\mathrm{dl}(S_1)=3\)}}

{\color{red}
\multiput(4,-12)(0,-2){1}{\line(0,-1){2}}
\put(4.8,-13){\makebox(0,0)[cc]{\(s_4=1\)}}
\put(3.9,-12.1){\framebox(0.2,0.2){}}
\put(3.95,-14.05){\framebox(0.1,0.1){}}
\put(4.1,-11.8){\makebox(0,0)[lb]{\(2;1\)}}
\put(4.1,-13.8){\makebox(0,0)[lb]{\(1;2\) (not coclass-settled)}}
\put(5.1,-14.8){\makebox(0,0)[lb]{\(3^{\text{rd}}\) bifurcation}}
}
\multiput(4,-14)(0,-4){1}{\line(-3,-2){3}}
\multiput(4,-14)(0,-4){1}{\line(-1,-1){2}}
\multiput(4,-14)(0,-4){1}{\line(-1,-2){1}}
\multiput(4,-14)(0,-2){1}{\line(0,-1){2}}
\multiput(0.9,-16.1)(1,0){4}{\framebox(0.2,0.2){}}
\put(1,-16.2){\makebox(0,0)[ct]{\(1;1\)}}
\put(2,-16.2){\makebox(0,0)[ct]{\(1;2\)}}
\put(3,-16.2){\makebox(0,0)[ct]{\(1;3\)}}
\put(4,-16.2){\makebox(0,0)[ct]{\(1;4\)}}
{\color{red}
\put(4,-14){\line(1,-2){2}}
\put(5.7,-16){\makebox(0,0)[cc]{\(s_3=2\)}}
}
\put(4,-14){\line(1,-4){1}}
\multiput(4.9,-18.1)(1,0){1}{\framebox(0.2,0.2){}}
\multiput(5,-18)(2,-6){1}{\circle*{0.1}}
\put(5,-18.2){\makebox(0,0)[ct]{\(2;2\)}}
\put(5,-18.6){\makebox(0,0)[ct]{\(S_2\)}}
\put(4,-18.6){\makebox(0,0)[rt]{\(\mathrm{dl}(S_2)=3\)}}

{\color{red}
\multiput(6,-18)(0,-2){1}{\line(0,-1){2}}
\put(6.8,-19){\makebox(0,0)[cc]{\(s_2=1\)}}
\put(5.9,-18.1){\framebox(0.2,0.2){}}
\put(5.95,-20.05){\framebox(0.1,0.1){}}
\put(6.1,-17.8){\makebox(0,0)[lb]{\(2;1\)}}
\put(6.1,-19.8){\makebox(0,0)[lb]{\(1;1\) (not coclass-settled)}}
\put(7.1,-20.8){\makebox(0,0)[lb]{\(4^{\text{th}}\) bifurcation}}
}
\multiput(6,-20)(0,-4){1}{\line(-3,-2){3}}
\multiput(6,-20)(0,-4){1}{\line(-1,-1){2}}
\multiput(6,-20)(0,-4){1}{\line(-1,-2){1}}
\multiput(6,-20)(0,-2){1}{\line(0,-1){2}}
\multiput(2.9,-22.1)(1,0){4}{\framebox(0.2,0.2){}}
\put(3,-22.2){\makebox(0,0)[ct]{\(1;1\)}}
\put(4,-22.2){\makebox(0,0)[ct]{\(1;2\)}}
\put(5,-22.2){\makebox(0,0)[ct]{\(1;3\)}}
\put(6,-22.2){\makebox(0,0)[ct]{\(1;4\)}}
\put(6,-20){\line(1,-2){2}}
{\color{red}
\put(6,-20){\line(1,-4){1}}
\put(7.7,-22){\makebox(0,0)[cc]{\(s_1=2\)}}
\multiput(6.9,-24.1)(1,0){1}{\framebox(0.2,0.2){}}
\multiput(7,-24)(2,-6){1}{\circle*{0.1}}
\put(7,-24.2){\makebox(0,0)[ct]{\(2;2\)}}
\put(7,-24.6){\makebox(0,0)[ct]{\(S_3\)}}
\put(6,-24.6){\makebox(0,0)[rt]{\(\mathrm{dl}(S_3)=4\)}}
}
\multiput(7.9,-24.1)(0,-2){1}{\framebox(0.2,0.2){}}
\put(8.1,-23.8){\makebox(0,0)[lb]{\(2;1\)}}

\end{picture}

\end{figure}
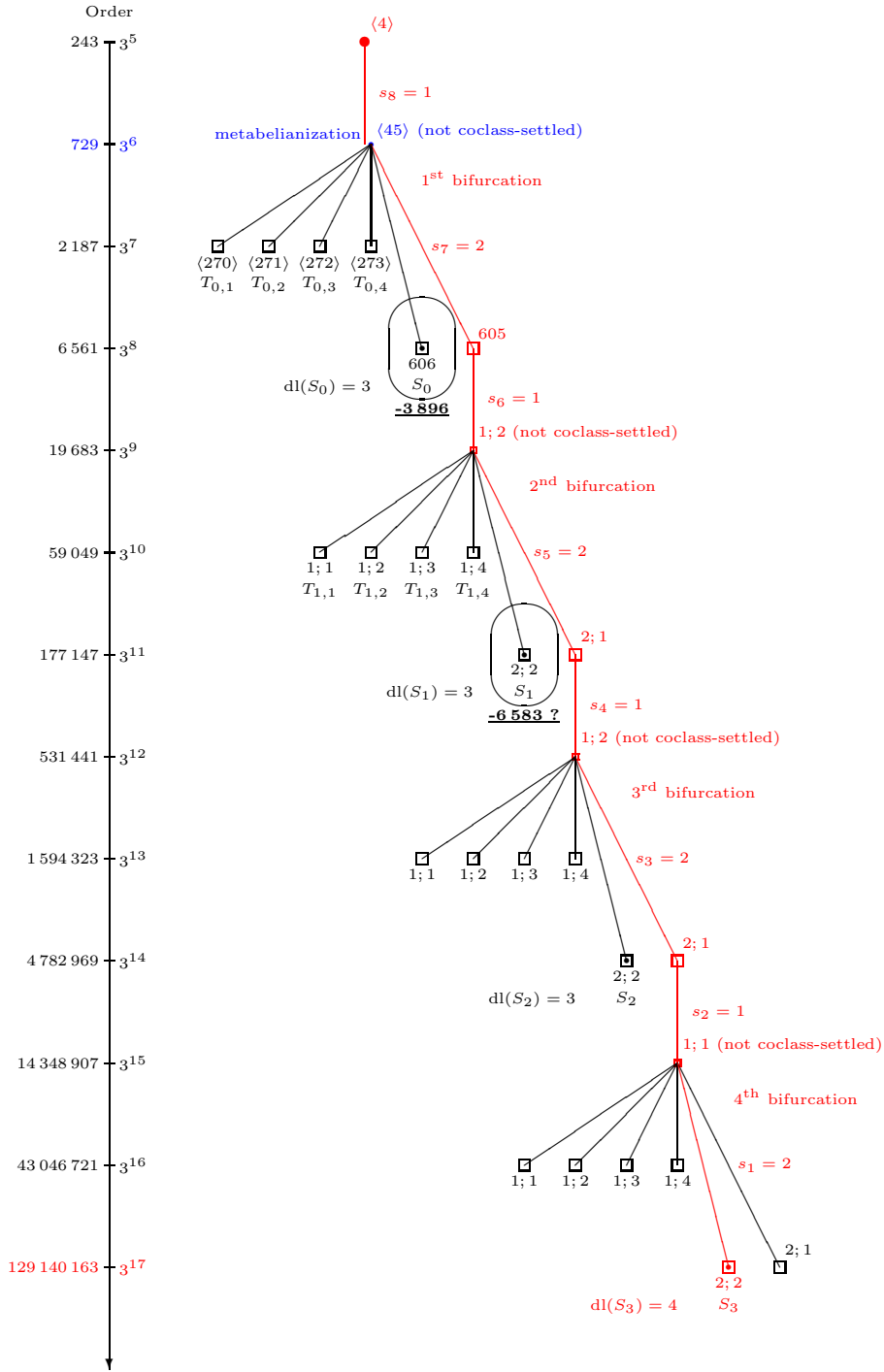

}

\newpage


\subsection{Infinite series of Schur \(\sigma\)-groups \(G\) with derived length three}
\label{ss:ThreeStage}
(\(\mathrm{dl}(G)=3\))

{\color{blue}
\begin{proposition}
\label{prp:SectionE}
Among the finite \(3\)-groups \(V\) with abelian quotient invariants \(V/V^\prime\simeq (3,3)\)
and transfer kernel type in Section \(\mathrm{E}\),
there exist \textbf{precisely six} \textbf{infinite sequences} of Schur \(\sigma\)-groups
\((S^{(t)}_n)_{n\ge 0}\), \(1\le t\le 6\),
sharing common derived length \(\mathrm{dl}(S^{(t)}_n)=3\).
Their logarithmic order, nilpotency class and coclass are given by the following laws:
\begin{equation}
\label{eqn:SectionE}
\mathrm{lo}(S^{(t)}_n)=8+3n, \quad \mathrm{cl}(S^{(t)}_n)=5+2n, \quad \mathrm{cc}(S^{(t)}_n)=3+n, \text{ for all } n\ge 0.
\end{equation}
The TKT is
\(\mathrm{E}.6\) for \(t=1\),
\(\mathrm{E}.14\) for \(t=2,3\),
\(\mathrm{E}.8\) for \(t=4\), and
\(\mathrm{E}.9\) for \(t=5,6\).
\end{proposition}
}

\begin{proof}
The existence and the constant derived length of the infinite sequences was proved by Newman and ourselves
\cite{Ma2018b}.
In fact, \(\mathrm{dl}=3\) is due to the inclusion of the (non-trivial) second derived subgroup in the centre.
The deterministic laws for invariants can be deduced in the same way as in the proof of Proposition
\ref{prp:SporadicSectionH},
replacing the root \(\langle 243,4\rangle\) of the purged descendant tree
by either \(R:=\langle 243,6\rangle\) or \(R:=\langle 243,8\rangle\).
As before, the tree is not a coclass tree, but it has an infinite main trunk
with strictly alternating step sizes \(s=1\) and \(s=2\),
and thus contains periodic bifurcations to higher coclass
\cite{Ma2015a,Ma2015c}.
Since \(S^{(t)}_n=R(-\#1;a_j-\#2;b_j)_{1\le j\le n+1}\) with certain \(1\le a_j\le 4\), \(1\le b_j\le 6\),
for each \(n\ge 0\),
we arrive at the same laws as in Proposition
\ref{prp:SporadicSectionH}.
\end{proof}


{\color{blue}
\begin{theorem}
\label{thm:SectionE}
The \(3\)-class field tower of an imaginary quadratic field
with \(3\)-capitulation type in section \(\mathrm{E}\)
is non-metabelian with precise length \(\ell_3(K)=3\).
The Galois group \(G\) of the ground state
satisfies the extremal property with
\(s_1=2\), \(s_2=1\) and \(s_3=2\).
Generally, the extremal path property of all excited states
is satisfied with strictly alternating step sizes \(s=2\) and \(s=1\).
\end{theorem}
}

\begin{proof}
The length was proved for the ground state by Boston, Bush and ourselves in
\cite{BuMa},
and for all excited states by Newman and ourselves in
\cite{Ma2018b}.
The extremal property is a consequence of Tables
\ref{tbl:E6},
resp.
\ref{tbl:E8},
and analogous tables for the remaining types in section E;
for type E.14 with \(\langle 6561,616\rangle\) replaced by \(\langle 6561,617\vert 618\rangle\),
resp.
for type E.9 with \(\langle 6561,622\rangle\) replaced by \(\langle 6561,620\vert 624\rangle\).
The extremal property for excited states was proved by ourselves in
\cite{Ma2015a,Ma2015c}.
\end{proof}


\renewcommand{\arraystretch}{1.1}

\begin{table}[ht]
\caption{Root path of \(G\) for the ground state of transfer kernel type E.6}
\label{tbl:E6}
\begin{center}
\begin{tabular}{|c|c||c|c|c|}
\hline
 Ancestor     & Id                          & \((\nu,\mu)\)                   & \((N_s/C_s)_{1\le s\le\nu}\)        & TKT  \\
\hline
 \(\pi^3(G)\) & \(\langle 27,3\rangle\)     & \(({\color{red}\mathbf{2,4}})\) & \((4/1,{\color{red}\mathbf{7/5}})\) & a.1  \\
 \(\pi^2(G)\) & \(\langle 243,6\rangle\)    & \(({\color{red}\mathbf{1,3}})\) & \(({\color{red}\mathbf{4/4}})\)     & c.18 \\
 \(\pi(G)\)   & \(\langle 729,49\rangle\)   & \(({\color{red}\mathbf{2,4}})\) & \((8/3,{\color{red}\mathbf{6/3}})\) & c.18 \\
 \(G\)        & \(\langle 6561,616\rangle\) & \(({\color{red}\mathbf{0,2}})\) &                                     & E.6  \\
\hline
\end{tabular}
\end{center}
\end{table}

\renewcommand{\arraystretch}{1.1}

\begin{table}[ht]
\caption{Root path of \(G\) for the ground state of transfer kernel type E.8}
\label{tbl:E8}
\begin{center}
\begin{tabular}{|c|c||c|c|c|}
\hline
 Ancestor     & Id                          & \((\nu,\mu)\)                   & \((N_s/C_s)_{1\le s\le\nu}\)        & TKT  \\
\hline
 \(\pi^3(G)\) & \(\langle 27,3\rangle\)     & \(({\color{red}\mathbf{2,4}})\) & \((4/1,{\color{red}\mathbf{7/5}})\) & a.1  \\
 \(\pi^2(G)\) & \(\langle 243,8\rangle\)    & \(({\color{red}\mathbf{1,3}})\) & \(({\color{red}\mathbf{4/4}})\)     & c.21 \\
 \(\pi(G)\)   & \(\langle 729,54\rangle\)   & \(({\color{red}\mathbf{2,4}})\) & \((8/3,{\color{red}\mathbf{6/3}})\) & c.21 \\
 \(G\)        & \(\langle 6561,622\rangle\) & \(({\color{red}\mathbf{0,2}})\) &                                     & E.8  \\
\hline
\end{tabular}
\end{center}
\end{table}

Root path of Schur \(\sigma\)-groups
\(S^{(t)}_0\), resp. \(S^{(t)}_1\), \(S^{(t)}_2\),
\(4\le t\le 6\),
is shown in Figure
\ref{fig:E8GS},
resp.
\ref{fig:E8ES1},
\ref{fig:E8ES2}.

\newpage


{\tiny

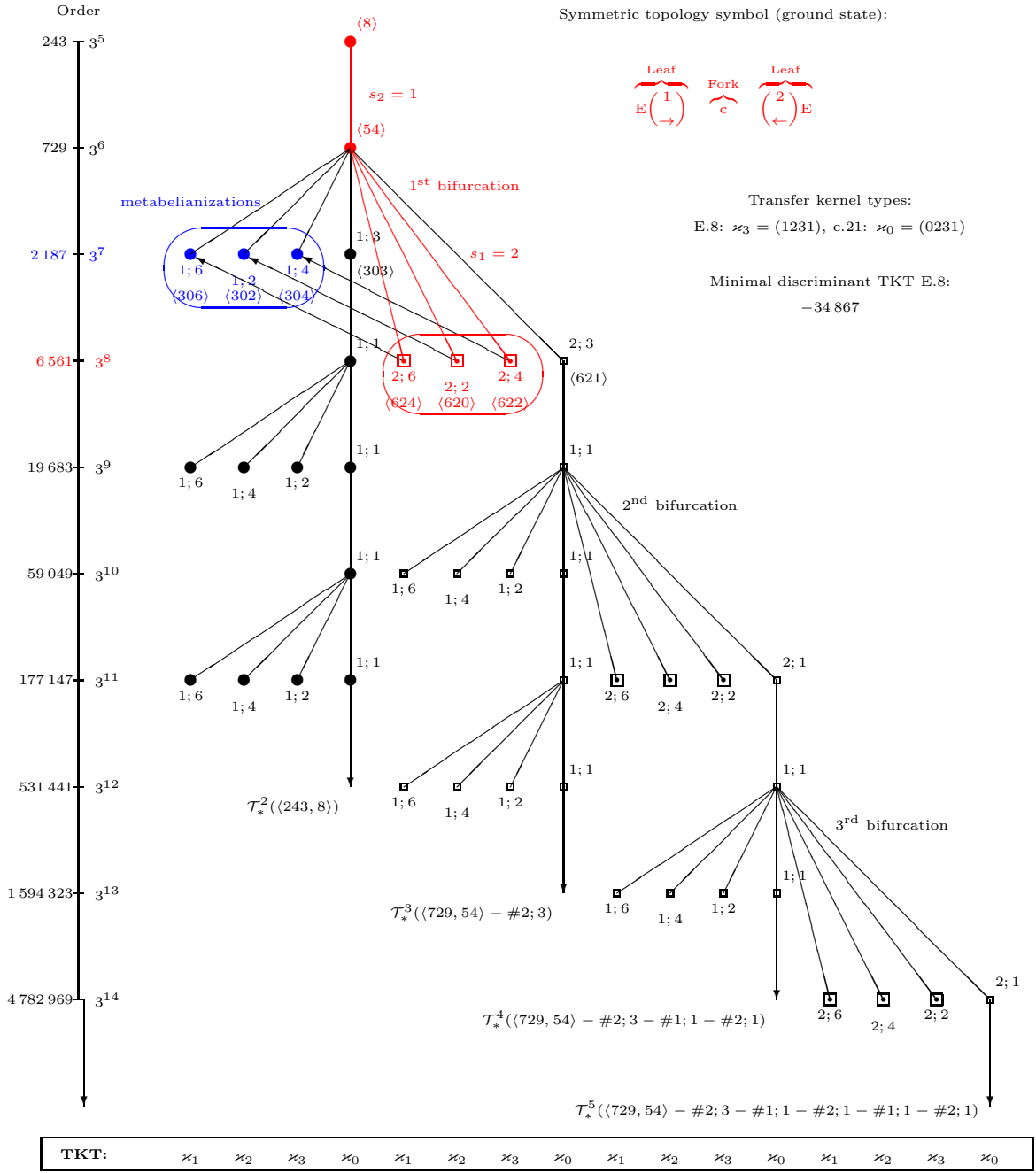
\begin{figure}[hb]
\caption{Extremal paths to Schur \(\sigma\)-groups, log ord \(8\), on purged tree \(\mathcal{T}_\ast(\langle 243,8\rangle)\)}
\label{fig:E8GS}


\setlength{\unitlength}{0.8cm}
\begin{picture}(18,22)(-6,-21)


\put(-5,0.5){\makebox(0,0)[cb]{Order}}
\put(-5,0){\line(0,-1){18}}
\multiput(-5.1,0)(0,-2){10}{\line(1,0){0.2}}
\put(-5.2,0){\makebox(0,0)[rc]{\(243\)}}
\put(-4.8,0){\makebox(0,0)[lc]{\(3^5\)}}
\put(-5.2,-2){\makebox(0,0)[rc]{\(729\)}}
\put(-4.8,-2){\makebox(0,0)[lc]{\(3^6\)}}
{\color{blue}
\put(-5.2,-4){\makebox(0,0)[rc]{\(2\,187\)}}
\put(-4.8,-4){\makebox(0,0)[lc]{\(3^7\)}}
}
{\color{red}
\put(-5.2,-6){\makebox(0,0)[rc]{\(6\,561\)}}
\put(-4.8,-6){\makebox(0,0)[lc]{\(3^8\)}}
}
\put(-5.2,-8){\makebox(0,0)[rc]{\(19\,683\)}}
\put(-4.8,-8){\makebox(0,0)[lc]{\(3^9\)}}
\put(-5.2,-10){\makebox(0,0)[rc]{\(59\,049\)}}
\put(-4.8,-10){\makebox(0,0)[lc]{\(3^{10}\)}}
\put(-5.2,-12){\makebox(0,0)[rc]{\(177\,147\)}}
\put(-4.8,-12){\makebox(0,0)[lc]{\(3^{11}\)}}
\put(-5.2,-14){\makebox(0,0)[rc]{\(531\,441\)}}
\put(-4.8,-14){\makebox(0,0)[lc]{\(3^{12}\)}}
\put(-5.2,-16){\makebox(0,0)[rc]{\(1\,594\,323\)}}
\put(-4.8,-16){\makebox(0,0)[lc]{\(3^{13}\)}}
\put(-5.2,-18){\makebox(0,0)[rc]{\(4\,782\,969\)}}
\put(-4.8,-18){\makebox(0,0)[lc]{\(3^{14}\)}}
\put(-5,-18){\vector(0,-1){2}}


\put(7,0.5){\makebox(0,0)[cc]{Symmetric topology symbol (ground state):}}
{\color{red}
\put(7,-1){\makebox(0,0)[cc]{
\(\overbrace{\mathrm{E}\binom{1}{\rightarrow}}^{\text{Leaf}}\quad
\overbrace{\mathrm{c}}^{\text{Fork}}\quad
\overbrace{\binom{2}{\leftarrow}\mathrm{E}}^{\text{Leaf}}\)
}}
}
\put(9,-3){\makebox(0,0)[cc]{Transfer kernel types:}}
\put(9,-3.5){\makebox(0,0)[cc]{E.8: \(\varkappa_3=(1231)\), c.21: \(\varkappa_0=(0231)\)}}
\put(9,-4.5){\makebox(0,0)[cc]{Minimal discriminant TKT E.8:}}
\put(9,-5){\makebox(0,0)[cc]{\(-34\,867\)}}

{\color{red}
\put(0.1,0.2){\makebox(0,0)[lb]{\(\langle 8\rangle\)}}
\put(0.8,-1){\makebox(0,0)[cc]{\(s_2=1\)}}
\put(0.1,-1.8){\makebox(0,0)[lb]{\(\langle 54\rangle\)}}
\put(1.1,-2.8){\makebox(0,0)[lb]{\(1^{\text{st}}\) bifurcation}}
}
\put(0.1,-3.8){\makebox(0,0)[lb]{\(1;3\)}}
\put(0.1,-4.2){\makebox(0,0)[lt]{\(\langle 303\rangle\)}}
\put(0.1,-5.8){\makebox(0,0)[lb]{\(1;1\)}}
\put(0.1,-7.8){\makebox(0,0)[lb]{\(1;1\)}}
\put(0.1,-9.8){\makebox(0,0)[lb]{\(1;1\)}}
\put(0.1,-11.8){\makebox(0,0)[lb]{\(1;1\)}}
{\color{red}
\put(0,0){\circle*{0.2}}
\multiput(0,0)(0,-2){1}{\line(0,-1){2}}
\put(0,-2){\circle*{0.2}}
}
\multiput(0,-4)(0,-2){5}{\circle*{0.2}}
\multiput(0,-2)(0,-2){5}{\line(0,-1){2}}
\put(0,-12){\vector(0,-1){2}}
\put(-0.2,-14.2){\makebox(0,0)[rt]{\(\mathcal{T}_\ast^2(\langle 243,8\rangle)\)}}

{\color{blue}
\put(-3,-3){\makebox(0,0)[cc]{metabelianizations}}
\put(-3,-4.2){\makebox(0,0)[ct]{\(1;6\)}}
\put(-3,-4.8){\makebox(0,0)[cc]{\(\langle 306\rangle\)}}
\multiput(-3,-4)(0,-4){1}{\circle*{0.2}}
}
\put(-3,-8.2){\makebox(0,0)[ct]{\(1;6\)}}
\put(-3,-12.2){\makebox(0,0)[ct]{\(1;6\)}}
\multiput(0,-2)(0,-4){3}{\line(-3,-2){3}}
\multiput(-3,-8)(0,-4){2}{\circle*{0.2}}

{\color{blue}
\put(-2,-4.4){\makebox(0,0)[ct]{\(1;2\)}}
\put(-2,-4.8){\makebox(0,0)[cc]{\(\langle 302\rangle\)}}
\multiput(-2,-4)(0,-4){1}{\circle*{0.2}}
}
\put(-2,-8.4){\makebox(0,0)[ct]{\(1;4\)}}
\put(-2,-12.4){\makebox(0,0)[ct]{\(1;4\)}}
\multiput(0,-2)(0,-4){3}{\line(-1,-1){2}}
\multiput(-2,-8)(0,-4){2}{\circle*{0.2}}

{\color{blue}
\put(-1,-4.2){\makebox(0,0)[ct]{\(1;4\)}}
\put(-1,-4.8){\makebox(0,0)[cc]{\(\langle 304\rangle\)}}
\multiput(-1,-4)(0,-4){1}{\circle*{0.2}}
}
\put(-1,-8.2){\makebox(0,0)[ct]{\(1;2\)}}
\put(-1,-12.2){\makebox(0,0)[ct]{\(1;2\)}}
\multiput(0,-2)(0,-4){3}{\line(-1,-2){1}}
\multiput(-1,-8)(0,-4){2}{\circle*{0.2}}



\multiput(1,-6)(1,0){3}{\vector(-2,1){3.9}}

\put(0,-2){\line(1,-1){4}}
{\color{red}
\put(2.7,-4){\makebox(0,0)[cc]{\(s_1=2\)}}
}

\put(4.1,-5.8){\makebox(0,0)[lb]{\(2;3\)}}
\put(4.1,-6.2){\makebox(0,0)[lt]{\(\langle 621\rangle\)}}
\put(4.1,-7.8){\makebox(0,0)[lb]{\(1;1\)}}
\put(5.1,-8.8){\makebox(0,0)[lb]{\(2^{\text{nd}}\) bifurcation}}
\put(4.1,-9.8){\makebox(0,0)[lb]{\(1;1\)}}
\put(4.1,-11.8){\makebox(0,0)[lb]{\(1;1\)}}
\put(4.1,-13.8){\makebox(0,0)[lb]{\(1;1\)}}
\multiput(3.95,-6.05)(0,-2){5}{\framebox(0.1,0.1){}}
\multiput(4,-6)(0,-2){4}{\line(0,-1){2}}
\put(4,-14){\vector(0,-1){2}}
\put(3.8,-16.2){\makebox(0,0)[rt]{\(\mathcal{T}_\ast^3(\langle 729,54\rangle-\#2;3)\)}}

{\color{red}
\put(1,-6.2){\makebox(0,0)[ct]{\(2;6\)}}
\put(1,-6.8){\makebox(0,0)[cc]{\(\langle 624\rangle\)}}
}
\put(1,-10.2){\makebox(0,0)[ct]{\(1;6\)}}
\put(1,-14.2){\makebox(0,0)[ct]{\(1;6\)}}
{\color{red}
\put(0,-2){\line(1,-4){1}}
\multiput(0.9,-6.1)(0,-4){1}{\framebox(0.2,0.2){}}
\multiput(1,-6)(0,-4){1}{\circle*{0.1}}
}
\multiput(4,-8)(0,-4){2}{\line(-3,-2){3}}
\multiput(0.95,-10.05)(0,-4){2}{\framebox(0.1,0.1){}}

{\color{red}
\put(2,-6.4){\makebox(0,0)[ct]{\(2;2\)}}
\put(2,-6.8){\makebox(0,0)[cc]{\(\langle 620\rangle\)}}
}
\put(2,-10.4){\makebox(0,0)[ct]{\(1;4\)}}
\put(2,-14.4){\makebox(0,0)[ct]{\(1;4\)}}
{\color{red}
\put(0,-2){\line(1,-2){2}}
\multiput(1.9,-6.1)(0,-4){1}{\framebox(0.2,0.2){}}
\multiput(2,-6)(0,-4){1}{\circle*{0.1}}
}
\multiput(4,-8)(0,-4){2}{\line(-1,-1){2}}
\multiput(1.95,-10.05)(0,-4){2}{\framebox(0.1,0.1){}}

{\color{red}
\put(3,-6.2){\makebox(0,0)[ct]{\(2;4\)}}
\put(3,-6.8){\makebox(0,0)[cc]{\(\langle 622\rangle\)}}
}
\put(3,-10.2){\makebox(0,0)[ct]{\(1;2\)}}
\put(3,-14.2){\makebox(0,0)[ct]{\(1;2\)}}
{\color{red}
\put(0,-2){\line(3,-4){3}}
\put(2.9,-6.1){\framebox(0.2,0.2){}}
\multiput(3,-6)(0,-4){1}{\circle*{0.1}}
}
\multiput(4,-8)(0,-4){2}{\line(-1,-2){1}}
\multiput(2.95,-10.05)(0,-4){2}{\framebox(0.1,0.1){}}

\put(4,-8){\line(1,-1){4}}

\put(8.1,-11.8){\makebox(0,0)[lb]{\(2;1\)}}
\put(8.1,-13.8){\makebox(0,0)[lb]{\(1;1\)}}
\put(9.1,-14.8){\makebox(0,0)[lb]{\(3^{\text{rd}}\) bifurcation}}
\put(8.1,-15.8){\makebox(0,0)[lb]{\(1;1\)}}
\multiput(7.95,-12.05)(0,-2){3}{\framebox(0.1,0.1){}}
\multiput(8,-12)(0,-2){2}{\line(0,-1){2}}
\put(8,-16){\vector(0,-1){2}}
\put(7.8,-18.2){\makebox(0,0)[rt]{\(\mathcal{T}_\ast^4(\langle 729,54\rangle-\#2;3-\#1;1-\#2;1)\)}}

\put(5,-12.2){\makebox(0,0)[ct]{\(2;6\)}}
\put(5,-16.2){\makebox(0,0)[ct]{\(1;6\)}}
\put(4,-8){\line(1,-4){1}}
\multiput(8,-14)(0,-4){1}{\line(-3,-2){3}}
\multiput(4.9,-12.1)(0,-4){1}{\framebox(0.2,0.2){}}
\multiput(5,-12)(0,-4){1}{\circle*{0.1}}
\multiput(4.95,-16.05)(0,-4){1}{\framebox(0.1,0.1){}}

\put(6,-12.4){\makebox(0,0)[ct]{\(2;4\)}}
\put(6,-16.4){\makebox(0,0)[ct]{\(1;4\)}}
\put(4,-8){\line(1,-2){2}}
\multiput(5.9,-12.1)(0,-4){1}{\framebox(0.2,0.2){}}
\multiput(6,-12)(0,-4){1}{\circle*{0.1}}
\multiput(8,-14)(0,-4){1}{\line(-1,-1){2}}
\multiput(5.95,-16.05)(0,-4){1}{\framebox(0.1,0.1){}}

\put(7,-12.2){\makebox(0,0)[ct]{\(2;2\)}}
\put(7,-16.2){\makebox(0,0)[ct]{\(1;2\)}}
\put(4,-8){\line(3,-4){3}}
\put(6.9,-12.1){\framebox(0.2,0.2){}}
\multiput(7,-12)(0,-4){1}{\circle*{0.1}}
\multiput(8,-14)(0,-4){1}{\line(-1,-2){1}}
\multiput(6.95,-16.05)(0,-4){1}{\framebox(0.1,0.1){}}

\put(8,-14){\line(1,-1){4}}

\put(12.1,-17.8){\makebox(0,0)[lb]{\(2;1\)}}
\multiput(11.95,-18.05)(0,-2){1}{\framebox(0.1,0.1){}}
\put(12,-18){\vector(0,-1){2}}
\put(11.8,-19.9){\makebox(0,0)[rt]{\(\mathcal{T}_\ast^5(\langle 729,54\rangle-\#2;3-\#1;1-\#2;1-\#1;1-\#2;1)\)}}

\put(9,-18.2){\makebox(0,0)[ct]{\(2;6\)}}
\put(8,-14){\line(1,-4){1}}
\multiput(8.9,-18.1)(0,-4){1}{\framebox(0.2,0.2){}}
\multiput(9,-18)(0,-4){1}{\circle*{0.1}}

\put(10,-18.4){\makebox(0,0)[ct]{\(2;4\)}}
\put(8,-14){\line(1,-2){2}}
\multiput(9.9,-18.1)(0,-4){1}{\framebox(0.2,0.2){}}
\multiput(10,-18)(0,-4){1}{\circle*{0.1}}

\put(11,-18.2){\makebox(0,0)[ct]{\(2;2\)}}
\put(8,-14){\line(3,-4){3}}
\put(10.9,-18.1){\framebox(0.2,0.2){}}
\multiput(11,-18)(0,-4){1}{\circle*{0.1}}

\put(-5,-20.9){\makebox(0,0)[cc]{\textbf{TKT:}}}
\put(-3,-21){\makebox(0,0)[cc]{\(\varkappa_1\)}}
\put(-2,-21){\makebox(0,0)[cc]{\(\varkappa_2\)}}
\put(-1,-21){\makebox(0,0)[cc]{\(\varkappa_3\)}}
\put(0,-21){\makebox(0,0)[cc]{\(\varkappa_0\)}}
\put(1,-21){\makebox(0,0)[cc]{\(\varkappa_1\)}}
\put(2,-21){\makebox(0,0)[cc]{\(\varkappa_2\)}}
\put(3,-21){\makebox(0,0)[cc]{\(\varkappa_3\)}}
\put(4,-21){\makebox(0,0)[cc]{\(\varkappa_0\)}}
\put(5,-21){\makebox(0,0)[cc]{\(\varkappa_1\)}}
\put(6,-21){\makebox(0,0)[cc]{\(\varkappa_2\)}}
\put(7,-21){\makebox(0,0)[cc]{\(\varkappa_3\)}}
\put(8,-21){\makebox(0,0)[cc]{\(\varkappa_0\)}}
\put(9,-21){\makebox(0,0)[cc]{\(\varkappa_1\)}}
\put(10,-21){\makebox(0,0)[cc]{\(\varkappa_2\)}}
\put(11,-21){\makebox(0,0)[cc]{\(\varkappa_3\)}}
\put(12,-21){\makebox(0,0)[cc]{\(\varkappa_0\)}}
\put(-5.8,-21.2){\framebox(18.6,0.6){}}

{\color{blue}
\put(-2,-4.25){\oval(3,1.5)}
}
{\color{red}
\put(2,-6.25){\oval(3,1.5)}
}

\end{picture}

\end{figure}

}

\newpage


{\tiny

\begin{figure}[ht]
\caption{Extremal paths to Schur \(\sigma\)-groups, log ord \(11\), on purged tree \(\mathcal{T}_\ast(\langle 243,8\rangle)\)}
\label{fig:E8ES1}


\setlength{\unitlength}{0.8cm}
\begin{picture}(18,22)(-6,-21)


\put(-5,0.5){\makebox(0,0)[cb]{Order}}
\put(-5,0){\line(0,-1){18}}
\multiput(-5.1,0)(0,-2){10}{\line(1,0){0.2}}
\put(-5.2,0){\makebox(0,0)[rc]{\(243\)}}
\put(-4.8,0){\makebox(0,0)[lc]{\(3^5\)}}
\put(-5.2,-2){\makebox(0,0)[rc]{\(729\)}}
\put(-4.8,-2){\makebox(0,0)[lc]{\(3^6\)}}
\put(-5.2,-4){\makebox(0,0)[rc]{\(2\,187\)}}
\put(-4.8,-4){\makebox(0,0)[lc]{\(3^7\)}}
\put(-5.2,-6){\makebox(0,0)[rc]{\(6\,561\)}}
\put(-4.8,-6){\makebox(0,0)[lc]{\(3^8\)}}
{\color{blue}
\put(-5.2,-8){\makebox(0,0)[rc]{\(19\,683\)}}
\put(-4.8,-8){\makebox(0,0)[lc]{\(3^9\)}}
}
\put(-5.2,-10){\makebox(0,0)[rc]{\(59\,049\)}}
\put(-4.8,-10){\makebox(0,0)[lc]{\(3^{10}\)}}
{\color{red}
\put(-5.2,-12){\makebox(0,0)[rc]{\(177\,147\)}}
\put(-4.8,-12){\makebox(0,0)[lc]{\(3^{11}\)}}
}
\put(-5.2,-14){\makebox(0,0)[rc]{\(531\,441\)}}
\put(-4.8,-14){\makebox(0,0)[lc]{\(3^{12}\)}}
\put(-5.2,-16){\makebox(0,0)[rc]{\(1\,594\,323\)}}
\put(-4.8,-16){\makebox(0,0)[lc]{\(3^{13}\)}}
\put(-5.2,-18){\makebox(0,0)[rc]{\(4\,782\,969\)}}
\put(-4.8,-18){\makebox(0,0)[lc]{\(3^{14}\)}}
\put(-5,-18){\vector(0,-1){2}}


\put(7,0.5){\makebox(0,0)[cc]{Symmetric topology symbol (\(1\)st excited state):}}
{\color{red}
\put(7,-1){\makebox(0,0)[cc]{
\(\overbrace{\mathrm{E}\binom{1}{\rightarrow}}^{\text{Leaf}}\quad
\overbrace{\left\lbrace\mathrm{c}\binom{1}{\rightarrow}\right\rbrace^{2}\ }^{\text{Mainline}}\quad
\overbrace{\mathrm{c}}^{\text{Fork}}\quad
\overbrace{\left\lbrace\binom{2}{\leftarrow}\mathrm{c}\binom{1}{\leftarrow}\mathrm{c}\right\rbrace\ }^{\text{Trunk}}\quad
\overbrace{\binom{2}{\leftarrow}\mathrm{E}}^{\text{Leaf}}\)
}}
}
\put(9,-3){\makebox(0,0)[cc]{Transfer kernel types:}}
\put(9,-3.5){\makebox(0,0)[cc]{E.8: \(\varkappa_3=(1231)\), c.21: \(\varkappa_0=(0231)\)}}
\put(9,-4.5){\makebox(0,0)[cc]{Minimal discriminant TKT E.8:}}
\put(9,-5){\makebox(0,0)[cc]{\(-370\,740\)}}

{\color{red}
\put(0.1,0.2){\makebox(0,0)[lb]{\(\langle 8\rangle\)}}
\put(0.8,-1){\makebox(0,0)[cc]{\(s_4=1\)}}
\put(0.1,-1.8){\makebox(0,0)[lb]{\(\langle 54\rangle\)}}
\put(1.1,-2.8){\makebox(0,0)[lb]{\(1^{\text{st}}\) bifurcation}}
}
\put(0.1,-3.8){\makebox(0,0)[lb]{\(1;3\)}}
\put(0.1,-5.8){\makebox(0,0)[lb]{\(1;1\)}}
\put(0.1,-7.8){\makebox(0,0)[lb]{\(1;1\)}}
\put(0.1,-9.8){\makebox(0,0)[lb]{\(1;1\)}}
\put(0.1,-11.8){\makebox(0,0)[lb]{\(1;1\)}}
{\color{red}
\put(0,0){\circle*{0.2}}
\put(0,0){\line(0,-1){2}}
\multiput(0,-2)(0,-2){1}{\circle*{0.2}}
}
\multiput(0,-4)(0,-2){5}{\circle*{0.2}}
\multiput(0,-2)(0,-2){5}{\line(0,-1){2}}
\put(0,-12){\vector(0,-1){2}}
\put(-0.2,-14.2){\makebox(0,0)[rt]{\(\mathcal{T}_\ast^2(\langle 243,8\rangle)\)}}

{\color{blue}
\put(-3,-7){\makebox(0,0)[cc]{metabelianizations}}
\multiput(-3,-8)(0,-4){1}{\circle*{0.2}}
\put(-3,-8.2){\makebox(0,0)[ct]{\(1;6\)}}
}
\put(-3,-4.2){\makebox(0,0)[ct]{\(1;6\)}}
\put(-3,-12.2){\makebox(0,0)[ct]{\(1;6\)}}
\multiput(0,-2)(0,-4){3}{\line(-3,-2){3}}
\multiput(-3,-4)(0,-8){2}{\circle*{0.2}}

{\color{blue}
\multiput(-2,-8)(0,-4){1}{\circle*{0.2}}
\put(-2,-8.4){\makebox(0,0)[ct]{\(1;4\)}}
}
\put(-2,-4.4){\makebox(0,0)[ct]{\(1;2\)}}
\put(-2,-12.4){\makebox(0,0)[ct]{\(1;4\)}}
\multiput(0,-2)(0,-4){3}{\line(-1,-1){2}}
\multiput(-2,-4)(0,-8){2}{\circle*{0.2}}

{\color{blue}
\put(-1,-8){\circle*{0.2}}
\put(-1,-8.2){\makebox(0,0)[ct]{\(1;2\)}}
}
\put(-1,-4.2){\makebox(0,0)[ct]{\(1;4\)}}
\put(-1,-12.2){\makebox(0,0)[ct]{\(1;2\)}}
\multiput(0,-2)(0,-8){2}{\line(-1,-2){1}}
\multiput(-1,-4)(0,-8){2}{\circle*{0.2}}
\put(0,-6){\line(-1,-2){1}}



\multiput(5,-12)(1,0){3}{\vector(-2,1){7.9}}

{\color{red}
\put(0,-2){\line(1,-1){4}}
\put(2.7,-4){\makebox(0,0)[cc]{\(s_3=2\)}}
}

{\color{red}
\put(4.1,-5.8){\makebox(0,0)[lb]{\(2;3\)}}
\put(4.8,-7){\makebox(0,0)[cc]{\(s_2=1\)}}
\put(4.1,-7.8){\makebox(0,0)[lb]{\(1;1\)}}
\put(5.1,-8.8){\makebox(0,0)[lb]{\(2^{\text{nd}}\) bifurcation}}
}
\put(4.1,-9.8){\makebox(0,0)[lb]{\(1;1\)}}
\put(4.1,-11.8){\makebox(0,0)[lb]{\(1;1\)}}
\put(4.1,-13.8){\makebox(0,0)[lb]{\(1;1\)}}
{\color{red}
\multiput(3.95,-6.05)(0,-2){2}{\framebox(0.1,0.1){}}
}
\multiput(3.95,-10.05)(0,-2){3}{\framebox(0.1,0.1){}}
{\color{red}
\put(4,-6){\line(0,-1){2}}
}
\multiput(4,-8)(0,-2){3}{\line(0,-1){2}}
\put(4,-14){\vector(0,-1){2}}
\put(3.8,-16.2){\makebox(0,0)[rt]{\(\mathcal{T}_\ast^3(\langle 729,54\rangle-\#2;3)\)}}

\put(1,-6.2){\makebox(0,0)[ct]{\(2;6\)}}
\put(1,-10.2){\makebox(0,0)[ct]{\(1;6\)}}
\put(1,-14.2){\makebox(0,0)[ct]{\(1;6\)}}
\put(0,-2){\line(1,-4){1}}
\multiput(4,-8)(0,-4){2}{\line(-3,-2){3}}
\multiput(0.9,-6.1)(0,-4){1}{\framebox(0.2,0.2){}}
\multiput(1,-6)(0,-4){1}{\circle*{0.1}}
\multiput(0.95,-10.05)(0,-4){2}{\framebox(0.1,0.1){}}

\put(2,-6.4){\makebox(0,0)[ct]{\(2;2\)}}
\put(2,-10.4){\makebox(0,0)[ct]{\(1;4\)}}
\put(2,-14.4){\makebox(0,0)[ct]{\(1;4\)}}
\put(0,-2){\line(1,-2){2}}
\multiput(4,-8)(0,-4){2}{\line(-1,-1){2}}
\multiput(1.9,-6.1)(0,-4){1}{\framebox(0.2,0.2){}}
\multiput(2,-6)(0,-4){1}{\circle*{0.1}}
\multiput(1.95,-10.05)(0,-4){2}{\framebox(0.1,0.1){}}

\put(3,-6.2){\makebox(0,0)[ct]{\(2;4\)}}
\put(3,-10.2){\makebox(0,0)[ct]{\(1;2\)}}
\put(3,-14.2){\makebox(0,0)[ct]{\(1;2\)}}
\put(0,-2){\line(3,-4){3}}
\multiput(4,-8)(0,-4){2}{\line(-1,-2){1}}
\put(2.9,-6.1){\framebox(0.2,0.2){}}
\multiput(3,-6)(0,-4){1}{\circle*{0.1}}
\multiput(2.95,-10.05)(0,-4){2}{\framebox(0.1,0.1){}}

\put(4,-8){\line(1,-1){4}}
{\color{red}
\put(6.7,-10){\makebox(0,0)[cc]{\(s_1=2\)}}
}

\put(8.1,-11.8){\makebox(0,0)[lb]{\(2;1\)}}
\put(8.1,-13.8){\makebox(0,0)[lb]{\(1;1\)}}
\put(9.1,-14.8){\makebox(0,0)[lb]{\(3^{\text{rd}}\) bifurcation}}
\put(8.1,-15.8){\makebox(0,0)[lb]{\(1;1\)}}
\multiput(7.95,-12.05)(0,-2){3}{\framebox(0.1,0.1){}}
\multiput(8,-12)(0,-2){2}{\line(0,-1){2}}
\put(8,-16){\vector(0,-1){2}}
\put(7.8,-18.2){\makebox(0,0)[rt]{\(\mathcal{T}_\ast^4(\langle 729,54\rangle-\#2;3-\#1;1-\#2;1)\)}}

{\color{red}
\put(5,-12.2){\makebox(0,0)[ct]{\(2;6\)}}
}
\put(5,-16.2){\makebox(0,0)[ct]{\(1;6\)}}
{\color{red}
\put(4,-8){\line(1,-4){1}}
\multiput(4.9,-12.1)(0,-4){1}{\framebox(0.2,0.2){}}
\multiput(5,-12)(0,-4){1}{\circle*{0.1}}
}
\multiput(8,-14)(0,-4){1}{\line(-3,-2){3}}
\multiput(4.95,-16.05)(0,-4){1}{\framebox(0.1,0.1){}}

{\color{red}
\put(6,-12.4){\makebox(0,0)[ct]{\(2;4\)}}
}
\put(6,-16.4){\makebox(0,0)[ct]{\(1;4\)}}
{\color{red}
\put(4,-8){\line(1,-2){2}}
\multiput(5.9,-12.1)(0,-4){1}{\framebox(0.2,0.2){}}
\multiput(6,-12)(0,-4){1}{\circle*{0.1}}
}
\multiput(8,-14)(0,-4){1}{\line(-1,-1){2}}
\multiput(5.95,-16.05)(0,-4){1}{\framebox(0.1,0.1){}}

{\color{red}
\put(7,-12.2){\makebox(0,0)[ct]{\(2;2\)}}
}
\put(7,-16.2){\makebox(0,0)[ct]{\(1;2\)}}
{\color{red}
\put(4,-8){\line(3,-4){3}}
\put(6.9,-12.1){\framebox(0.2,0.2){}}
\multiput(7,-12)(0,-4){1}{\circle*{0.1}}
}
\multiput(8,-14)(0,-4){1}{\line(-1,-2){1}}
\multiput(6.95,-16.05)(0,-4){1}{\framebox(0.1,0.1){}}

\put(8,-14){\line(1,-1){4}}

\put(12.1,-17.8){\makebox(0,0)[lb]{\(2;1\)}}
\multiput(11.95,-18.05)(0,-2){1}{\framebox(0.1,0.1){}}
\put(12,-18){\vector(0,-1){2}}
\put(11.8,-19.9){\makebox(0,0)[rt]{\(\mathcal{T}_\ast^5(\langle 729,54\rangle-\#2;3-\#1;1-\#2;1-\#1;1-\#2;1)\)}}

\put(9,-18.2){\makebox(0,0)[ct]{\(2;6\)}}
\put(8,-14){\line(1,-4){1}}
\multiput(8.9,-18.1)(0,-4){1}{\framebox(0.2,0.2){}}
\multiput(9,-18)(0,-4){1}{\circle*{0.1}}

\put(10,-18.4){\makebox(0,0)[ct]{\(2;4\)}}
\put(8,-14){\line(1,-2){2}}
\multiput(9.9,-18.1)(0,-4){1}{\framebox(0.2,0.2){}}
\multiput(10,-18)(0,-4){1}{\circle*{0.1}}

\put(11,-18.2){\makebox(0,0)[ct]{\(2;2\)}}
\put(8,-14){\line(3,-4){3}}
\put(10.9,-18.1){\framebox(0.2,0.2){}}
\multiput(11,-18)(0,-4){1}{\circle*{0.1}}

\put(-5,-20.9){\makebox(0,0)[cc]{\textbf{TKT:}}}
\put(-3,-21){\makebox(0,0)[cc]{\(\varkappa_1\)}}
\put(-2,-21){\makebox(0,0)[cc]{\(\varkappa_2\)}}
\put(-1,-21){\makebox(0,0)[cc]{\(\varkappa_3\)}}
\put(0,-21){\makebox(0,0)[cc]{\(\varkappa_0\)}}
\put(1,-21){\makebox(0,0)[cc]{\(\varkappa_1\)}}
\put(2,-21){\makebox(0,0)[cc]{\(\varkappa_2\)}}
\put(3,-21){\makebox(0,0)[cc]{\(\varkappa_3\)}}
\put(4,-21){\makebox(0,0)[cc]{\(\varkappa_0\)}}
\put(5,-21){\makebox(0,0)[cc]{\(\varkappa_1\)}}
\put(6,-21){\makebox(0,0)[cc]{\(\varkappa_2\)}}
\put(7,-21){\makebox(0,0)[cc]{\(\varkappa_3\)}}
\put(8,-21){\makebox(0,0)[cc]{\(\varkappa_0\)}}
\put(9,-21){\makebox(0,0)[cc]{\(\varkappa_1\)}}
\put(10,-21){\makebox(0,0)[cc]{\(\varkappa_2\)}}
\put(11,-21){\makebox(0,0)[cc]{\(\varkappa_3\)}}
\put(12,-21){\makebox(0,0)[cc]{\(\varkappa_0\)}}
\put(-5.8,-21.2){\framebox(18.6,0.6){}}

{\color{blue}
\put(-2,-8.25){\oval(3,1.5)}
}
{\color{red}
\put(6,-12.25){\oval(3,1.5)}
}

\end{picture}

\end{figure}
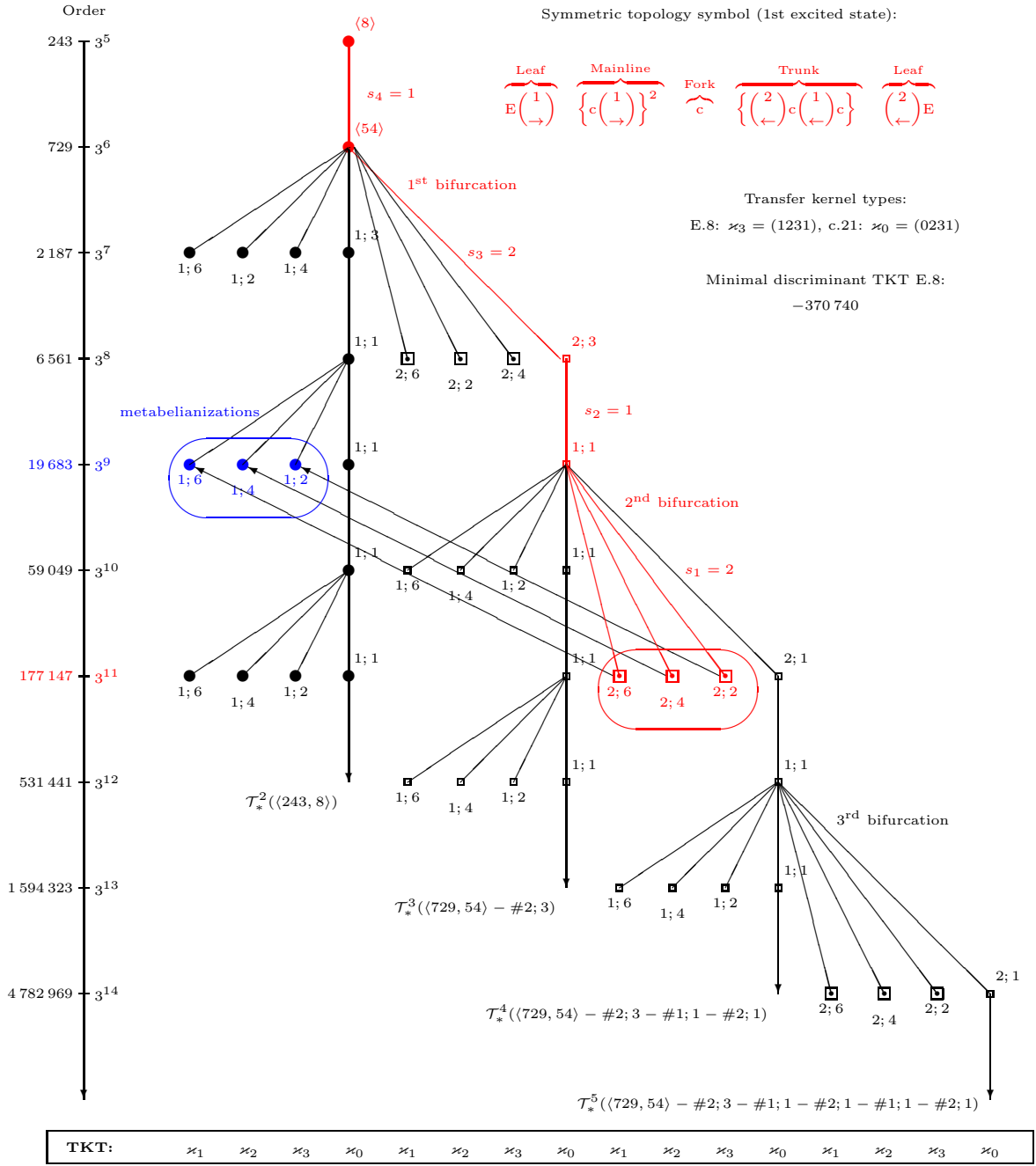

}

\newpage


{\tiny

\begin{figure}[ht]
\caption{Extremal paths to Schur \(\sigma\)-groups, log ord \(14\), on purged tree \(\mathcal{T}_\ast(\langle 243,8\rangle)\)}
\label{fig:E8ES2}

\input{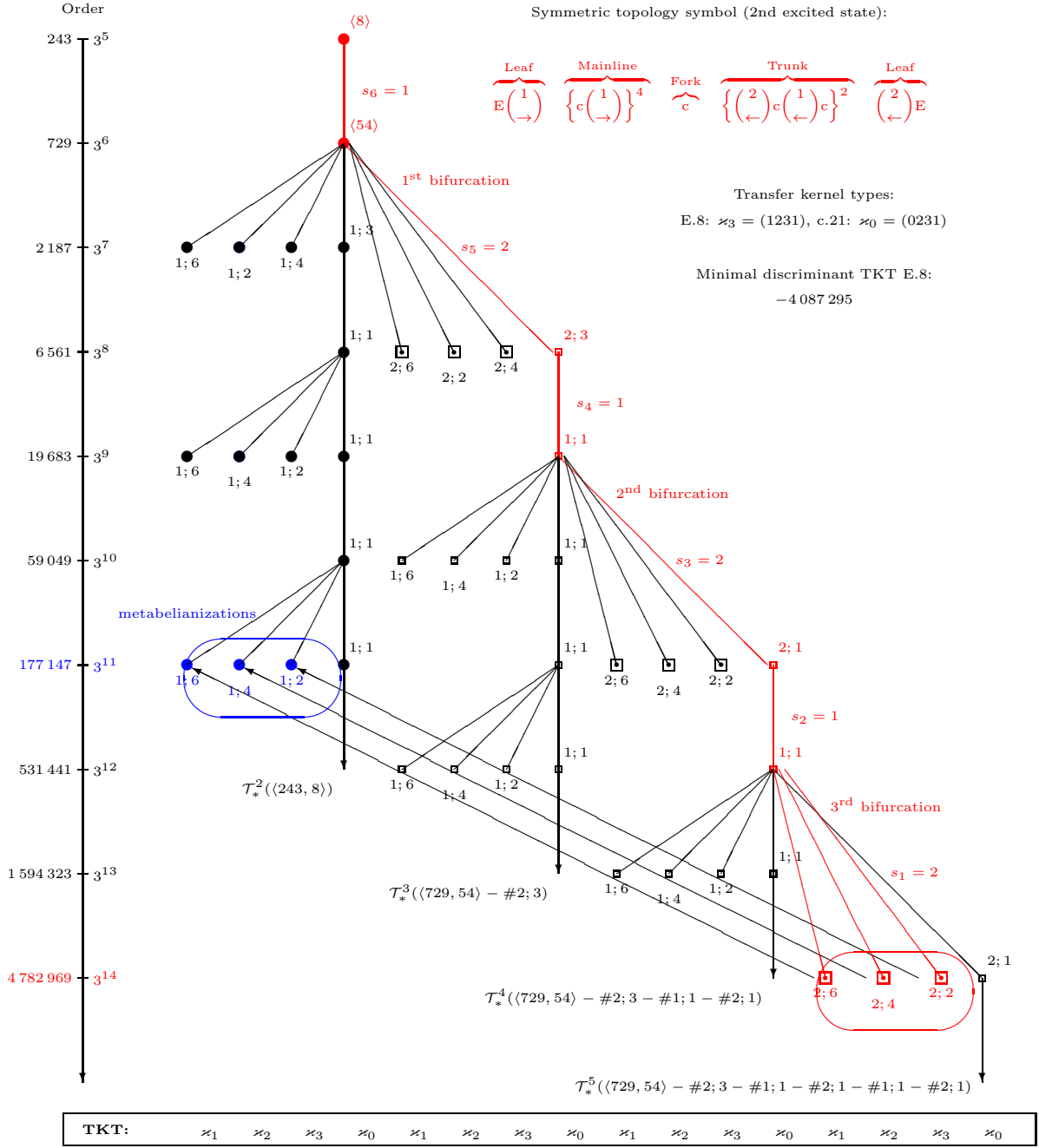}

\end{figure}

}

\newpage


\subsection{Challenges for finding non-metabelian Schur \(\sigma\)-groups}
\label{ss:NonMetabelian}

Why did it take mathematicians so long to discover \textit{non-metabelian} Schur \(\sigma\)-groups?
In \S\
\ref{ss:Metabelian},
we saw that there exist precisely two \textit{metabelian} Schur \(\sigma\)-groups \(V\)
with abelian quotient invariants \(V/V^\prime\simeq (3,3)\).
They are of logarithmic order \(\mathrm{lo}(V)=5\),
nilpotency class \(\mathrm{cl}(V)=3\), coclass \(\mathrm{cc}(V)=2\),
and have transfer kernel types in Section D.
These two groups were known to G. Bagnera
\cite{Bg}
in \(1898\) already,
and they were rediscovered by means of computer aided investigations
in the Ph.D. thesis of J. A. Ascione
\cite{AHL,As1,As2}
in \(1979\) and in the doctoral dissertation of B. Nebelung
\cite{Ne1,Ne2}
in \(1989\).
The reason for these discoveries
was Bagnera's search for all \(p\)-groups of order \(p^5\),
Ascione's investigation of \(3\)-groups with coclass two and order up to \(3^8\),
and Nebelung's complete classification of all metabelian \(3\)-groups with abelianization \((3,3)\).
None of these authors was aware of the balanced presentations.

In \S\S\
\ref{ss:BartholdiBush}
and
\ref{ss:ThreeStage},
we enumerated all \textit{non-metabelian} Schur \(\sigma\)-groups \(V\)
with abelian quotient invariants \(V/V^\prime\simeq (3,3)\)
and transfer kernel type either in Section E
or H.4 under the additional requirement of transfer target type \(\tau(V)=((1^3)^3,21)\).
We saw that they exist for every logarithmic order \(\mathrm{lo}(V)\equiv 2\,(\mathrm{mod}\,3)\)
setting in with log ord \(8\),
and they arise from periodic bifurcations in the (pruned) descendant trees
\(\mathcal{T}_\ast(\langle 729,i\rangle)\)
with either \(i\in\lbrace 49,54\rbrace\) or \(i=45\).
Since the first bifurcation already causes a transition from coclass two to coclass three,
it is clear that non-metabelian Schur \(\sigma\)-groups
were outside of the scope of all three above mentioned authors.

Due to the genesis of Schur \(\sigma\)-groups \(V\) by bifurcations,
they never belong to coclass trees,
they are rather \textit{isolated orphans} without proper parents in the \textit{sporadic part} of coclass forests.
They always have nuclear rank \(\nu(V)=0\) and consequently no descendants.


The following theorem
was known to Boston, Bush, Hajir
\cite{BBH}
in \(2012\) already,
at least partially up to log ord \(11\).

{\color{blue}
\begin{theorem}
\label{thm:SchurSigmaUpToLo14}
The exact counts of
Schur \(\sigma\)-groups \(V\)
of order \(\mathrm{ord}(V)\) a power of \(3\)
with abelian quotient invariants \(V/V^\prime\simeq (3,3)\)
and logarithmic order \(\mathrm{lo}(V)\in\lbrace 5,8,11,14\rbrace\)
are given in the following way,
classified by transfer kernel types (TKT):
\begin{enumerate}
\item
There are \(\mathbf{2}\) groups of order \(3^5=243\): \\
one of TKT \(\mathrm{D}.10\), \(\langle 243,5\rangle\), and one of TKT \(\mathrm{D}.5\), \(\langle 243,7\rangle\) (Fig.
\ref{fig:SporCc2}).
\item
There are \(\mathbf{7}\) groups of order \(3^8=6\,561\): \\
one of TKT \(\mathrm{H}.4\), \(\langle 6561,606\rangle\) (Fig.
\ref{fig:H4S0}), \\
one of TKT \(\mathrm{E}.6\), \(\langle 6561,616\rangle\), two of TKT \(\mathrm{E}.14\), \(\langle 6561,i\rangle\) with \(i\in\lbrace 617,618\rbrace\), \\
one of TKT \(\mathrm{E}.8\), \(\langle 6561,622\rangle\), and two of TKT \(\mathrm{E}.9\), \(\langle 6561,i\rangle\) with \(i\in\lbrace 620,624\rbrace\)
(Fig.
\ref{fig:E8GS}).
\item
There are \(\mathbf{15}\) groups of order \(3^{11}=177\,147\): (Fig.
\ref{fig:H4S1},
\ref{fig:E8ES1},
\ref{fig:G19}) \\
three of TKT \(\mathrm{H}.4\), \(\langle 6561,605\rangle-\#1;2-\#2;2\), \(\langle 6561,i\rangle-\#1;2-\#2;2\) with \(i\in\lbrace 614,615\rbrace\), \\
one of TKT \(\mathrm{E}.6\), \(\langle 6561,613\rangle-\#1;1-\#2;4\), \\
two of TKT \(\mathrm{E}.14\), \(\langle 6561,613\rangle-\#1;1-\#2;j\) with \(j\in\lbrace 5,6\rbrace\), \\
one of TKT \(\mathrm{E}.8\), \(\langle 6561,621\rangle-\#1;1-\#2;2\), \\
two of TKT \(\mathrm{E}.9\), \(\langle 6561,621\rangle-\#1;1-\#2;j\) with \(j\in\lbrace 4,6\rbrace\), \\
two of TKT \(\mathrm{G}.16\), \(\langle 6561,i\rangle-\#1;4-\#2;1\) with \(i\in\lbrace 619,623\rbrace\), \\ 
and four of TKT \(\mathrm{G}.19\), \(\langle 6561,i\rangle-\#1;2-\#2;j\) with \(i\in\lbrace 625,629\rbrace\), \(j\in\lbrace 1,2\rbrace\).
\item
There are \(\mathbf{23}\) groups of order \(3^{14}=4\,782\,969\):
three of TKT \(\mathrm{H}.4\), \\
one of TKT \(\mathrm{E}.6\), two of TKT \(\mathrm{E}.14\),
one of TKT \(\mathrm{E}.8\), two of TKT \(\mathrm{E}.9\), \\
two of TKT \(\mathrm{G}.16\) and twelve of TKT \(\mathrm{G}.19\)
(Fig.
\ref{fig:H4S2},
\ref{fig:E8ES2},
\ref{fig:G19}).
\end{enumerate}
\end{theorem}
}

\begin{proof}
Item (1) was proved in Proposition
\ref{prp:SectionD},
Item (2), and parts of items (3) and (4),
were proved in Propositions
\ref{prp:SporadicSectionH}
and
\ref{prp:SectionE}.
The remainder of items (3) and (4)
will be proved in Proposition
\ref{prp:G19}.
\end{proof}


Currently, the number of Schur \(\sigma\)-groups \(V\)
with order \(\mathrm{ord}(V)=3^{17}\)
and abelian quotient invariants \(V/V^\prime\simeq (3,3)\)
seems to be unknown.
In Theorem
\ref{thm:SporadicCc4CoverSs4},
however,
we shall see that the number of Schur \(\sigma\)-groups \(V\)
with order \(\mathrm{ord}(V)=3^{20}\)
and merely types in Section \(\mathrm{F}\) becomes astronomic,
compared to the numbers in Theorem
\ref{thm:SchurSigmaUpToLo14}.

\newpage


\subsection{Schur \(\sigma\)-groups with TKT in Section F}
\label{ss:SectionF}

Now we come to Schur \(\sigma\)-groups with root paths of \textbf{high complexity}
and orders bigger than \(3^8=6561\), outside of the SmallGroups Library.
In this case,
it is convenient to characterize descendants of a parent \(P\)
by the step size \(s\) and a counter \(c\)
in the form \(P-\#s;c\),
which is used by the ANUPQ package
\cite{GNO},
implemented in GAP
\cite{GAP}
and Magma
\cite{MAGMA}.
In Table
\ref{tbl:F7Simple},
we present a complete root path containing edges with step size \(s=4\),
and in Table
\ref{tbl:F13Complex},
we communicate a \textbf{rudimentary} root path with edges up to step size \(s=8\),
which would continue with an edge of step size \(s=12\)
if this were within the reach of actual computations.
(Of course, it isn't.)
In both tables, \(\mathrm{lo}\) (or log ord) denotes the logarithmic order.


We must illuminate the origin of such complicated root paths
with inspirations which came from algebraic number theory.
Therefore, we start by gathering arithmetical information.

In Table
\ref{tbl:ModerateATI2}
and
\ref{tbl:ExtremeATI2},
we present arithmetical information about iterated index-\(p\) abelianization data (IPADs),
\(\tau^{(2)}{K}=
\lbrack\mathrm{Cl}_3{K};(\mathrm{Cl}_3{E_i};(\mathrm{Cl}_3{E^\prime})_{E^\prime\in\mathrm{Lyr}_1{E_i}})_{1\le i\le 4}\rbrack\),
of second order
for \textit{imaginary} quadratic fields \(K=\mathbb{Q}(\sqrt{d})\)
with \(3\)-class group \(\mathrm{Cl}_3{K}\) of type \(1^2\hat{=}(3,3)\),
TKT \(\mathrm{F}\),
and a second \(3\)-class group \(\mathrm{G}_3^2{K}\) of coclass \(4\),
which occur in the range
\(-5\cdot 10^5<d<0\) of fundamental discriminants.
The negative discriminants were taken from the lower half range of
\cite[Tbl. 3, p. 497]{Ma2012},
but they were separated into the four TKTs in June \(2016\).
The IPAD of first order of such a field has the form
\(\tau^{(1)}{K}=
\lbrack\mathrm{Cl}_3{K};(\mathrm{Cl}_3{E_i})_{1\le i\le 4}\rbrack
=\lbrack 1^2;(32,32,1^3,1^3)\rbrack\),
according to
\cite[Thm. 4.5, pp. 444--445, and Tbl. 6.10, p. 455]{Ma2014}.
We point out that we use \textit{logarithmic type invariants} throughout this article, e.g.,
\(32\hat{=}(27,9)\) and \(1^3\hat{=}(3,3,3)\).
Since the Hilbert \(3\)-class field of the fields \(K\) under investigation
has \(3\)-class group \(\mathrm{Cl}_3{\mathrm{F}_3^1{K}}\simeq 2^31\hat{=}(9,9,9,3)\),
the iterated IPAD of second order of \(K\) has the shape
\(\tau^{(2)}{K}=\lbrack 1^2;(32;2^31,T_1),(32;2^31,T_2),(1^3;2^31,T_3),(1^3;2^31,T_4)\rbrack\),
where the families \(T_1,T_2\), resp. \(T_3,T_4\), consist of \(3\), resp. \(12\), remaining components.
Exceptional entries are printed in \textbf{boldface} font.
The concept of IPADs is due to Boston, Bush and Hajir
\cite{BBH}.
In arithmetical context, IPADs are collections of \textit{abelian type invariants} (ATI)
of \(3\)-class groups \(\mathrm{Cl}_3{F}\) of number fields \(F\).

\bigskip
In Table
\ref{tbl:ModerateATI2}
and
\ref{tbl:ExtremeATI2},
we uniformly have metabelianizations
with coclass \(\mathrm{cc}(\mathrm{G}_3^2{K})=4\).



\renewcommand{\arraystretch}{1.2}

\begin{table}[hb]
\caption{\textbf{Moderate} ATI of second order for imaginary \(K=\mathbb{Q}(\sqrt{d})\)}
\label{tbl:ModerateATI2}
\begin{center}
\begin{tabular}{|r|c|c|c|c|}
\hline
         Type & \multicolumn{4}{|c|}{\(\tau^{(2)}{K}=\lbrack 1^2;(32;2^31,T_1),(32;2^31,T_2),(1^3;2^31,T_3),(1^3;2^31,T_4)\rbrack\)} \\
       \(-d\) & \(T_1\) & \(T_2\) & \(T_3\) & \(T_4\) \\
\hline
 \(\mathrm{F}.7\) & \multicolumn{4}{|c|}{} \\
\hline
 \(225\,299\) & \((31^3)^3\)           & \((31^3)^3\)           & \((2^21)^3,(21^2)^9\)                     & \((2^21)^3,(21^2)^9\)                     \\
 \(343\,380\) & \((31^3)^3\)           & \((31^3)^3\)           & \((2^21)^3,(21^2)^9\)                     & \((2^21)^3,(21^2)^9\)                     \\
 \(423\,476\) & \((31^3)^3\)           & \((31^3)^3\)           & \((2^21)^3,(21^2)^9\)                     & \((2^21)^3,(21^2)^9\)                     \\
 \(486\,264\) & \((31^3)^3\)           & \((31^3)^3\)           & \((2^21)^3,(21^2)^9\)                     & \((2^21)^3,(21^2)^9\)                     \\
\hline
 \(\mathrm{F}.11\) & \multicolumn{4}{|c|}{} \\
\hline
  \(27\,156\) & \((41^3)^3\)           & \((31^3)^3\)           & \((2^21)^3,(21^2)^9\)                     & \((2^21)^3,(21^2)^9\)                     \\
 \(241\,160\) & \((41^3)^3\)           & \((31^3)^3\)           & \((2^21)^3,(21^2)^9\)                     & \((2^21)^3,(21^2)^9\)                     \\
 \(477\,192\) & \((41^3)^3\)           & \((31^3)^3\)           & \((2^21)^3,(21^2)^9\)                     & \((2^21)^3,(21^2)^9\)                     \\
 \(484\,804\) & \((41^3)^3\)           & \((31^3)^3\)           & \((2^21)^3,(21^2)^9\)                     & \((2^21)^3,(21^2)^9\)                     \\
\hline
 \(\mathrm{F}.12\) & \multicolumn{4}{|c|}{} \\
\hline
 \(291\,220\) & \((31^3)^3\)           & \((31^3)^3\)           & \((2^21)^3,(21^2)^9\)                     & \((2^21)^3,(21^2)^9\)                     \\
\hline
 \(\mathrm{F}.13\) & \multicolumn{4}{|c|}{} \\
\hline
 \(167\,064\) & \((41^3)^3\)           & \((31^3)^3\)           & \((2^21)^3,(21^2)^9\)                     & \((2^21)^3,(21^2)^9\)                     \\
 \(296\,407\) & \((41^3)^3\)           & \((31^3)^3\)           & \((2^21)^3,(21^2)^9\)                     & \((2^21)^3,(21^2)^9\)                     \\
 \(317\,747\) & \((41^3)^3\)           & \((31^3)^3\)           & \((2^21)^3,(21^2)^9\)                     & \((2^21)^3,(21^2)^9\)                     \\
 \(401\,603\) & \((41^3)^3\)           & \((31^3)^3\)           & \((2^21)^3,(21^2)^9\)                     & \((2^21)^3,(21^2)^9\)                     \\
\hline
\end{tabular}
\end{center}
\end{table}


\noindent
The following theorem was proved in
\cite{MaNm}
and shows that Schur \(\sigma\)-groups in the descendant tree of \(\langle 243,3\rangle\),
resp. \(P_7:=\langle 2187,64\rangle\),
set in extensively at logarithmic order \(20\).

{\color{blue}
\begin{theorem}
\label{thm:SporadicCc4CoverSs4}
Let \(M:=P_7-\#2;m\) be a sporadic metabelian \(3\)-group \(M\)
with type in Section \(\mathrm{F}\) and coclass \(\mathrm{cc}(M)=4\).
The following counters concern \(\mathbf{1359}\) pairwise non-isomorphic Schur \(\sigma\)-groups \(S\)
of logarithmic order \(\mathrm{lo}(S)=20\) and nilpotency class \(\mathrm{cl}(S)=9\)
such that \(S/S^{\prime\prime}\simeq M\).
\begin{enumerate}
\item
For type \(\mathrm{F}.7\), there exist \(171\), in more detail, 
\(81\), resp. \(45\), \(45\), Schur \(\sigma\)-groups \(S\) satisfying
\begin{equation}
\label{eqn:SchurSigmaLo20F7IPAD}
\tau^{(2)}(S)=\lbrack 1^2;(32;2^31,(31^3)^3)^2,(1^3;2^31,(2^21)^3,(21^2)^9)^2\rbrack
\end{equation}
with metabelianization \(m=55\), resp. \(56\), \(58\).
They all have \(\#\mathrm{Aut}(S)=2\cdot 3^{25}\).

\item
For type \(\mathrm{F}.11\), there exist \(108+324\), in more detail,
\begin{enumerate}
\item
\(54\), resp. \(54\), Schur \(\sigma\)-groups \(S\) satisfying
\begin{equation}
\label{eqn:SchurSigmaLo20F11IPAD}
\tau^{(2)}(S)=\lbrack 1^2;(32;2^31,(41^3)^3),(32;2^31,(31^3)^3),(1^3;2^31,(2^21)^3,(21^2)^9)^2\rbrack
\end{equation}
with \(\#\mathrm{Aut}(S)=2\cdot 3^{25}\)
and metabelianization \(m=36\), resp. \(38\);
\item
\(162\), resp. \(162\), Schur \(\sigma\)-groups \(S\) satisfying Formula
\eqref{eqn:SchurSigmaLo20F11IPAD}
with \(\#\mathrm{Aut}(S)=2\cdot 3^{26}\)
and metabelianization \(m=36\), resp. \(38\).
\end{enumerate}

\item
For type \(\mathrm{F}.12\), there exist \(216+162\), in more detail,
\begin{enumerate}
\item
\(54\), resp. \(54\), \(54\), \(54\), Schur \(\sigma\)-groups \(S\) satisfying
\begin{equation}
\label{eqn:SchurSigmaLo20F12IPAD}
\tau^{(2)}(S)=\lbrack 1^2;(32;2^31,(n1^3)^3),(32;2^31,(31^3)^3),(1^3;2^31,(2^21)^3,(21^2)^9)^2\rbrack,
\end{equation}
with \(n=3\)
and metabelianization \(m=43\), resp. \(46\), \(51\), \(53\);
\item
\(54\), resp. \(54\), \(27\), \(27\), Schur \(\sigma\)-groups \(S\) satisfying Formula
\eqref{eqn:SchurSigmaLo20F12IPAD}
with \(n=4\)
and metabelianization \(m=43\), resp. \(46\), \(51\), \(53\).
\end{enumerate}
They all have \(\#\mathrm{Aut}(S)=2\cdot 3^{25}\).

\item
For type \(\mathrm{F}.13\), there exist \(216+162\), in more detail,
\begin{enumerate}
\item
\(54\), resp. \(54\), \(54\), \(54\), Schur \(\sigma\)-groups \(S\) satisfying
\begin{equation}
\label{eqn:SchurSigmaLo20F13IPAD}
\tau^{(2)}(S)=\lbrack 1^2;(32;2^31,(n1^3)^3),(32;2^31,(31^3)^3),(1^3;2^31,(2^21)^3,(21^2)^9)^2\rbrack,
\end{equation}
with \(n=3\)
and metabelianization \(m=41\), resp. \(47\), \(50\), \(52\);
\item
\(54\), resp. \(54\), \(27\), \(27\), Schur \(\sigma\)-groups \(S\) satisfying Formula
\eqref{eqn:SchurSigmaLo20F13IPAD}
with \(n=4\)
and metabelianization \(m=41\), resp. \(47\), \(50\), \(52\).
\end{enumerate}
They all have \(\#\mathrm{Aut}(S)=2\cdot 3^{25}\).
\end{enumerate}
\end{theorem}
}

\begin{proof}
This impressive result goes back to some sparkling emails
which we received from Professor Mike F. Newman
in January and February 2013.
These communications, though being cryptic without explicit identifiers,
unambiguously illuminated the way to the rigorous proof.
So Theorem
\ref{thm:SporadicCc4CoverSs4}
without doubt is a joint achievement by Prof. Newman and ourselves.
The detailed justification constitutes the dominating part of
\cite{MaNm}.
\end{proof}


{\color{red}
\begin{conjecture}
\label{cnj:SchurSigmaLo20}
(Tower ground state)
The imaginary quadratic fields \(K=\mathbb{Q}(\sqrt{d})\) with fundamental discriminants
\(d\in\lbrace -225\,299, -343\,380, -423\,476, -486\,264\rbrace\) of type \(\mathrm{F}.7\),
resp. \(d\in\lbrace -27\,156, -241\,160, -477\,192, -484\,804\rbrace\) of type \(\mathrm{F}.11\),
resp. \(d= -291\,220\) of type \(\mathrm{F}.12\),
resp. \(d\in\lbrace -167\,064, -296\,407, -317\,747, -401\,603\rbrace\) of type \(\mathrm{F}.13\),
have \(3\)-class field towers of \textbf{exact length} \(\mathbf{\ell_3(K)=3}\) with a suitable Schur \(\sigma\)-group in Theorem
\ref{thm:SporadicCc4CoverSs4}.

For all types, \(\mathrm{F}.7\), \(\mathrm{F}.11\), \(\mathrm{F}.12\), \(\mathrm{F}.13\),
the tower group \(S=\mathrm{G}_3^\infty{K}\) has
\(\mathrm{lo}(S)=20\), \(\mathrm{cl}(S)=9\), \(\mathrm{cc}(S)=11\), \(\zeta_1{S}=(9,9)\) or \((9,3,3)\),
\(\gamma_2^2{S}=(27,27,9,3,3,3)\) or \((27,9,9,9,3,3)\), 
and usually \(\#\mathrm{Aut}(S)=2\cdot 3^{25}\), rarely \(2\cdot 3^{26}\).
\end{conjecture}
}


The extremal root path to the \(3\)-class field tower group
\(G=\mathrm{Gal}(\mathrm{F}_3^\infty(K)/K)\simeq S\)
of \(K=\mathbb{Q}(\sqrt{-225\,299})\)
is described in Table
\ref{tbl:F7Simple}
and drawn in Figure
\ref{fig:F7Lo20}.
As opposed to the figures concerning TKT H.4 and Section E,
we do not know to which coclass trees the capable vertices on the path belong.
We assume \(M=S/S^{\prime\prime}\simeq P_7-\#2;55\), \(S\simeq P_7-\#4;196-\#2;31-\#4;1-\#1;2-\#2;1\).


\renewcommand{\arraystretch}{1.2}

\begin{table}[ht]
\caption{\textbf{Extreme} ATI of second order for imaginary \(K=\mathbb{Q}(\sqrt{d})\)}
\label{tbl:ExtremeATI2}
\begin{center}
\begin{tabular}{|r|c|c|c|c|}
\hline
         Type & \multicolumn{4}{|c|}{\(\tau^{(2)}{K}=\lbrack 1^2;(32;2^31,T_1),(32;2^31,T_2),(1^3;2^31,T_3),(1^3;2^31,T_4)\rbrack\)} \\
       \(-d\) & \(T_1\) & \(T_2\) & \(T_3\) & \(T_4\) \\
\hline
 \(\mathrm{F}.7\) & \multicolumn{4}{|c|}{} \\
\hline
 \(124\,363\) & \(\mathbf{(421^2)^3}\) & \(\mathbf{(321^2)^3}\) & \(\mathbf{(2^31)^3,(1^6)^3,(2^21^2)^6}\)  & \(\mathbf{(2^21^3)^3,(21^4)^3,(2^21^2)^6}\) \\
 \(260\,515\) & \(\mathbf{(3^221)^3}\) & \(\mathbf{(3^221)^3}\) & \(\mathbf{(32^21)^3,(2^31)^3,(1^6)^3,(2^21^2)^3}\) & \(\mathbf{(32^21)^3,(1^6)^3,(21^4)^3,(2^21^2)^3}\) \\
\hline
 \(\mathrm{F}.12\) & \multicolumn{4}{|c|}{} \\
\hline
 \(160\,403\) & \(\mathbf{(321^2)^3}\) & \(\mathbf{(321^2)^3}\) & \(\mathbf{(2^4)^3,(21^4)^6,(2^21^2)^3}\)  & \(\mathbf{(2^21^2)^3,(1^5)^6,(21^3)^3}\)  \\
\hline
 \(\mathrm{F}.13\) & \multicolumn{4}{|c|}{} \\
\hline
 \(224\,580\) & \(\mathbf{(321^2)^3}\) & \(\mathbf{(321^2)^3}\) & \(\mathbf{(21^4)^3,(1^5)^6,(2^21^2)^3}\)  & \(\mathbf{(1^5)^3,(2^21^2)^3,(21^3)^6}\)  \\
\hline
\end{tabular}
\end{center}
\end{table}

\newpage


\subsection{SENSATION: First \(3\)-class field towers with at least four stages}
\label{ss:FourStageTower}
(\(\mathrm{dl}(G)\ge 4\))


{\color{blue}
\begin{theorem}
\label{thm:FourStageTowerF13}
An imaginary quadratic field \(K\)
with elementary \(3\)-class group \(\mathrm{Cl}_3(K)\simeq (3,3)\) of rank two,
\(3\)-capitulation type \(\mathrm{F}.13\), \(\varkappa=(3143)\),
and abelian type invariants of second order, \(\tau^{(2)}{K}\)
\begin{equation}
\label{eqn:ATI224580}
\lbrack 1^2;(32;2^31,(321^2)^3)^2,(1^3;2^31,(21^4)^3,(1^5)^6,(2^21^2)^3),(1^3;2^31,(1^5)^3,(2^21^2)^3,(21^3)^6)\rbrack
\end{equation}
possesses a \(3\)-class field tower with length \(\mathbf{\ell_3(K)\ge 4}\).

The field \(K=\mathbb{Q}(\sqrt{d})\)
with fundamental discriminant \(d=-224\,580\)
(Tbl.
\ref{tbl:ExtremeATI2})
is such an example.
\end{theorem}
}

\begin{proof}
Let \(P_7:=\langle 2187,64\rangle\) denote the common fork vertex
for all \(\sigma\)-groups with transfer kernel types in Section F.
We use the \(p\)-group generation algorithm by Newman
\cite{Nm}
and O'Brien
\cite{Ob}
(see also
\cite{HEO})
in order to construct the relevant part of the descendant tree of the root \(P_7\).

Since the capitulation of the field \(K\) is of type F.13, \(\varkappa=(3143)\),
(containing a transposition but no fixed point)
with three different kernels, \(1\), \(3\) and \(4\),
and the logarithmic abelian type invariants of
the four unramified cyclic cubic extensions \(E_t/K\), \(1\le t\le 4\),
are given by \(\tau=(32,32,1^3,1^3)\),
there are precisely four possibilities for the metabelianization \(M=G/G^{\prime\prime}\),
namely \(M\simeq M_{9,m}:=P_7-\#2;m\) with \(m\in\lbrace 41,47,50,52\rbrace\),
and four corresponding possibilities for ancestors of the tower group \(G\) itself,
\(V_{11,\ell}:=P_7-\#4;\ell\) with \(\ell\in\lbrace 9,15,18,20\rbrace\),
in the same order, e.g. \(M_{9,41}\simeq V_{11,9}/V_{11,9}^{\prime\prime}\).
(So the common fork of the root paths of \(M_{9,m}\) and \(V_{11,\ell}\) is always \(P_7\).)
These are the unique four groups
for which all abelian quotients of subgroups of index \(9\) possess \(3\)-rank \(4\):
\begin{equation}
\label{eqn:AQIofLO11}
\lbrack 1^2;(32;\mathbf{2^31,(31^3)^3})^2,(1^3;\mathbf{2^31,(21^3)^3,(1^4)^9})^2\rbrack.
\end{equation}
This will be required in order to arrive at the intended AQI in Formula
\eqref{eqn:ATI224580}.

Each of the four \(V_{11,\ell}\) has derived length \(3\) and possesses
a unique \(\sigma\)-descendant among \(41\) descendants \(V_{11,\ell}-\#2;k\) of step size \(s=2\),
namely \(V_{13,\ell}:=V_{11,\ell}-\#2;10\) for \(\ell\in\lbrace 9,15\rbrace\) (i.e. \(k=10\)),
and \(V_{13,\ell}:=V_{11,\ell}-\#2;2\) for \(\ell\in\lbrace 18,20\rbrace\) (i.e. \(k=2\)).

The further search yields the following results
(see Table
\ref{tbl:F13Complex}
for the case \(\ell=20\)):
Each of the four \(V_{13,\ell}\) has \(729\) descendants of step size \(s=4\).
For \(\ell\in\lbrace 9,15\rbrace\), resp. \(\ell\in\lbrace 18,20\rbrace\),
we find \(81\), resp. \(243\),
pairs \((j,i)\) with \(1\le j\le 729\) and \(1\le i\le 3281\)
such that \(V_{22,\ell,j}:=V_{17,\ell,j}-\#5;i\) is
the unique \(\sigma\)-descendant among \(3281\) descendants of step size \(s=5\)
of \(V_{17,\ell,j}:=V_{13,\ell}-\#4;j\).
Among these groups,
a veritable abundance of slight variations of
the abelian quotient invariants of second order arises,
covering all cases in Table
\ref{tbl:ExtremeATI2}.

Each of the groups \(V_{22,\ell,j}\) has \(19\,683\) descendants of step size \(s=8\).
The main difficulty in this proof
was the construction of the \(265\,721\) descendants of step size \(s=7\)
for each of the groups \(V_{30,\ell,j,h}:=V_{22,\ell,j}-\#8;h\) with \(1\le h\le 19\,683\),
which are still of soluble length \(\mathrm{dl}(V_{30,\ell,j,h})=3\).
The storage for such a batch of \(265\,721\) so-called compact presentations
in the computational algebra system Magma
\cite{BCP,BCFS,MAGMA}
occupied \(120\,\)GB of RAM,
which blew up to \(170\,\)GB when additional group theoretic operations were performed on these descendants.

Summarizing the results of this proof,
it turned out that all the groups of order \(3^{37}\),
constructed in the last computationally possible step,
were of the common shape \(V_{37,\ell,j,h,g}:=\)
\begin{equation}
\label{eqn:FourStageTower}
P_7-\#4;\ell-\#2;k-\#4;j-\#5;i-\#8;h-\#7;g, \text{ with } \mathrm{dl}(V_{37})=4,\ \nu(V_{37})=12,\ \mu(V_{37})=14.
\end{equation}
(Recall that \(k\) is determined uniquely by \(\ell\)
and \(i\) is determined uniquely by \(j\).)
This implies that the desired Schur \(\sigma\)-group \(S\simeq G=\mathrm{Gal}(\mathrm{F}_3^\infty(K)/K)\)
is either infinite or a descendant of one of the groups \(V_{37}\),
and thus must have derived length \(\mathrm{dl}(S)\ge 4\).
So the \(3\)-class field tower of \(K\) has certainly at least four stages.
\end{proof}


\begin{remark}
\label{rmk:FourStageTower}
In September \(2017\), Professor Mike F. Newman
sent us a somewhat cryptic message
in which he indicated a rough sketch of the preceding proof,
up to logarithmic order \(30\) and nilpotency class \(9\).
In particular, he pointed out,
that if the next step with nilpotency class \(10\)
yielded groups of soluble length \(4\),
then \(\mathbb{Q}(\sqrt{-224\,580})\)
would certainly have a \(3\)-class field tower with at least four stages.

The details of Prof. Newman's guidelines
were carried out by ourselves in March \(2020\)
and culminated in a final verification of his prophetic presentiment
on \(06\) April \(2020\).
The conduction of the preceding proof posed
extraordinary demands on our computational resources.
We had to employ our most powerful workstation
with a CPU consisting of two hyperthreaded Intel Xeon eight-core processors
(i.e. a total of 32 independent threads)
and \(256\) Giga Byte RAM storage.

These difficulties could be avoided
when individual descendants of an assigned finite \(p\)-group
could be constructed
without computing the entire set of all children having a fixed step size.
In September \(2017\), Prof. Newman mentioned that he knows such a method,
but unfortunately he did not communicate this technique to us.
\end{remark}


\begin{hypothesis}
\label{hth:FourStageTower}
The \textbf{rudimentary} extremal root path to the \(3\)-class field tower group
\(G=\mathrm{Gal}(\mathrm{F}_3^\infty(K)/K)\simeq S\)
of \(K=\mathbb{Q}(\sqrt{-224\,580})\)
is described in Table
\ref{tbl:F13Complex}
and drawn in Figure
\ref{fig:F13Lo37}.
We assume \(M=S/S^{\prime\prime}\simeq P_7-\#2;41\) and
\(\pi^{c-9}S\simeq P_7-\#4;9-\#2;10-\#4;144-\#5;2516-\#8;1\),
where \(c=\mathrm{cl}(S)\) denotes the unknown nilpotency class of \(S\).
The illustration terminates at the last vertex with soluble length \(\mathrm{dl}=3\).
For the following two vertices
\(\pi^{c-10}S=\pi^{c-9}S-\#7;1\)
and (hypothetically, because outside of current computational reach)
\(\pi^{c-11}S=\pi^{c-10}S-\#12;1\),
the desired soluble length \(\mathrm{dl}=4\) sets in.
Anyway, the proven fact that \(\mathrm{dl}(\pi^{c-10}S)=4\)
rigorously justifies \(\mathrm{dl}(S)\ge 4\) and thus the claimed four-stage tower, QED.
\end{hypothesis}


{\color{blue}
\begin{theorem}
\label{thm:FourStageTowersF7F12}
Each of the imaginary quadratic fields \(K\) in Table
\ref{tbl:ExtremeATI2}
possesses a \(3\)-class field tower with \textbf{at least four stages}.
\end{theorem}
}

\begin{proof}
For transfer kernel type F.7
there are three possibilities for the metabelianization \(M=G/G^{\prime\prime}\),
namely \(M\simeq M_{9,m}:=P_7-\#2;m\) with \(m\in\lbrace 55,56,58\rbrace\),
and three possibilities for the corresponding sibling on the path to \(G\),
\(V_{11,\ell}:=P_7-\#4;\ell\) with \(\ell\in\lbrace 23,24,26\rbrace\).

For transfer kernel type F.12,
there are four possibilities for the metabelianization \(M=G/G^{\prime\prime}\),
namely \(M\simeq M_{9,m}:=P_7-\#2;m\) with \(m\in\lbrace 43,46,51,53\rbrace\),
and four possibilities for the corresponding sibling on the path to \(G\),
\(V_{11,\ell}:=P_7-\#4;\ell\) with \(\ell\in\lbrace 11,14,19,21\rbrace\).

And, hypothetically, for transfer kernel type F.11,
although no suitable field has been detected up to now,
there are only two possibilities for the metabelianization \(M=G/G^{\prime\prime}\),
namely \(M\simeq M_{9,m}:=P_7-\#2;m\) with \(m\in\lbrace 36,38\rbrace\),
and two possibilities for the corresponding sibling on the path to \(G\),
\(V_{11,\ell}:=P_7-\#4;\ell\) with \(\ell\in\lbrace 4,6\rbrace\).
\end{proof}

\newpage


\renewcommand{\arraystretch}{1.1}

\begin{table}[hb]
\caption{Root path of \(G\), log ord 20, for the simplest case of transfer kernel type F.7}
\label{tbl:F7Simple}
\begin{center}
\begin{tabular}{|c|c|c||c|c|c|}
\hline
 Ancestor     & Vertex                     & lo & \((\nu,\mu)\)                   & \((N_s/C_s)_{1\le s\le\nu}\)                            & TKT  \\
\hline
 \(\pi^7(G)\) & \(\langle 27,3\rangle\)    &  3 & \(({\color{red}\mathbf{2,4}})\) & \((4/1,{\color{red}\mathbf{7/5}})\)                     & a.1  \\
 \(\pi^6(G)\) & \(\langle 243,3\rangle\)   &  5 & \(({\color{red}\mathbf{2,4}})\) & \((10/6,{\color{red}\mathbf{15/15}})\)                  & b.10 \\
 \(\pi^5(G)\) & \(\langle 2187,64\rangle\) &  7 & \(({\color{red}\mathbf{4,6}})\) & \((33/2,453/84,918/540,{\color{red}\mathbf{198/198}})\) & b.10 \\
 \(\pi^4(G)\) & \(-\#4;196\)               & 11 & \(({\color{red}\mathbf{2,4}})\) & \((20/20,{\color{red}\mathbf{41/41}})\)                 & F.7  \\
 \(\pi^3(G)\) & \(-\#2;31\)                & 13 & \(({\color{red}\mathbf{4,6}})\) & \((44/0,204/3,180/24,{\color{red}\mathbf{27/27}})\)     & F.7  \\
 \(\pi^2(G)\) & \(-\#4;1\)                 & 17 & \(({\color{red}\mathbf{1,3}})\) & \(({\color{red}\mathbf{5/5}})\)                         & F.7  \\
 \(\pi(G)\)   & \(-\#1;2\)                 & 18 & \(({\color{red}\mathbf{2,4}})\) & \((4/0,{\color{red}\mathbf{1/0}})\)                     & F.7  \\
 \(G\)        & \(-\#2;1\)                 & 20 & \(({\color{red}\mathbf{0,2}})\) &                                                         & F.7  \\
\hline
\end{tabular}
\end{center}
\end{table}


\renewcommand{\arraystretch}{1.1}

\begin{table}[hb]
\caption{Incomplete root path of \(G\) for an extreme case of transfer kernel type F.13}
\label{tbl:F13Complex}
\begin{center}
\begin{tabular}{|c|c|c||c|c|c|c|c|}
\hline
 Ancestor          & Vertex                     & lo & cl & dl & \((\nu,\mu)\)                     & \((N_s/C_s)_{1\le s\le\nu}\)                                   & TKT  \\
\hline
 \(\pi^{c-2}(G)\)  & \(\langle 27,3\rangle\)    &  3 &  2 &  2 & \(({\color{red}\mathbf{2,4}})\)   & \((4/1,{\color{red}\mathbf{7/5}})\)                            & a.1  \\
 \(\pi^{c-3}(G)\)  & \(\langle 243,3\rangle\)   &  5 &  3 &  2 & \(({\color{red}\mathbf{2,4}})\)   & \((10/6,{\color{red}\mathbf{15/15}})\)                         & b.10 \\
 \(\pi^{c-4}(G)\)  & \(\langle 2187,64\rangle\) &  7 &  4 &  2 & \(({\color{red}\mathbf{4,6}})\)   & \((33/2,453/84,918/540,{\color{red}\mathbf{198/198}})\)        & b.10 \\
 \(\pi^{c-5}(G)\)  & \(-\#4;9\)                 & 11 &  5 &  3 & \(({\color{red}\mathbf{2,4}})\)   & \((20/20,{\color{red}\mathbf{41/41}})\)                        & F.13 \\
 \(\pi^{c-6}(G)\)  & \(-\#2;10\)                & 13 &  6 &  3 & \(({\color{red}\mathbf{4,6}})\)   & \((108/54,1674/1674,3564/3564,{\color{red}\mathbf{729/729}})\) & F.13 \\
 \(\pi^{c-7}(G)\)  & \(-\#4;144\)               & 17 &  7 &  3 & \(({\color{red}\mathbf{5,7}})\)   & \((227/173,\ldots,{\color{red}\mathbf{3281/3281}})\)           & F.13 \\
 \(\pi^{c-8}(G)\)  & \(-\#5;2516\)              & 22 &  8 &  3 & \(({\color{red}\mathbf{8,10}})\)  & \((3444/54,\ldots,{\color{red}\mathbf{19683/19683}})\)         & F.13 \\
 \(\pi^{c-9}(G)\)  & \(-\#8;1\)                 & 30 &  9 &  3 & \(({\color{red}\mathbf{7,9}})\)   & \((\ldots,{\color{red}\mathbf{265721/265721}})\)               & F.13 \\
 \(\pi^{c-10}(G)\) & \(-\#7;1\)                 & 37 & 10 &  4 & \(({\color{red}\mathbf{12,14}})\) & \((\ldots,{\color{red}\mathbf{?/?}})\)                         & F.13 \\
\hline
\end{tabular}
\end{center}
\end{table}

\newpage


{\tiny

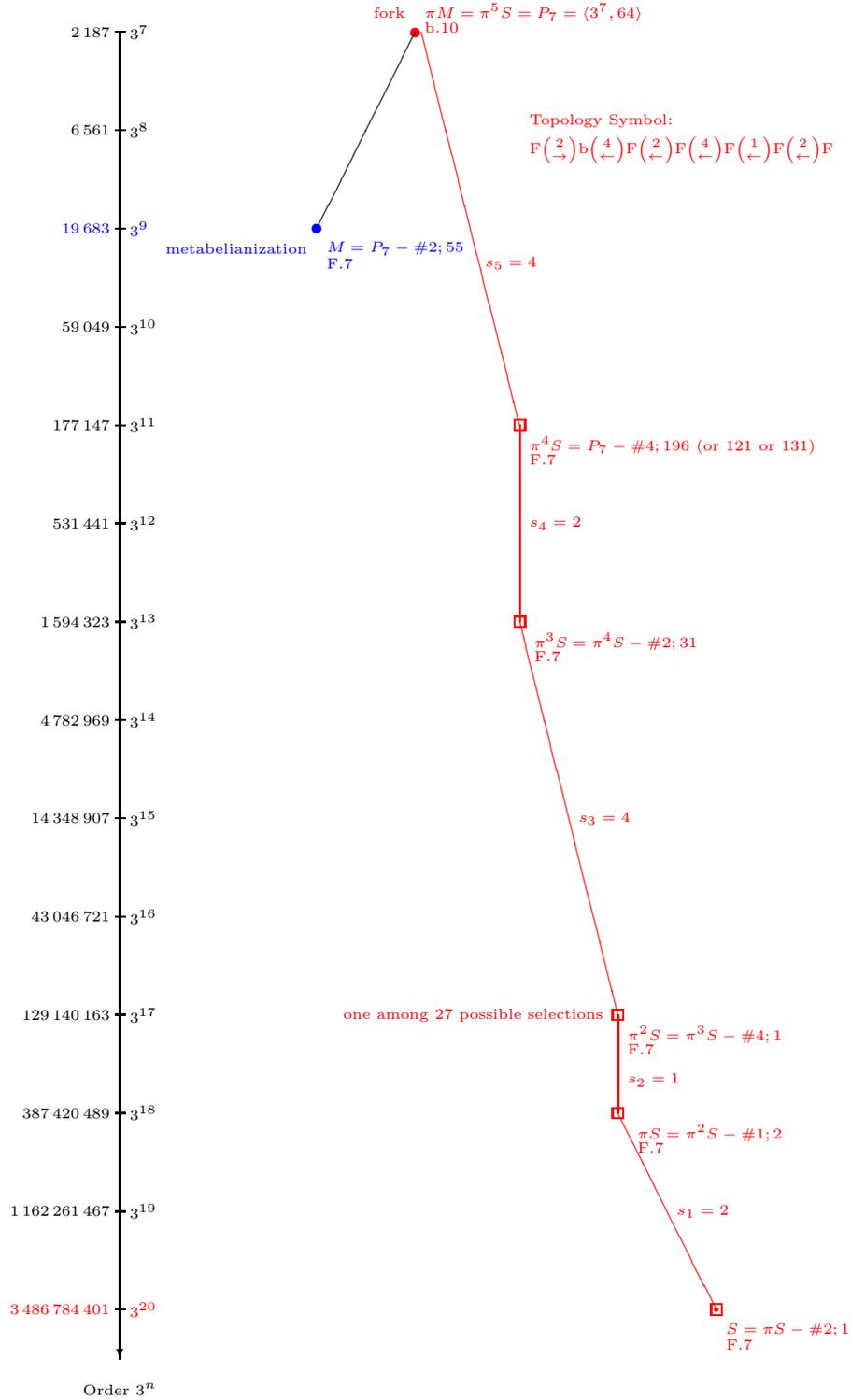
\begin{figure}[hb]
\caption{Extremal path to Schur \(\sigma\)-group, log ord \(20\), with TKT \(\mathrm{F}.7\)}
\label{fig:F7Lo20}


{\tiny

\setlength{\unitlength}{0.7cm}
\begin{picture}(14,29)(0,-28)


\put(0,0){\line(0,-1){26}}
\multiput(-0.1,0)(0,-2){14}{\line(1,0){0.2}}

\put(-0.2,0){\makebox(0,0)[rc]{\(2\,187\)}}
\put(0.2,0){\makebox(0,0)[lc]{\(3^7\)}}
\put(-0.2,-2){\makebox(0,0)[rc]{\(6\,561\)}}
\put(0.2,-2){\makebox(0,0)[lc]{\(3^8\)}}
{\color{blue}
\put(-0.2,-4){\makebox(0,0)[rc]{\(19\,683\)}}
\put(0.2,-4){\makebox(0,0)[lc]{\(3^9\)}}
}
\put(-0.2,-6){\makebox(0,0)[rc]{\(59\,049\)}}
\put(0.2,-6){\makebox(0,0)[lc]{\(3^{10}\)}}
\put(-0.2,-8){\makebox(0,0)[rc]{\(177\,147\)}}
\put(0.2,-8){\makebox(0,0)[lc]{\(3^{11}\)}}
\put(-0.2,-10){\makebox(0,0)[rc]{\(531\,441\)}}
\put(0.2,-10){\makebox(0,0)[lc]{\(3^{12}\)}}
\put(-0.2,-12){\makebox(0,0)[rc]{\(1\,594\,323\)}}
\put(0.2,-12){\makebox(0,0)[lc]{\(3^{13}\)}}
\put(-0.2,-14){\makebox(0,0)[rc]{\(4\,782\,969\)}}
\put(0.2,-14){\makebox(0,0)[lc]{\(3^{14}\)}}
\put(-0.2,-16){\makebox(0,0)[rc]{\(14\,348\,907\)}}
\put(0.2,-16){\makebox(0,0)[lc]{\(3^{15}\)}}
\put(-0.2,-18){\makebox(0,0)[rc]{\(43\,046\,721\)}}
\put(0.2,-18){\makebox(0,0)[lc]{\(3^{16}\)}}
\put(-0.2,-20){\makebox(0,0)[rc]{\(129\,140\,163\)}}
\put(0.2,-20){\makebox(0,0)[lc]{\(3^{17}\)}}
\put(-0.2,-22){\makebox(0,0)[rc]{\(387\,420\,489\)}}
\put(0.2,-22){\makebox(0,0)[lc]{\(3^{18}\)}}
\put(-0.2,-24){\makebox(0,0)[rc]{\(1\,162\,261\,467\)}}
\put(0.2,-24){\makebox(0,0)[lc]{\(3^{19}\)}}
{\color{red}
\put(-0.2,-26){\makebox(0,0)[rc]{\(3\,486\,784\,401\)}}
\put(0.2,-26){\makebox(0,0)[lc]{\(3^{20}\)}}
}

\put(0,-26){\vector(0,-1){1}}
\put(0,-27.5){\makebox(0,0)[ct]{Order \(3^n\)}}


{\color{red}
\put(5.8,0.4){\makebox(0,0)[rc]{fork}}
\put(6.2,0.4){\makebox(0,0)[lc]{\(\pi{M}=\pi^5{S}=P_7=\langle 3^7,64\rangle\)}}
\put(6.2,0.1){\makebox(0,0)[lc]{\(\mathrm{b}.10\)}}

\put(6,0){\circle*{0.2}}
}

\put(6,0){\line(-1,-2){2}}

{\color{blue}
\put(4,-4){\circle*{0.2}}
\put(3.8,-4.4){\makebox(0,0)[rc]{metabelianization}}
\put(4.2,-4.4){\makebox(0,0)[lc]{\(M=P_7-\#2;55\)}}
\put(4.2,-4.7){\makebox(0,0)[lc]{\(\mathrm{F}.7\)}}
}

{\color{red}
\put(6,0){\line(1,-4){2}}
\put(7.3,-4.7){\makebox(0,0)[lc]{\(s_5=4\)}}

\put(7.9,-8.1){\framebox(0.2,0.2){}}
\put(8.2,-8.4){\makebox(0,0)[lc]{\(\pi^4{S}=P_7-\#4;196\) (or 121 or 131)}}
\put(8.2,-8.7){\makebox(0,0)[lc]{\(\mathrm{F}.7\)}}

\put(8,-8){\line(0,-1){4}}
\put(8.2,-10){\makebox(0,0)[lc]{\(s_4=2\)}}

\put(7.9,-12.1){\framebox(0.2,0.2){}}
\put(8.3,-12.4){\makebox(0,0)[lc]{\(\pi^3{S}=\pi^4{S}-\#2;31\)}}
\put(8.3,-12.7){\makebox(0,0)[lc]{\(\mathrm{F}.7\)}}

\put(8,-12){\line(1,-4){2}}
\put(9.2,-16){\makebox(0,0)[lc]{\(s_3=4\)}}

\put(9.7,-20){\makebox(0,0)[rc]{one among \(27\) possible selections}}
\put(9.9,-20.1){\framebox(0.2,0.2){}}
\put(10.2,-20.4){\makebox(0,0)[lc]{\(\pi^2{S}=\pi^3{S}-\#4;1\)}}
\put(10.2,-20.7){\makebox(0,0)[lc]{\(\mathrm{F}.7\)}}

\put(10,-20){\line(0,-1){2}}
\put(10.2,-21.2){\makebox(0,0)[lt]{\(s_2=1\)}}

\put(9.9,-22.1){\framebox(0.2,0.2){}}
\put(10.4,-22.4){\makebox(0,0)[lc]{\(\pi{S}=\pi^2{S}-\#1;2\)}}
\put(10.4,-22.7){\makebox(0,0)[lc]{\(\mathrm{F}.7\)}}

\put(10,-22){\line(1,-2){2}}
\put(11.2,-24){\makebox(0,0)[lc]{\(s_1=2\)}}

\put(11.9,-26.1){\framebox(0.2,0.2){}}
\put(12,-26){\circle*{0.1}}
\put(12.2,-26.4){\makebox(0,0)[lc]{\(S=\pi{S}-\#2;1\)}}
\put(12.2,-26.7){\makebox(0,0)[lc]{\(\mathrm{F}.7\)}}

\put(8.2,-1.8){\makebox(0,0)[lc]{Topology Symbol:}}
\put(8.2,-2.4){\makebox(0,0)[lc]{\(\mathrm{F}\binom{2}{\rightarrow}\mathrm{b}\binom{4}{\leftarrow}\mathrm{F}\binom{2}{\leftarrow}\mathrm{F}\binom{4}{\leftarrow}\mathrm{F}\binom{1}{\leftarrow}\mathrm{F}\binom{2}{\leftarrow}\mathrm{F}\)}}
}


\end{picture}

}

\end{figure}

}

\newpage


{\tiny

\begin{figure}[hb]
\caption{\textbf{Incomplete} path, \textbf{four-stage} Schur \(\sigma\)-group, log ord \(>37\), TKT \(\mathrm{F}.13\)}
\label{fig:F13Lo37}


{\tiny

\setlength{\unitlength}{0.4cm}
\begin{picture}(17,49)(0,-48)


\put(-2,0){\line(0,-1){46}}
\multiput(-2.1,0)(0,-2){23}{\line(1,0){0.2}}

\put(-2.2,0){\makebox(0,0)[rc]{\(2\,187\)}}
\put(-1.8,0){\makebox(0,0)[lc]{\(3^7\)}}
\put(-2.2,-2){\makebox(0,0)[rc]{\(6\,561\)}}
\put(-1.8,-2){\makebox(0,0)[lc]{\(3^8\)}}
{\color{blue}
\put(-2.2,-4){\makebox(0,0)[rc]{\(19\,683\)}}
\put(-1.8,-4){\makebox(0,0)[lc]{\(3^9\)}}
}
\put(-2.2,-6){\makebox(0,0)[rc]{\(59\,049\)}}
\put(-1.8,-6){\makebox(0,0)[lc]{\(3^{10}\)}}
\put(-2.2,-8){\makebox(0,0)[rc]{\(177\,147\)}}
\put(-1.8,-8){\makebox(0,0)[lc]{\(3^{11}\)}}
\put(-2.2,-10){\makebox(0,0)[rc]{\(531\,441\)}}
\put(-1.8,-10){\makebox(0,0)[lc]{\(3^{12}\)}}
\put(-2.2,-12){\makebox(0,0)[rc]{\(1\,594\,323\)}}
\put(-1.8,-12){\makebox(0,0)[lc]{\(3^{13}\)}}
\put(-2.2,-14){\makebox(0,0)[rc]{\(4\,782\,969\)}}
\put(-1.8,-14){\makebox(0,0)[lc]{\(3^{14}\)}}
\put(-2.2,-16){\makebox(0,0)[rc]{\(14\,348\,907\)}}
\put(-1.8,-16){\makebox(0,0)[lc]{\(3^{15}\)}}
\put(-2.2,-18){\makebox(0,0)[rc]{\(43\,046\,721\)}}
\put(-1.8,-18){\makebox(0,0)[lc]{\(3^{16}\)}}
\put(-2.2,-20){\makebox(0,0)[rc]{\(129\,140\,163\)}}
\put(-1.8,-20){\makebox(0,0)[lc]{\(3^{17}\)}}
\put(-2.2,-22){\makebox(0,0)[rc]{\(387\,420\,489\)}}
\put(-1.8,-22){\makebox(0,0)[lc]{\(3^{18}\)}}
\put(-2.2,-24){\makebox(0,0)[rc]{\(1\,162\,261\,467\)}}
\put(-1.8,-24){\makebox(0,0)[lc]{\(3^{19}\)}}
\put(-2.2,-26){\makebox(0,0)[rc]{\(3\,486\,784\,401\)}}
\put(-1.8,-26){\makebox(0,0)[lc]{\(3^{20}\)}}
\put(-2.2,-28){\makebox(0,0)[rc]{\(10\,460\,353\,203\)}}
\put(-1.8,-28){\makebox(0,0)[lc]{\(3^{21}\)}}
\put(-2.2,-30){\makebox(0,0)[rc]{\(31\,381\,059\,609\)}}
\put(-1.8,-30){\makebox(0,0)[lc]{\(3^{22}\)}}
\put(-2.2,-32){\makebox(0,0)[rc]{\(94\,143\,178\,827\)}}
\put(-1.8,-32){\makebox(0,0)[lc]{\(3^{23}\)}}
\put(-2.2,-34){\makebox(0,0)[rc]{\(282\,429\,536\,481\)}}
\put(-1.8,-34){\makebox(0,0)[lc]{\(3^{24}\)}}
\put(-2.2,-36){\makebox(0,0)[rc]{\(847\,288\,609\,443\)}}
\put(-1.8,-36){\makebox(0,0)[lc]{\(3^{25}\)}}
\put(-2.2,-38){\makebox(0,0)[rc]{\(2\,541\,865\,828\,329\)}}
\put(-1.8,-38){\makebox(0,0)[lc]{\(3^{26}\)}}
\put(-2.2,-40){\makebox(0,0)[rc]{\(7\,625\,597\,484\,987\)}}
\put(-1.8,-40){\makebox(0,0)[lc]{\(3^{27}\)}}
\put(-2.2,-42){\makebox(0,0)[rc]{\(22\,876\,792\,454\,961\)}}
\put(-1.8,-42){\makebox(0,0)[lc]{\(3^{28}\)}}
\put(-2.2,-44){\makebox(0,0)[rc]{\(68\,630\,377\,364\,883\)}}
\put(-1.8,-44){\makebox(0,0)[lc]{\(3^{29}\)}}
\put(-2.2,-46){\makebox(0,0)[rc]{\(205\,891\,132\,094\,649\)}}
\put(-1.8,-46){\makebox(0,0)[lc]{\(3^{30}\)}}

\put(-2,-46){\vector(0,-1){1}}
\put(-2,-47.5){\makebox(0,0)[ct]{Order \(3^n\)}}


{\color{red}
\put(5.8,0.4){\makebox(0,0)[rc]{fork}}
\put(6.2,0.4){\makebox(0,0)[lc]{\(\pi{M}=\pi^{c-4}{S}=P_7=\langle 3^7,64\rangle\)}}
\put(6.2,-0.1){\makebox(0,0)[lc]{\(\mathrm{b}.10\)}}

\put(6,0){\circle*{0.2}}
}

\put(6,0){\line(-1,-2){2}}

{\color{blue}
\put(4,-4){\circle*{0.2}}
\put(3.8,-4.4){\makebox(0,0)[rc]{metabelianization}}
\put(4.2,-4.4){\makebox(0,0)[lc]{\(M=P_7-\#2;41\)}}
\put(4.2,-4.9){\makebox(0,0)[lc]{\(\mathrm{F}.13\)}}
}

{\color{red}
\put(6,0){\line(1,-4){2}}
\put(7.3,-4.9){\makebox(0,0)[lc]{\(s_{c-4}=4\)}}

\put(7.9,-8.1){\framebox(0.2,0.2){}}
\put(8.2,-8.4){\makebox(0,0)[lc]{\(\pi^{c-5}{S}=P_7-\#4;9\) (or 15 or 18 or 20)}}
\put(8.2,-8.9){\makebox(0,0)[lc]{\(\mathrm{F}.13\)}}

\put(8,-8){\line(0,-1){4}}
\put(8.2,-10){\makebox(0,0)[lc]{\(s_{c-5}=2\)}}

\put(7.9,-12.1){\framebox(0.2,0.2){}}
\put(8.4,-12.4){\makebox(0,0)[lc]{\(\pi^{c-6}{S}=\pi^{c-5}{S}-\#2;10\)}}
\put(8.4,-12.9){\makebox(0,0)[lc]{\(\mathrm{F}.13\)}}

\put(8,-12){\line(1,-4){2}}
\put(9.2,-16){\makebox(0,0)[lc]{\(s_{c-6}=4\)}}

\put(9.7,-20){\makebox(0,0)[rc]{one among \(12\) possible selections}}
\put(9.9,-20.1){\framebox(0.2,0.2){}}
\put(10.2,-20.4){\makebox(0,0)[lc]{\(\pi^{c-7}{S}=\pi^{c-6}{S}-\#4;144\)}}
\put(10.2,-20.9){\makebox(0,0)[lc]{\(\mathrm{F}.13\)}}

\put(10,-20){\line(0,-1){10}}
\put(10.2,-25){\makebox(0,0)[lc]{\(s_{c-7}=5\)}}

\put(9.9,-30.1){\framebox(0.2,0.2){}}
\put(10.4,-30.4){\makebox(0,0)[lc]{\(\pi^{c-8}{S}=\pi^{c-7}{S}-\#5;2516\)}}
\put(10.4,-30.9){\makebox(0,0)[lc]{\(\mathrm{F}.13\)}}

\put(10,-30){\line(1,-4){4}}
\put(12.2,-38){\makebox(0,0)[lc]{\(s_{c-8}=8\)}}

\put(13.9,-46.1){\framebox(0.2,0.2){}}
\put(14.2,-46.4){\makebox(0,0)[lc]{\(\pi^{c-9}{S}=\pi^{c-8}{S}-\#8;1\)}}
\put(14.2,-46.9){\makebox(0,0)[lc]{\(\mathrm{F}.13\)}}

\put(14,-46){\vector(0,-1){2}}
\put(14,-48.3){\makebox(0,0)[cc]{\(\mathrm{dl}=4\)}}


\put(8.2,-1.5){\makebox(0,0)[lc]{Topology Symbol:}}
\put(8.2,-2.4){\makebox(0,0)[lc]{\(\mathrm{F}\binom{2}{\rightarrow}\mathrm{b}\binom{4}{\leftarrow}\mathrm{F}\binom{2}{\leftarrow}\mathrm{F}\binom{4}{\leftarrow}\mathrm{F}\binom{5}{\leftarrow}\mathrm{F}\binom{8}{\leftarrow}\mathrm{F}\binom{7}{\leftarrow}\mathrm{F}\binom{12}{\leftarrow}\mathrm{F}\cdots\)}}
}


\end{picture}

}

\end{figure}
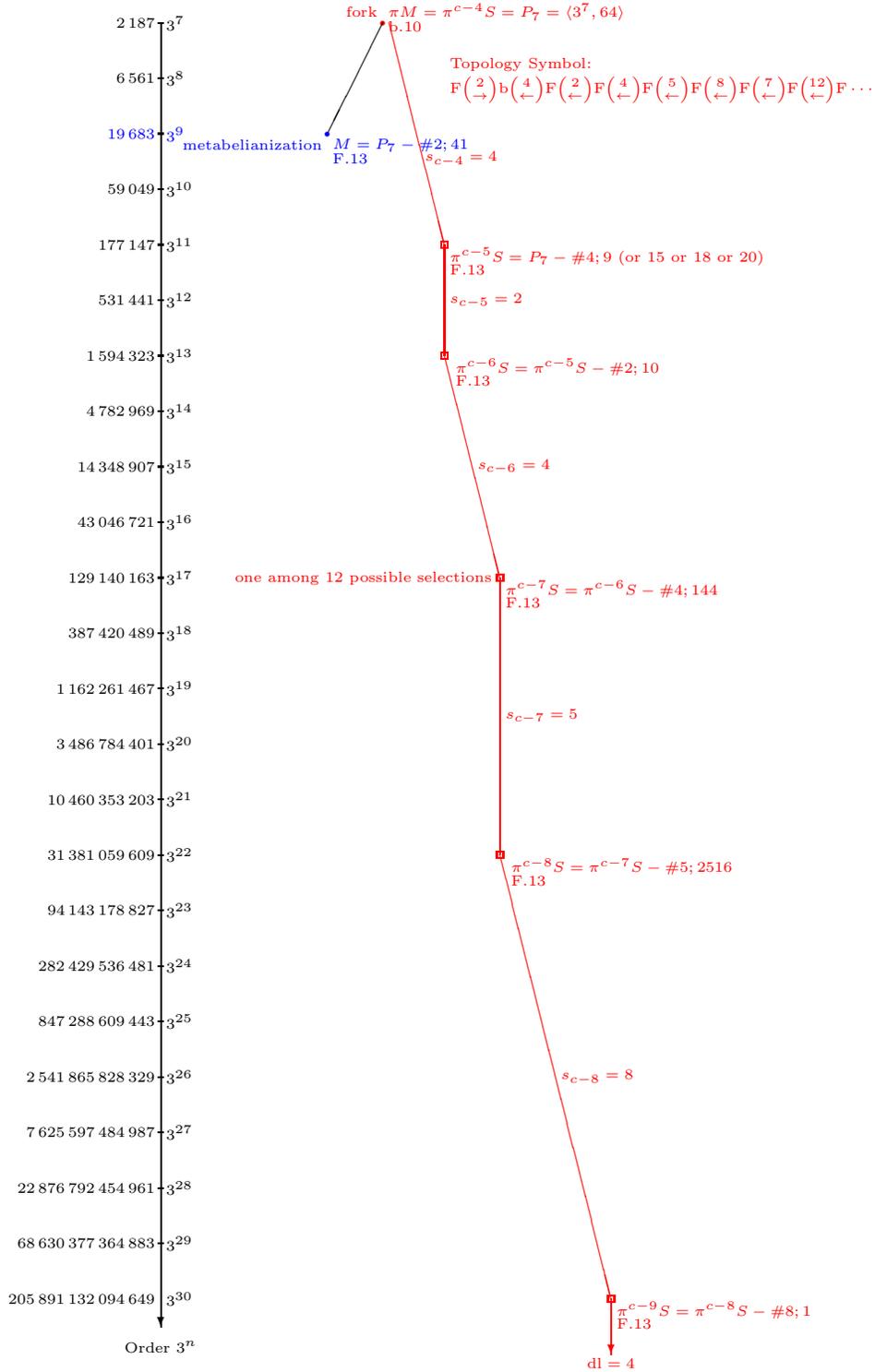

}

\newpage


\subsection{Schur \(\sigma\)-groups with TKT in other Sections}
\label{ss:SectionG}

We have seen that the seven \(3\)-groups with order \(3^5=243\), class \(3\) and coclass \(2\),
namely \(\langle 243,i\rangle\) with \(3\le i\le 9\),
which form the stem \(\Phi_6(0)\) of Hall's isoclinism family \(\Phi_6\),
are crucial roots of various descendant trees or even coclass trees
accommodating either Schur \(\sigma\)-groups, e.g. for \(i=4\) in Proposition
\ref{prp:SporadicSectionH},
or metabelianizations of Schur \(\sigma\)-groups, e.g. for \(i\in\lbrace 6,8\rbrace\) in Proposition
\ref{prp:SectionE},
whereas for \(i\in\lbrace 5,7\rbrace\)
they are isolated Schur \(\sigma\)-groups themselves (Proposition
\ref{prp:SectionD}).
All these vertices are connected with their common parent,
the extra special group \(\langle 27,3\rangle\),
by edges of step size \(s=2\).
We also realized that the descendant \(\langle 2187,64\rangle\)
of \(\langle 243,3\rangle\) with step size \(s=2\)
is the fork which is responsible for all \(\sigma\)-groups with types in Section F,
but also for groups with types H.4 and G.16.

Up to this point,
we did not consider the root \(\langle 243,9\rangle\) with TKT G.19.
It gives rise to a rather complicated descendant tree containing bifurcations.
To our knowledge it is unknown if the tree is infinite
and associated with a limit group like \(\langle 243,4\rangle\).
At the first bifurcation \(W:=\langle 729,57\rangle\),
there arise six descendants \(\langle 6561,i\rangle\), \(625\le i\le 630\), with step size \(s=2\)
having very different properties,
as Figure
\ref{fig:G19}
shows.


{\color{blue}
\begin{proposition}
\label{prp:G19}
There exist precisely four Schur \(\sigma\)-groups with transfer kernel type \(\mathrm{G}.19\)
and log ord \(11\).
They are the smallest groups with these properties and can be given in the shape
\(\langle 6561,i\rangle-\#1;2-\#2;j\) with \(i\in\lbrace 625,629\rbrace\) and \(1\le j\le 2\).
The descendant trees of \(\langle 6561,i\rangle\) with \(i\in\lbrace 625,629\rbrace\)
are finite with depth two.

There exist precisely twelve Schur \(\sigma\)-groups with transfer kernel type \(\mathrm{G}.19\)
and log ord \(14\).
They can be given in the shape
\(\langle 6561,627\rangle-\#1;2-\#2;j-\#1;k-\#2;\ell\) with \((j,k)\in\lbrace (1,3),(2,2)\rbrace\)
and \(1\le\ell\le 3\),
and six similar descendants of \(\langle 6561,628\rangle\).
The descendant tree of \(\langle 6561,627\rangle\) 
is finite with depth four.
\end{proposition}
}

\begin{proof}
This is the result of a construction of the descendant tree of \(\langle 243,9\rangle\)
up to logarithmic order \(14\)
by means of Magma
\cite{MAGMA}.
See Figure
\ref{fig:G19},
where the descendant subtrees of \(\langle 6561,i\rangle\) with even \(i\) are not drawn.
The subtree of \(\langle 6561,626\rangle\) does not contain many \(\sigma\)-groups,
in particular no Schur \(\sigma\)-group with log ord up to \(14\).
The subtree of \(\Psi:=\langle 6561,628\rangle\) starts similar as
the subtree of \(\Phi:=\langle 6561,627\rangle\).
The subtree of \(\langle 6561,630\rangle\) is the wildest of the six,
but does not contain Schur \(\sigma\)-groups with log ord up to \(14\).
\end{proof}


\begin{example}
\label{exm:G19}
The imaginary quadratic field \(K=\mathbb{Q}(\sqrt{d})\) with discriminant \(d=-114\,936\)
possesses a \(3\)-class field tower with precisely three stages, \(\ell_3(K)=3\),
having automorphism group \(G=\mathrm{Gal}(\mathrm{F}_3^\infty(K)/K)\simeq\langle 6561,629\rangle-\#1;2-\#2;j\)
with \(1\le j\le 2\).
This was proved by ourselves with the aid of abelian type invariants of second order in
\cite{Ma2017}.
Path see Table
\ref{tbl:G19Simple}
and Figure
\ref{fig:G19}.
\end{example}


\begin{hypothesis}
\label{hth:G19}
The imaginary quadratic field \(K=\mathbb{Q}(\sqrt{d})\) with discriminant \(d=-12\,067\)
(the smallest absolute discriminant with TKT G.19)
has probably a \(3\)-class field tower with exactly three stages, \(\ell_3(K)=3\),
and Galois group \(G=\mathrm{Gal}(\mathrm{F}_3^\infty(K)/K)\simeq\langle 6561,625\rangle-\#1;2-\#2;j\)
with \(1\le j\le 2\).
This is only a conjecture,
motivated by the bigger probability of candidates with smallest possible order,
but the abelian type invariants of second order of this field
also occur for other Schur \(\sigma\)-groups,
which causes uncertainty.
\end{hypothesis}

\newpage


{\tiny

\begin{figure}[hb]
\caption{Extremal paths to Schur \(\sigma\)-groups, log ord \(11\), purged tree \(\mathcal{T}_\ast(\langle 243,9\rangle)\)}
\label{fig:G19}


\setlength{\unitlength}{0.8cm}
\begin{picture}(15,21)(0,-20)

\put(0,0.5){\makebox(0,0)[cb]{Order \(3^n\)}}

\put(0,0){\line(0,-1){18}}
\multiput(-0.1,0)(0,-2){10}{\line(1,0){0.2}}

\put(-0.2,0){\makebox(0,0)[rc]{\(243\)}}
\put(0.2,0){\makebox(0,0)[lc]{\(3^5\)}}
{\color{blue}
\put(-0.2,-2){\makebox(0,0)[rc]{\(729\)}}
\put(0.2,-2){\makebox(0,0)[lc]{\(3^6\)}}
}
\put(-0.2,-4){\makebox(0,0)[rc]{\(2\,187\)}}
\put(0.2,-4){\makebox(0,0)[lc]{\(3^7\)}}
\put(-0.2,-6){\makebox(0,0)[rc]{\(6\,561\)}}
\put(0.2,-6){\makebox(0,0)[lc]{\(3^8\)}}
\put(-0.2,-8){\makebox(0,0)[rc]{\(19\,683\)}}
\put(0.2,-8){\makebox(0,0)[lc]{\(3^9\)}}
\put(-0.2,-10){\makebox(0,0)[rc]{\(59\,049\)}}
\put(0.2,-10){\makebox(0,0)[lc]{\(3^{10}\)}}
{\color{red}
\put(-0.2,-12){\makebox(0,0)[rc]{\(177\,147\)}}
\put(0.2,-12){\makebox(0,0)[lc]{\(3^{11}\)}}
}
\put(-0.2,-14){\makebox(0,0)[rc]{\(531\,441\)}}
\put(0.2,-14){\makebox(0,0)[lc]{\(3^{12}\)}}
\put(-0.2,-16){\makebox(0,0)[rc]{\(1\,594\,323\)}}
\put(0.2,-16){\makebox(0,0)[lc]{\(3^{13}\)}}
\put(-0.2,-18){\makebox(0,0)[rc]{\(4\,782\,969\)}}
\put(0.2,-18){\makebox(0,0)[lc]{\(3^{14}\)}}

\put(0,-18){\vector(0,-1){2}}



{\color{red}
\put(1.8,0.2){\makebox(0,0)[rc]{\(\langle 9\rangle\)}}
\put(2.2,0.2){\makebox(0,0)[lc]{root}}
\put(2,0){\circle*{0.2}}

\put(2,0){\line(0,-1){2}}

\put(2.2,-1.8){\makebox(0,0)[lc]{\(1^{\text{st}}\) bifurcation}}
}
{\color{blue}
\put(2.2,-1.3){\makebox(0,0)[lc]{metabelianization}}
\put(1.8,-1.3){\makebox(0,0)[rc]{\(W=\)}}
\put(1.8,-1.8){\makebox(0,0)[rc]{\(\langle 57\rangle\)}}
\put(2,-2){\circle*{0.1}}
}

\put(2,-2){\line(0,-1){2}}

\put(1.8,-3.8){\makebox(0,0)[rc]{\(\langle 311\rangle\)}}
\put(1.9,-4.1){\framebox(0.2,0.2){}}





{\color{red}
\put(2,-2){\line(1,-2){2}}

\put(4.2,-5.8){\makebox(0,0)[lc]{\(\langle 625\rangle\)}}
\put(3.9,-6.1){\framebox(0.2,0.2){}}


\put(4,-6){\line(0,-1){2}}

\put(4.2,-7.8){\makebox(0,0)[lc]{\(1;2\)}}
\put(3.8,-7.8){\makebox(0,0)[rc]{\(2^{\text{nd}}\) bifurcations}}
\put(3.95,-8.05){\framebox(0.1,0.1){}}
}

\put(4,-8){\line(1,-2){1}}

\put(5.2,-9.8){\makebox(0,0)[lc]{\(1;1\)}}
\put(4.9,-10.1){\framebox(0.2,0.2){}}

{\color{red}
\put(4,-8){\line(0,-1){4}}

\put(4.2,-11.8){\makebox(0,0)[lc]{\(2;1..2\)}}
\put(3.8,-11.8){\makebox(0,0)[rc]{\(2\ast\)}}
\put(3.9,-12.1){\framebox(0.2,0.2){}}
\put(4,-12){\circle*{0.1}}

\put(4.2,-12){\oval(2.6,1.3)}
\put(4,-12.9){\makebox(0,0)[cc]{\underbar{\textbf{-12\,067\ ?}}}}
}


\put(2,-2){\line(1,-1){4}}

\put(6.2,-5.8){\makebox(0,0)[lc]{\(\langle 626\rangle\)}}
\put(5.9,-6.1){\framebox(0.2,0.2){}}


\put(2,-2){\line(3,-2){6}}

\put(8.2,-5.8){\makebox(0,0)[lc]{\(\langle 627\rangle\)}}
\put(8.2,-6.1){\makebox(0,0)[lc]{\(=\Phi\)}}
\put(7.9,-6.1){\framebox(0.2,0.2){}}

\put(8,-6){\line(0,-1){2}}

\put(8.2,-7.8){\makebox(0,0)[lc]{\(1;2\)}}
\put(7.95,-8.05){\framebox(0.1,0.1){}}

\put(8,-8){\line(1,-1){2}}

\put(10.2,-9.8){\makebox(0,0)[lc]{\(1;1\)}}
\put(9.9,-10.1){\framebox(0.2,0.2){}}

\put(8,-8){\line(0,-1){4}}

\put(8.2,-11.8){\makebox(0,0)[lc]{\(2;1\)}}
\put(7.9,-12.1){\framebox(0.2,0.2){}}

\put(8,-12){\line(0,-1){2}}

\put(8.2,-13.8){\makebox(0,0)[lc]{\(1;3\)}}
\put(7.8,-13.8){\makebox(0,0)[rc]{\(3^{\text{rd}}\) bifurcations}}
\put(7.95,-14.05){\framebox(0.1,0.1){}}

\put(8,-14){\line(1,-2){1}}

\put(9.1,-15.8){\makebox(0,0)[lc]{\(1;1..2\)}}
\put(8.8,-15.8){\makebox(0,0)[rc]{\(2\ast\)}}
\put(8.9,-16.1){\framebox(0.2,0.2){}}

\put(8,-14){\line(0,-1){4}}

\put(8.2,-17.8){\makebox(0,0)[lc]{\(2;1..3\)}}
\put(7.8,-17.8){\makebox(0,0)[rc]{\(3\ast\)}}
\put(7.9,-18.1){\framebox(0.2,0.2){}}
\put(8,-18){\circle*{0.1}}

\put(8,-8){\line(1,-2){2}}

\put(10.2,-11.8){\makebox(0,0)[lc]{\(2;2\)}}
\put(9.9,-12.1){\framebox(0.2,0.2){}}

\put(10,-12){\line(0,-1){2}}

\put(10.2,-13.8){\makebox(0,0)[lc]{\(1;2\)}}
\put(9.95,-14.05){\framebox(0.1,0.1){}}

\put(10,-14){\line(1,-2){1}}

\put(11.2,-15.8){\makebox(0,0)[lc]{\(1;1\)}}
\put(10.9,-16.1){\framebox(0.2,0.2){}}

\put(10,-14){\line(0,-1){4}}

\put(10.2,-17.8){\makebox(0,0)[lc]{\(2;1..3\)}}
\put(9.8,-17.8){\makebox(0,0)[rc]{\(3\ast\)}}
\put(9.9,-18.1){\framebox(0.2,0.2){}}
\put(10,-18){\circle*{0.1}}


\put(2,-2){\line(2,-1){8}}

\put(10.2,-5.8){\makebox(0,0)[lc]{\(\langle 628\rangle\)}}
\put(10.2,-6.1){\makebox(0,0)[lc]{\(=\Psi\)}}
\put(9.9,-6.1){\framebox(0.2,0.2){}}

\put(8,-6){\line(0,-1){2}}

\put(8.2,-7.8){\makebox(0,0)[lc]{\(1;2\)}}
\put(7.95,-8.05){\framebox(0.1,0.1){}}



{\color{red}
\put(2,-2){\line(5,-2){10}}

\put(12.2,-5.8){\makebox(0,0)[lc]{\(\langle 629\rangle\)}}
\put(12.2,-6.1){\makebox(0,0)[lc]{\(=Y\)}}
\put(11.9,-6.1){\framebox(0.2,0.2){}}

\put(12,-6){\line(0,-1){2}}

\put(12.2,-7.8){\makebox(0,0)[lc]{\(1;2\)}}
\put(11.95,-8.05){\framebox(0.1,0.1){}}
}

\put(12,-8){\line(1,-2){1}}

\put(13.2,-9.8){\makebox(0,0)[lc]{\(1;1\)}}
\put(12.9,-10.1){\framebox(0.2,0.2){}}

{\color{red}
\put(12,-8){\line(0,-1){4}}

\put(12.2,-11.8){\makebox(0,0)[lc]{\(2;1..2\)}}
\put(11.8,-11.8){\makebox(0,0)[rc]{\(2\ast\)}}
\put(11.9,-12.1){\framebox(0.2,0.2){}}
\put(12,-12){\circle*{0.1}}
}



{\color{red}
\put(12.2,-12){\oval(2.6,1.3)}
\put(12,-12.9){\makebox(0,0)[cc]{\underbar{\textbf{-114\,936}}}}
}


\put(2,-2){\line(3,-1){12}}

\put(14.2,-5.8){\makebox(0,0)[lc]{\(\langle 630\rangle\)}}
\put(14.2,-6.1){\makebox(0,0)[lc]{\(=Z\)}}
\put(13.95,-6.05){\framebox(0.1,0.1){}}

\put(14,-6.0){\vector(0,-1){2}}


\end{picture}

\end{figure}
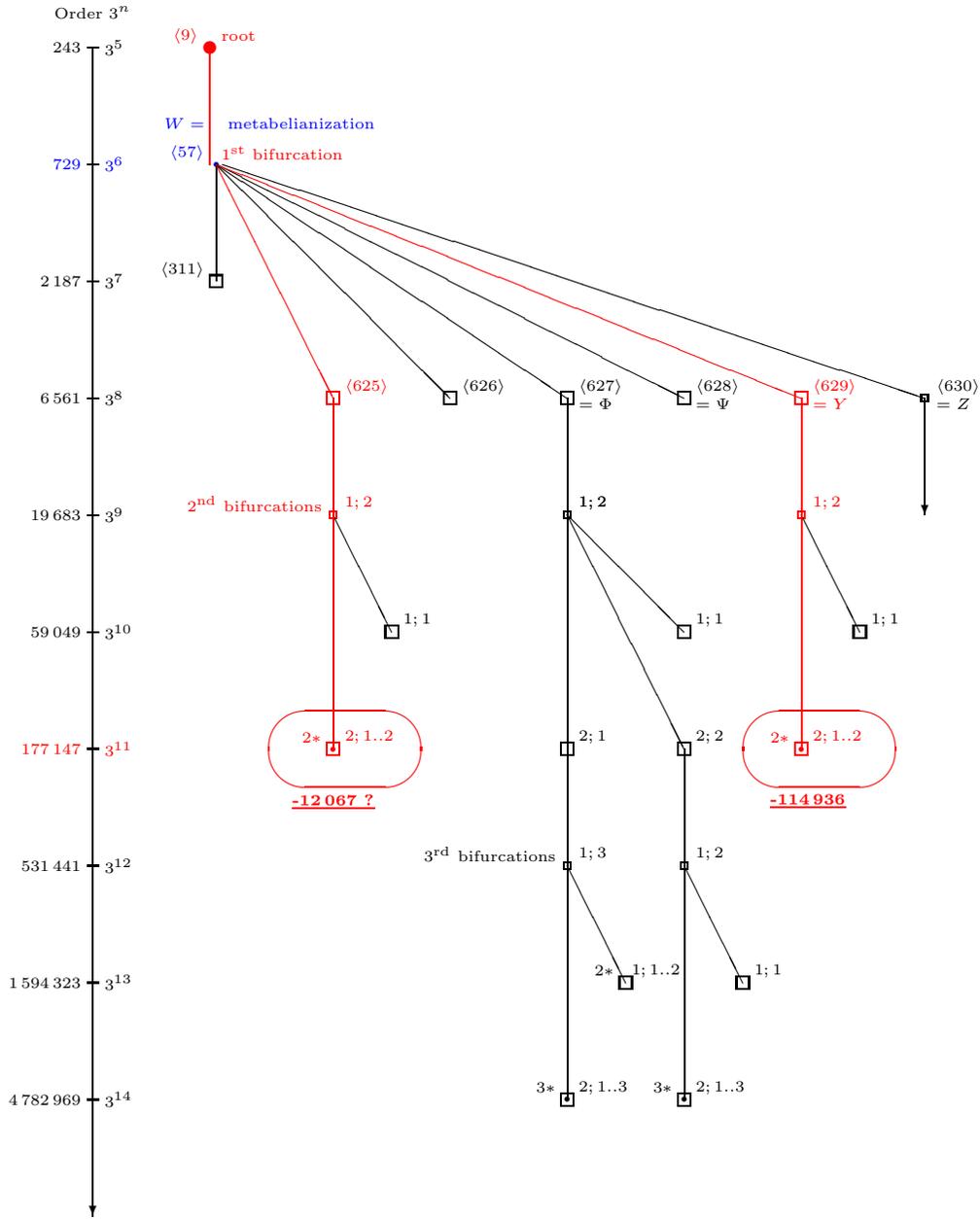

}

\newpage


\noindent

\renewcommand{\arraystretch}{1.1}

\begin{table}[hb]
\caption{Root path of \(G\) for the simplest case of transfer kernel type G.19}
\label{tbl:G19Simple}
\begin{center}
\begin{tabular}{|c|c|c||c|c|c|}
\hline
 Ancestor     & Vertex                      & lo & \((\nu,\mu)\)                   & \((N_s/C_s)_{1\le s\le\nu}\)                            & TKT  \\
\hline
 \(\pi^5(G)\) & \(\langle 27,3\rangle\)     &  3 & \(({\color{red}\mathbf{2,4}})\) & \((4/1,{\color{red}\mathbf{7/5}})\)                     & a.1  \\
 \(\pi^4(G)\) & \(\langle 243,9\rangle\)    &  5 & \(({\color{red}\mathbf{1,3}})\) & \(({\color{red}\mathbf{2/2}})\)                         & G.19 \\
 \(\pi^3(G)\) & \(\langle 729,57\rangle\)   &  6 & \(({\color{red}\mathbf{2,4}})\) & \((1/0,{\color{red}\mathbf{6/6}})\)                     & G.19 \\
 \(\pi^2(G)\) & \(\langle 6561,629\rangle\) &  8 & \(({\color{red}\mathbf{1,3}})\) & \(({\color{red}\mathbf{2/2}})\)                         & G.19 \\
 \(\pi(G)\)   & \(-\#1;2\)                  &  9 & \(({\color{red}\mathbf{2,4}})\) & \((1/0,{\color{red}\mathbf{2/0}})\)                     & G.19 \\
 \(G\)        & \(-\#2;1\)                  & 11 & \(({\color{red}\mathbf{0,2}})\) &                                                         & G.19 \\
\hline
\end{tabular}
\end{center}
\end{table}


\renewcommand{\arraystretch}{1.1}

\begin{table}[hb]
\caption{Root path of \(G\), log ord 20, for the irregular case of transfer kernel type H.4}
\label{tbl:H4Irregular}
\begin{center}
\begin{tabular}{|c|c|c||c|c|c|}
\hline
 Ancestor     & Vertex                     & lo & \((\nu,\mu)\)                   & \((N_s/C_s)_{1\le s\le\nu}\)                            & TKT  \\
\hline
 \(\pi^7(G)\) & \(\langle 27,3\rangle\)    &  3 & \(({\color{red}\mathbf{2,4}})\) & \((4/1,{\color{red}\mathbf{7/5}})\)                     & a.1  \\
 \(\pi^6(G)\) & \(\langle 243,3\rangle\)   &  5 & \(({\color{red}\mathbf{2,4}})\) & \((10/6,{\color{red}\mathbf{15/15}})\)                  & b.10 \\
 \(\pi^5(G)\) & \(\langle 2187,64\rangle\) &  7 & \(({\color{red}\mathbf{4,6}})\) & \((33/2,453/84,918/540,{\color{red}\mathbf{198/198}})\) & b.10 \\
 \(\pi^4(G)\) & \(-\#4;111\)               & 11 & \(({\color{red}\mathbf{2,4}})\) & \((20/20,{\color{red}\mathbf{41/41}})\)                 & H.4  \\
 \(\pi^3(G)\) & \(-\#2;38\)                & 13 & \(({\color{red}\mathbf{4,6}})\) & \((54/0,324/27,270/108,{\color{red}\mathbf{27/27}})\)   & H.4  \\
 \(\pi^2(G)\) & \(-\#4;1\)                 & 17 & \(({\color{red}\mathbf{1,3}})\) & \(({\color{red}\mathbf{5/5}})\)                         & H.4  \\
 \(\pi(G)\)   & \(-\#1;5\)                 & 18 & \(({\color{red}\mathbf{2,4}})\) & \((4/0,{\color{red}\mathbf{1/0}})\)                     & H.4  \\
 \(G\)        & \(-\#2;1\)                 & 20 & \(({\color{red}\mathbf{0,2}})\) &                                                         & H.4  \\
\hline
\end{tabular}
\end{center}
\end{table}


\begin{hypothesis}
\label{hth:H4Irregular}
A possible extremal root path to the \(3\)-class tower group
\(G=\mathrm{Gal}(\mathrm{F}_3^\infty(K)/K)\simeq S\)
of \(K=\mathbb{Q}(\sqrt{-186\,483})\) with \textbf{irregular} type H.4
is described in Table
\ref{tbl:H4Irregular}
and drawn in Figure
\ref{fig:H4IrrLo20}.
Except for the additional edge with step size \(s=1\) from the metabelianization \(M\),
the graph is isomorphic to the graph concerning type F.7 in Figure
\ref{fig:F7Lo20}.
Similarly as in this preceding figure,
we do not know to which coclass trees the capable vertices on the path belong.
Here, we assume \(M=S/S^{\prime\prime}\simeq P_7-\#2;34-\#1;7\) and
\(S\simeq P_7-\#4;111-\#2;38-\#4;1-\#1;5-\#2;1\).
The type H.4 is \textbf{irregular}
in the sense of B. Nebelung
\cite{Ne1,Ne2},
since the commutator subgroup \(M^\prime\), resp. \(S^\prime\),
has abelian type, resp. quotient, invariants \((2^4)\hat{=}(9,9,9,9)\)
instead of the regular \((32^21)\hat{=}(27,9,9,3)\).

This paradigm shows that Schur \(\sigma\)-groups \(S\) of log ord \(20\)
are expected not only for types in Section F
but also for types G.16, G.19, H.4, both, regular and irregular,
always under assumption of a sporadic metabelianization \(M=S/S^{\prime\prime}\)
of coclass \(\mathrm{cc}(M)=4\), either of log ord \(9\) or \(10\). 
\end{hypothesis}

\newpage


{\tiny

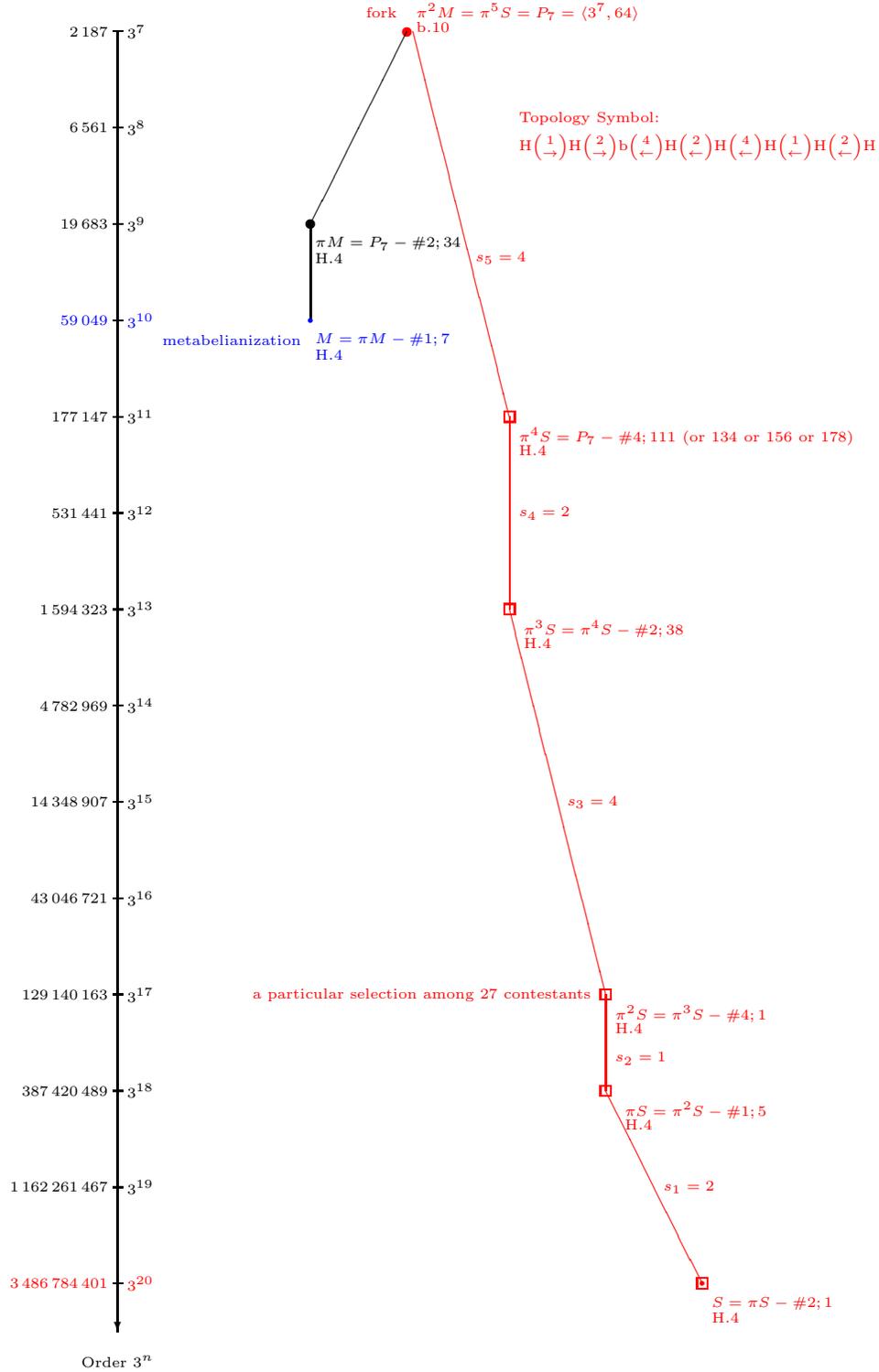
\begin{figure}[hb]
\caption{Extremal path to Schur \(\sigma\)-group, log ord \(20\), with irregular TKT \(\mathrm{H}.4\)}
\label{fig:H4IrrLo20}


{\tiny

\setlength{\unitlength}{0.7cm}
\begin{picture}(14,29)(0,-28)


\put(0,0){\line(0,-1){26}}
\multiput(-0.1,0)(0,-2){14}{\line(1,0){0.2}}

\put(-0.2,0){\makebox(0,0)[rc]{\(2\,187\)}}
\put(0.2,0){\makebox(0,0)[lc]{\(3^7\)}}
\put(-0.2,-2){\makebox(0,0)[rc]{\(6\,561\)}}
\put(0.2,-2){\makebox(0,0)[lc]{\(3^8\)}}
\put(-0.2,-4){\makebox(0,0)[rc]{\(19\,683\)}}
\put(0.2,-4){\makebox(0,0)[lc]{\(3^9\)}}
{\color{blue}
\put(-0.2,-6){\makebox(0,0)[rc]{\(59\,049\)}}
\put(0.2,-6){\makebox(0,0)[lc]{\(3^{10}\)}}
}
\put(-0.2,-8){\makebox(0,0)[rc]{\(177\,147\)}}
\put(0.2,-8){\makebox(0,0)[lc]{\(3^{11}\)}}
\put(-0.2,-10){\makebox(0,0)[rc]{\(531\,441\)}}
\put(0.2,-10){\makebox(0,0)[lc]{\(3^{12}\)}}
\put(-0.2,-12){\makebox(0,0)[rc]{\(1\,594\,323\)}}
\put(0.2,-12){\makebox(0,0)[lc]{\(3^{13}\)}}
\put(-0.2,-14){\makebox(0,0)[rc]{\(4\,782\,969\)}}
\put(0.2,-14){\makebox(0,0)[lc]{\(3^{14}\)}}
\put(-0.2,-16){\makebox(0,0)[rc]{\(14\,348\,907\)}}
\put(0.2,-16){\makebox(0,0)[lc]{\(3^{15}\)}}
\put(-0.2,-18){\makebox(0,0)[rc]{\(43\,046\,721\)}}
\put(0.2,-18){\makebox(0,0)[lc]{\(3^{16}\)}}
\put(-0.2,-20){\makebox(0,0)[rc]{\(129\,140\,163\)}}
\put(0.2,-20){\makebox(0,0)[lc]{\(3^{17}\)}}
\put(-0.2,-22){\makebox(0,0)[rc]{\(387\,420\,489\)}}
\put(0.2,-22){\makebox(0,0)[lc]{\(3^{18}\)}}
\put(-0.2,-24){\makebox(0,0)[rc]{\(1\,162\,261\,467\)}}
\put(0.2,-24){\makebox(0,0)[lc]{\(3^{19}\)}}
{\color{red}
\put(-0.2,-26){\makebox(0,0)[rc]{\(3\,486\,784\,401\)}}
\put(0.2,-26){\makebox(0,0)[lc]{\(3^{20}\)}}
}

\put(0,-26){\vector(0,-1){1}}
\put(0,-27.5){\makebox(0,0)[ct]{Order \(3^n\)}}


{\color{red}
\put(5.8,0.4){\makebox(0,0)[rc]{fork}}
\put(6.2,0.4){\makebox(0,0)[lc]{\(\pi^2{M}=\pi^5{S}=P_7=\langle 3^7,64\rangle\)}}
\put(6.2,0.1){\makebox(0,0)[lc]{\(\mathrm{b}.10\)}}

\put(6,0){\circle*{0.2}}
}

\put(6,0){\line(-1,-2){2}}

\put(4,-4){\circle*{0.2}}
\put(4.1,-4.4){\makebox(0,0)[lc]{\(\pi{M}=P_7-\#2;34\)}}
\put(4.1,-4.7){\makebox(0,0)[lc]{\(\mathrm{H}.4\)}}

\put(4,-4){\line(0,-1){2}}

{\color{blue}
\put(4,-6){\circle*{0.1}}
\put(3.8,-6.4){\makebox(0,0)[rc]{metabelianization}}
\put(4.1,-6.4){\makebox(0,0)[lc]{\(M=\pi{M}-\#1;7\)}}
\put(4.1,-6.7){\makebox(0,0)[lc]{\(\mathrm{H}.4\)}}
}

{\color{red}
\put(6,0){\line(1,-4){2}}
\put(7.3,-4.7){\makebox(0,0)[lc]{\(s_5=4\)}}

\put(7.9,-8.1){\framebox(0.2,0.2){}}
\put(8.2,-8.4){\makebox(0,0)[lc]{\(\pi^4{S}=P_7-\#4;111\) (or 134 or 156 or 178)}}
\put(8.2,-8.7){\makebox(0,0)[lc]{\(\mathrm{H}.4\)}}

\put(8,-8){\line(0,-1){4}}
\put(8.2,-10){\makebox(0,0)[lc]{\(s_4=2\)}}

\put(7.9,-12.1){\framebox(0.2,0.2){}}
\put(8.3,-12.4){\makebox(0,0)[lc]{\(\pi^3{S}=\pi^4{S}-\#2;38\)}}
\put(8.3,-12.7){\makebox(0,0)[lc]{\(\mathrm{H}.4\)}}

\put(8,-12){\line(1,-4){2}}
\put(9.2,-16){\makebox(0,0)[lc]{\(s_3=4\)}}

\put(9.7,-20){\makebox(0,0)[rc]{a particular selection among \(27\) contestants}}
\put(9.9,-20.1){\framebox(0.2,0.2){}}
\put(10.2,-20.4){\makebox(0,0)[lc]{\(\pi^2{S}=\pi^3{S}-\#4;1\)}}
\put(10.2,-20.7){\makebox(0,0)[lc]{\(\mathrm{H}.4\)}}

\put(10,-20){\line(0,-1){2}}
\put(10.2,-21.2){\makebox(0,0)[lt]{\(s_2=1\)}}

\put(9.9,-22.1){\framebox(0.2,0.2){}}
\put(10.4,-22.4){\makebox(0,0)[lc]{\(\pi{S}=\pi^2{S}-\#1;5\)}}
\put(10.4,-22.7){\makebox(0,0)[lc]{\(\mathrm{H}.4\)}}

\put(10,-22){\line(1,-2){2}}
\put(11.2,-24){\makebox(0,0)[lc]{\(s_1=2\)}}

\put(11.9,-26.1){\framebox(0.2,0.2){}}
\put(12,-26){\circle*{0.1}}
\put(12.2,-26.4){\makebox(0,0)[lc]{\(S=\pi{S}-\#2;1\)}}
\put(12.2,-26.7){\makebox(0,0)[lc]{\(\mathrm{H}.4\)}}

\put(8.2,-1.8){\makebox(0,0)[lc]{Topology Symbol:}}
\put(8.2,-2.4){\makebox(0,0)[lc]{\(\mathrm{H}\binom{1}{\rightarrow}\mathrm{H}\binom{2}{\rightarrow}\mathrm{b}\binom{4}{\leftarrow}\mathrm{H}\binom{2}{\leftarrow}\mathrm{H}\binom{4}{\leftarrow}\mathrm{H}\binom{1}{\leftarrow}\mathrm{H}\binom{2}{\leftarrow}\mathrm{H}\)}}
}


\end{picture}

}

\end{figure}

}

\newpage


\section{Real quadratic counterexamples}
\label{s:CounterExamples}

\noindent
According to Formula
\eqref{eqn:Shafarevich},
the less restrictive bounds \(2\le d_2(G)\le 3\) for the relation rank
of the Galois group \(G=\mathrm{Gal}(\mathrm{F}_3^\infty(K)/K)\)
of the \(3\)-class field tower
of a real quadratic number field \(K=\mathbb{Q}(\sqrt{d})\), \(d>0\),
with signature \((r_1,r_2)=(2,0)\) and Dirichlet unit rank \(r=r_1+r_2-1=1\)
give rise to a \textbf{totally different behavior},
in comparison to imaginary quadratic fields \(K\)
with elementary \(3\)-class group \(\mathrm{Cl}_3(K)\simeq C_3\times C_3\) of rank two.

Although the Schur \(\sigma\)-groups \(G\) with \(d_2(G)=2\) are still admissible for discriminants \(d>0\),
there arises a heavy competition by the much more numerous \(\sigma\)-groups \(G\) with relation rank \(d_2(G)=3\),
which also have the advantage of considerably lower order and higher probability.


{\color{red}
\begin{theorem}
\label{thm:CounterExampleCc1}
The Galois group \(G=\mathrm{Gal}(\mathrm{F}_3^\infty(K)/K)\) of the \(3\)-class field tower
of nearly \(90\%\) of all the \(34\,631\) real quadratic fields \(K=\mathbb{Q}(\sqrt{d})\)
with fundamental discriminants in the range \(0<d<10^8\)
and elementary \(3\)-class group \(\mathrm{Cl}_3(K)\simeq C_3\times C_3\) of rank two
is metabelian (i.e. two-stage) of coclass \(\mathrm{cc}(G)=1\),
and thus the root path of \(G\) does not follow the first possible bifurcation
at the common class-\(2\) quotient \(G/\gamma_3(G)\simeq\langle 27,3\rangle\) with nuclear rank \(\nu=2\).
Consequently the dominating part of real quadratic fields
\textbf{does not satisfy} the extremal path condition Conj.
\ref{cnj:ExtremalPath}.
\end{theorem}
}

\begin{proof}
The first edge \(G\to\pi(G)\) with step size \(s=1<\nu\) of the root path,
which \textbf{violates the extremal path condition},
is drawn in red color
for the ground state (GS) of
type a.2 (\(G\simeq\langle 81,10\rangle\) with fixed point TKT, for \(7\,104\) fields),
type a.3\(\ast\) (\(G\simeq\langle 81,7\rangle\) with irregular TTT, for \(10\,244\) fields),
and type a.3 (\(G\simeq\langle 81,8\rangle\) with regular TTT, for \(10\,514\) fields)
in Figure
\ref{fig:PeriCc1}.
Comparison of the percentage of all tower groups \(G\) with coclass \(1\)
(ground state and excited states)
corresponding to ranges of discriminants \(d\) with increasing upper bound
\(10^6\), \(10^7\), \(10^8\) in Table
\ref{tbl:RQ3x3}
reveals a growth from \(89.3\%\) over \(89.40\%\) and \(89.769\%\)
towards clearly dominating \(90\%\),
\textbf{quod erat demonstrandum}.
\end{proof}

The information in Table
\ref{tbl:RQ3x3}
and Figure
\ref{fig:PeriCc1}
has been computed by ourselves in \(2016\),
published in
\cite{Ma2016b}
presented in key note lectures
\cite{Ma2016c,Ma2018},
and refined in
\cite{Ma2018a}.



\renewcommand{\arraystretch}{1.1}

\begin{table}[ht]
\caption{Smallest \(G/G^{(2)}\simeq G\) for real quadratic fields of type \((3,3)\)}
\label{tbl:RQ3x3}
\begin{center}
\begin{tabular}{|c|r||r|r|r|r||r|l|}
\hline
 Discriminant     & Total\#     & total            & fixed point      & regular          & irregular         & \(\Sigma\)  & \(\%\)       \\
\hline
 \(0<d<10^6\)     &     \(149\) &     \(11\)       &     \(35\)       &      \(53\)      &      \(34\)       &     \(133\) & \(89.3\%\)   \\
 \(0<d<10^7\)     &  \(2\,576\) &    \(151\)       &       \(\)       &  \(1\,454\)      &     \(698\)       &  \(2\,303\) & \(89.40\%\)  \\
 \(0<d<10^8\)     & \(34\,631\) & \(2\,241\)       & \(7\,356\)       & \(11\,247\)      & \(10\,244\)       & \(31\,088\) & \(89.769\%\) \\
\hline
 \(\tau(G)\), GS  &             & \((32,(1^2)^3)\) & \((21,(1^2)^3)\) & \((21,(1^2)^3)\) & \((1^3,(1^2)^3)\) &             &  \\
 \(\varkappa(G)\) &             & \((0000)\)       & \((1000)\)       & \((2000)\)       & \((2000)\)        &             &  \\
 TKT              &             & a.1              & a.2              & a.3              & a.3\(\ast\)       &             &  \\
 GS: \(G\simeq\)  & & \(\langle 729,99\rangle\) & \(\langle 81,10\rangle\) & \(\langle 81,8\rangle\) & \(\langle 81,7\rangle\) & & \\
\hline
\end{tabular}
\end{center}
\end{table}


\begin{counterexample}
\label{cex:G19}
Even when the primary bifurcation at \(\langle 27,3\rangle\) is followed
by the root path of \(G\) for a real quadratic field \(K=\mathbb{Q}(\sqrt{d})\),
e.g. by the edge \(\langle 243,9\rangle\to\langle 27,3\rangle\) from TKT G.19,
the \(1^{\text{st}}\), resp. \(2^{\text{nd}}\), bifurcation
in the purged tree \(\mathcal{T}_\ast(\langle 243,9\rangle)\) 
may be violated 
by edges of step size \(s=1<\nu=2\),
as illustrated by Figure
\ref{fig:RealG19}.
For the field with discriminant \(d=+214\,712\)
\cite{Ma2012,Ma2014},
which is the smallest without any total capitulation
(discovered by ourselves in January 2006 already),
we have \(G\simeq\langle 2187,311\rangle\) of considerably lower log ord \(7\) instead of \(11\) in Figure
\ref{fig:G19},
and for the field with discriminant \(d=+21\,974\,161\)
\cite{Ma2016b,Ma2016c}
in much higher discriminantal ranges,
the edge to \(G\simeq Y_1-\#1;1\)
is abbreviated from log ord \(11\) to \(10\)
at \(Y_1\).
\end{counterexample}

\newpage

{\tiny

\begin{figure}[hb]
\caption{Non-extremal paths to unbalanced \(\sigma\)-groups, ord \(81\), coclass tree \(\mathcal{T}^1(\langle 9,2\rangle)\)}
\label{fig:PeriCc1}

\input{PeriCc1}

\end{figure}

}

\newpage


{\tiny

\begin{figure}[hb]
\caption{Non-extremal paths to unbalanced \(\sigma\)-groups, log ord \(7\) and \(10\)}
\label{fig:RealG19}


\setlength{\unitlength}{0.8cm}
\begin{picture}(15,21)(0,-20)

\put(0,0.5){\makebox(0,0)[cb]{Order \(3^n\)}}

\put(0,0){\line(0,-1){18}}
\multiput(-0.1,0)(0,-2){10}{\line(1,0){0.2}}

\put(-0.2,0){\makebox(0,0)[rc]{\(243\)}}
\put(0.2,0){\makebox(0,0)[lc]{\(3^5\)}}
{\color{blue}
\put(-0.2,-2){\makebox(0,0)[rc]{\(729\)}}
\put(0.2,-2){\makebox(0,0)[lc]{\(3^6\)}}
}
{\color{red}
\put(-0.2,-4){\makebox(0,0)[rc]{\(2\,187\)}}
\put(0.2,-4){\makebox(0,0)[lc]{\(3^7\)}}
}
\put(-0.2,-6){\makebox(0,0)[rc]{\(6\,561\)}}
\put(0.2,-6){\makebox(0,0)[lc]{\(3^8\)}}
\put(-0.2,-8){\makebox(0,0)[rc]{\(19\,683\)}}
\put(0.2,-8){\makebox(0,0)[lc]{\(3^9\)}}
{\color{red}
\put(-0.2,-10){\makebox(0,0)[rc]{\(59\,049\)}}
\put(0.2,-10){\makebox(0,0)[lc]{\(3^{10}\)}}
}
\put(-0.2,-12){\makebox(0,0)[rc]{\(177\,147\)}}
\put(0.2,-12){\makebox(0,0)[lc]{\(3^{11}\)}}
\put(-0.2,-14){\makebox(0,0)[rc]{\(531\,441\)}}
\put(0.2,-14){\makebox(0,0)[lc]{\(3^{12}\)}}
\put(-0.2,-16){\makebox(0,0)[rc]{\(1\,594\,323\)}}
\put(0.2,-16){\makebox(0,0)[lc]{\(3^{13}\)}}
\put(-0.2,-18){\makebox(0,0)[rc]{\(4\,782\,969\)}}
\put(0.2,-18){\makebox(0,0)[lc]{\(3^{14}\)}}

\put(0,-18){\vector(0,-1){2}}


\put(1.8,0.2){\makebox(0,0)[rc]{\(\langle 9\rangle\)}}
\put(2.2,0.2){\makebox(0,0)[lc]{root}}
\put(2,0){\circle*{0.2}}

\put(2,0){\line(0,-1){2}}

{\color{blue}
\put(2.2,-1.3){\makebox(0,0)[lc]{metabelianization}}
\put(1.8,-1.3){\makebox(0,0)[rc]{\(W=\)}}
\put(1.8,-1.8){\makebox(0,0)[rc]{\(\langle 57\rangle\)}}
\put(2,-2){\circle*{0.1}}
}

{\color{red}
\put(2.2,-1.8){\makebox(0,0)[lc]{\(1^{\text{st}}\) bifurcation}}

\put(2,-2){\line(0,-1){2}}
\put(1.8,-3.8){\makebox(0,0)[rc]{\(\langle 311\rangle\)}}
\put(1.9,-4.1){\framebox(0.2,0.2){}}

\put(1.8,-4){\oval(1.6,1.3)}
\put(2,-4.9){\makebox(0,0)[cc]{\underbar{\textbf{+214\,712}}}}
}




\put(2,-2){\line(1,-2){2}}

\put(4.2,-5.8){\makebox(0,0)[lc]{\(\langle 625\rangle\)}}
\put(3.9,-6.1){\framebox(0.2,0.2){}}


\put(4,-6){\line(0,-1){2}}

\put(4.2,-7.8){\makebox(0,0)[lc]{\(1;2\)}}

{\color{red}
\put(3.8,-7.8){\makebox(0,0)[rc]{\(2^{\text{nd}}\) bifurcations}}
}

\put(3.95,-8.05){\framebox(0.1,0.1){}}

\put(4,-8){\line(1,-2){1}}

\put(5.2,-9.8){\makebox(0,0)[lc]{\(1;1\)}}
\put(4.9,-10.1){\framebox(0.2,0.2){}}

\put(4,-8){\line(0,-1){4}}

\put(4.2,-11.8){\makebox(0,0)[lc]{\(2;1..2\)}}
\put(3.8,-11.8){\makebox(0,0)[rc]{\(2\ast\)}}
\put(3.9,-12.1){\framebox(0.2,0.2){}}

\put(4.2,-12){\oval(2.6,1.3)}
\put(4,-12.9){\makebox(0,0)[cc]{\underbar{\textbf{-12\,067\ ?}}}}


\put(2,-2){\line(1,-1){4}}

\put(6.2,-5.8){\makebox(0,0)[lc]{\(\langle 626\rangle\)}}
\put(6.2,-6.1){\makebox(0,0)[lc]{\(=\Phi\)}}
\put(5.9,-6.1){\framebox(0.2,0.2){}}


\put(2,-2){\line(3,-2){6}}

\put(8.2,-5.8){\makebox(0,0)[lc]{\(\langle 627\rangle\)}}
\put(7.9,-6.1){\framebox(0.2,0.2){}}

\put(8,-6){\line(0,-1){2}}

\put(8.2,-7.8){\makebox(0,0)[lc]{\(1;2\)}}
\put(7.95,-8.05){\framebox(0.1,0.1){}}

\put(8,-8){\line(1,-1){2}}

\put(10.2,-9.8){\makebox(0,0)[lc]{\(1;1\)}}
\put(9.9,-10.1){\framebox(0.2,0.2){}}

\put(8,-8){\line(0,-1){4}}

\put(8.2,-11.8){\makebox(0,0)[lc]{\(2;1\)}}
\put(7.9,-12.1){\framebox(0.2,0.2){}}

\put(8,-12){\line(0,-1){2}}

\put(8.2,-13.8){\makebox(0,0)[lc]{\(1;3\)}}
\put(7.8,-13.8){\makebox(0,0)[rc]{\(3^{\text{rd}}\) bifurcations}}
\put(7.95,-14.05){\framebox(0.1,0.1){}}

\put(8,-14){\line(1,-2){1}}

\put(9.1,-15.8){\makebox(0,0)[lc]{\(1;1..2\)}}
\put(8.8,-15.8){\makebox(0,0)[rc]{\(2\ast\)}}
\put(8.9,-16.1){\framebox(0.2,0.2){}}

\put(8,-14){\line(0,-1){4}}

\put(8.2,-17.8){\makebox(0,0)[lc]{\(2;1..3\)}}
\put(7.8,-17.8){\makebox(0,0)[rc]{\(3\ast\)}}
\put(7.9,-18.1){\framebox(0.2,0.2){}}

\put(8,-8){\line(1,-2){2}}

\put(10.2,-11.8){\makebox(0,0)[lc]{\(2;2\)}}
\put(9.9,-12.1){\framebox(0.2,0.2){}}

\put(10,-12){\line(0,-1){2}}

\put(10.2,-13.8){\makebox(0,0)[lc]{\(1;2\)}}
\put(9.95,-14.05){\framebox(0.1,0.1){}}

\put(10,-14){\line(1,-2){1}}

\put(11.2,-15.8){\makebox(0,0)[lc]{\(1;1\)}}
\put(10.9,-16.1){\framebox(0.2,0.2){}}

\put(10,-14){\line(0,-1){4}}

\put(10.2,-17.8){\makebox(0,0)[lc]{\(2;1..3\)}}
\put(9.8,-17.8){\makebox(0,0)[rc]{\(3\ast\)}}
\put(9.9,-18.1){\framebox(0.2,0.2){}}


\put(2,-2){\line(2,-1){8}}

\put(10.2,-5.8){\makebox(0,0)[lc]{\(\langle 628\rangle\)}}
\put(10.2,-6.1){\makebox(0,0)[lc]{\(=\Psi\)}}
\put(9.9,-6.1){\framebox(0.2,0.2){}}

\put(8,-6){\line(0,-1){2}}

\put(8.2,-7.8){\makebox(0,0)[lc]{\(1;2\)}}
\put(7.95,-8.05){\framebox(0.1,0.1){}}


\put(2,-2){\line(5,-2){10}}

\put(12.2,-5.8){\makebox(0,0)[lc]{\(\langle 629\rangle\)}}
\put(12.2,-6.1){\makebox(0,0)[lc]{\(=Y\)}}
\put(11.9,-6.1){\framebox(0.2,0.2){}}

\put(12,-6){\line(0,-1){2}}

{\color{red}
\put(12.2,-7.8){\makebox(0,0)[lc]{\(1;2\)}}
\put(12.2,-8.1){\makebox(0,0)[lc]{\(=Y_1\)}}
\put(11.95,-8.05){\framebox(0.1,0.1){}}

\put(12,-8){\line(1,-2){1}}

\put(13.2,-9.8){\makebox(0,0)[lc]{\(1;1\)}}
\put(12.9,-10.1){\framebox(0.2,0.2){}}
}

\put(12,-8){\line(0,-1){4}}

\put(12.2,-11.8){\makebox(0,0)[lc]{\(2;1..2\)}}
\put(11.8,-11.8){\makebox(0,0)[rc]{\(2\ast\)}}
\put(11.9,-12.1){\framebox(0.2,0.2){}}


{\color{red}
\put(13,-9.8){\oval(1.6,1.3)}
\put(13,-10.7){\makebox(0,0)[cc]{\underbar{\textbf{+21\,974\,161}}}}
}

\put(12.2,-12){\oval(2.6,1.3)}
\put(12,-12.9){\makebox(0,0)[cc]{\underbar{\textbf{-114\,936}}}}


\put(2,-2){\line(3,-1){12}}

\put(14.2,-5.8){\makebox(0,0)[lc]{\(\langle 630\rangle\)}}
\put(14.2,-6.1){\makebox(0,0)[lc]{\(=Z\)}}
\put(13.95,-6.05){\framebox(0.1,0.1){}}

\put(14,-6.0){\vector(0,-1){2}}


\end{picture}

\end{figure}
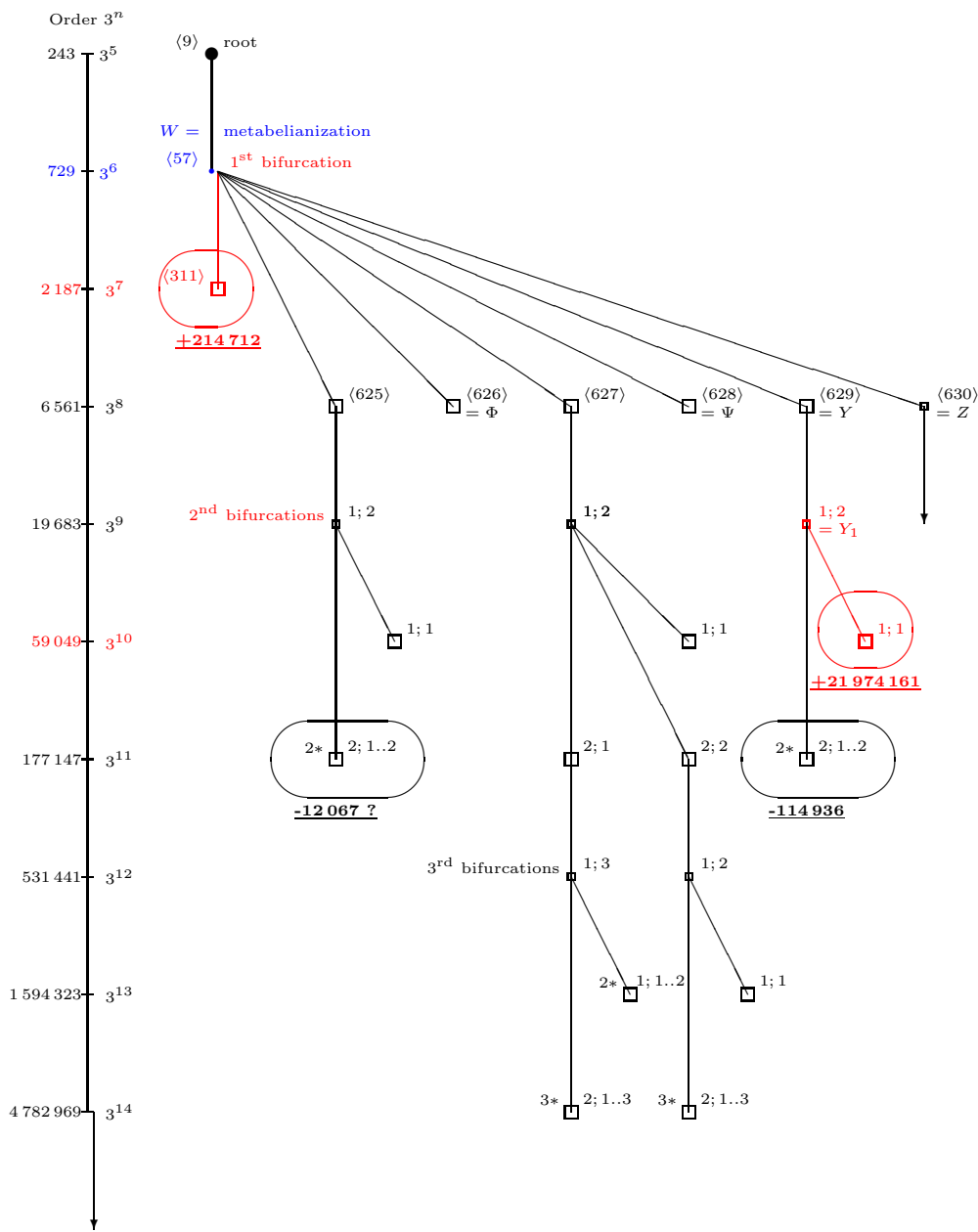

}

\noindent
After this terminating figure of a \textbf{descendant tree},
we would like to emphasize the esthetical beauty and the immense value
of this ostensive kind of graphical information
from the perspective of learning psychology.
We have developed this method of representing
ancestor-descendant relations in a detailed series of papers
\cite{Ma2013,Ma2015a,Ma2016a,Ma2018b,Ma2018c}.

\newpage


\section{Conclusion}
\label{s:Conclusion}

\noindent
In this paper,
we have proved the \textit{extremal root path property} of
the Galois group of finite \(3\)-class field towers,
stated in the Main Conjecture
\ref{cnj:ExtremalPath},
for \textit{infinite} series of parametrized Schur \(\sigma\)-groups.
The derived length of one of these series was \textit{unbounded},
whereas the other series had soluble length precisely equal to three.
Within the frame of Section F of \(3\)-capitulation types,
we provided evidence of root paths with high complexity,
and of the first \(3\)-class field tower with length \(\ell_3(K)\ge 4\).
No counter-examples against the Main Conjecture are known up to now.
Our experience suggests the extension of the Main Conjecture
to Schur \(\sigma\)-groups with abelianizations different from \((3,3)\),
and also to primes \(p\ge 5\).

In searches for \(p\)-groups with assigned Artin pattern
by means of the \textit{strategy of pattern recognition via Artin transfers},
considerable amounts of CPU time and RAM storage can be saved
by restricting the search paths to the \textit{maximal possible step sizes},
when the desired group is a Schur \(\sigma\)-group.
Frequently, the descendant number for the maximal step size
is smaller than for other step sizes,
which may be outside of the reach of current implementations already.

As an impressive application of
the technique of extremal root paths,
we succeeded in conducting a rigorous proof
of the first \(3\)-class field towers with at least four stages
over imaginary quadratic number fields.
Even when these towers should turn out to be infinite
this would be a striking novelty,
because it was believed that imaginary quadratic fields
with \(3\)-class groups of rank two
have towers with finite length.



\section{Acknowledgement}
\label{s:Acknowledgement}

\noindent
The
author gratefully acknowledges that his research was supported by the
Austrian Science Fund (FWF): P 26008-N25,
and by the Research Executive Agency of the European Union (EUREA).

Deep indebtedness is expressed to Professor Mike F. Newman
from the Australian National University in Canberra,
Australian Capital Territory.



\section{Historical remarks concerning TKTs in Section F}
\label{s:HistorySectionF}

\noindent
\textit{Imaginary} quadratic fields \(K=\mathbb{Q}(\sqrt{d})\)
with \(3\)-class group \(\mathrm{Cl}_3{K}\simeq C_3\times C_3\)
and transfer kernel type (TKT) in Section \(\mathrm{F}\)
have been detected by Brink in \(1984\)
\cite{Br}.
The absolute values of their fundamental discriminants \(d\)
set in with \(27\,156\), outside of the ranges investigated by
Scholz and Taussky in \(1934\)
\cite{SoTa},
and by Heider and Schmithals in \(1982\)
\cite{HeSm}.
However, the computational results in Brink's Thesis
\cite[Appendix A, pp. 96--113]{Br}
were unknown to us until we got a copy via ProQuest in \(2006\).
Their actual extent is not mentioned explicitly in the official paper
\cite{BrGo}
by Brink and his academic advisor Gold.
Therefore, we previously believed to have the priority in discovering the discriminant
\(d=-27\,156\) of a field \(K\) with type \(\mathrm{F}.11\) in \(1989\)
\cite[Tbl., p. 84]{Ma1991},
and the discriminants \(d=-31\,908\), \(-67\,480\), \(-124\,363\)
of fields \(K\) with types \(\mathrm{F}.12\), \(\mathrm{F}.13\), \(\mathrm{F}.7\) in \(2003\)
\cite[Tbl. 3, p. 497]{Ma2012},
all of them with second \(3\)-class groups \(\mathrm{G}_3^2{K}\) of coclass \(4\).
In \(2006\), it turned out that our claim must be restricted to \(d=-124\,363\),
which after nearly \(20\) years eventually provided the first example for type \(\mathrm{F}.7\),
called the unique undiscovered type by Brink
\cite[\S\ 7.2, p. 91]{Br}.

It required further \(10\) years until we had the courage to study the \(3\)-class tower
of number fields with transfer kernel type \(\mathrm{F}\),
based on abelian type invariants of second order, as developed in
\cite{Ma2015b,Ma2017},
and inspired by repeated sparkling ideas of Professor Mike F. Newman in 2013 and 2017.

As opposed to coclass \(4\),
we can definitely claim priority in discovering the discriminant
\(d=-423\,640\) of a \textit{complex} quadratic field \(K=\mathbb{Q}(\sqrt{d})\) with type \(\mathrm{F}.12\) in \(2010\)
\cite[Tbl. 3, p. 497]{Ma2012},
and the discriminants \(d=-1\,677\,768\), \(-2\,383\,059\), \(-4\,838\,891\)
of fields \(K\) with types \(\mathrm{F}.7\), \(\mathrm{F}.13\), \(\mathrm{F}.11\) in \(2016\),
all of them with second \(3\)-class groups \(\mathrm{G}_3^2{K}\) of coclass \(6\).

Similarly,
we were the first who found the discriminant
\(d=8\,321\,505\) of a \textit{real} quadratic field \(K=\mathbb{Q}(\sqrt{d})\) with type \(\mathrm{F}.13\) in \(2010\)
\cite[Tbl. 4, p. 498]{Ma2012},
and the discriminants \(d=10\,165\,597\), \(22\,937\,941\), \(66\,615\,244\)
of fields \(K\) with types \(\mathrm{F}.7\), \(\mathrm{F}.12\), \(\mathrm{F}.11\) in \(2016\)
\cite[Tbl. 4, p. 1291]{Ma2016b},
\cite{Ma2016c},
all of which possess second \(3\)-class groups \(\mathrm{G}_3^2{K}\) of coclass \(4\).

\newpage


\section{Personal historical remarks}
\label{s:PersonalHistory}

\noindent
When I began to investigate
the \(3\)-\textit{capitulation} in unramified cyclic cubic extensions \(E_i/K\)
of \textit{imaginary} quadratic number fields \(K\) with \(3\)-class rank two
(it was in autumn \(1989\), that is, thirty years ago),
I had a sound foundation in algebraic number theory
and elements of class field theory,
but only a very basic knowledge in \(p\)-group theory.
My usage of the concept \textit{metabelian group of class two}
in the resulting paper
\cite[pp. 73, 79, 84, 86]{Ma1991}
shows that I was thinking within Helmut Hasse's scope of
\textit{metabelian groups with multiple stages},
which would be called \textit{solvable} or \textit{soluble} nowadays.
With \textit{class two} I rather meant \textit{two stages},
which is simply \textit{metabelian}
(non-abelian with abelian commutator subgroup, or derived length two)
in modern mathematical language.

In
\cite{Ma1991},
my view of
the \textit{capitulation type} \(\varkappa(K)\) of imaginary quadratic fields \(K\)
was just a nice new invariant,
in addition to class number \(h_K\) and \(p\)-class rank \(\mathrm{rk}_p(K)\),
for various primes \(p\).
Although I knew that Scholz and Taussky
\cite{SoTa}
intended to derive information about
the metabelian \textit{second \(3\)-class group} \(M:=\mathrm{Gal}(\mathrm{F}_3^2(K)/K)\)
or even the pro-\(3\) Galois group \(G:=\mathrm{Gal}(\mathrm{F}_3^\infty(K)/K)\)
of the entire \(3\)-\textit{class field tower} of \(K\)
from the type \(\varkappa(K)\) of \(3\)-capitulation,
I could not get a firm grasp of their claim
that the \textit{annihilator ideal} of the main commutator \(s_2=\lbrack y,x\rbrack\) of \(M=\langle x,y\rangle\)
together with two Schreier polynomials
determines the metabelian group \(M\) uniquely (see
\cite{Ma2018d}).

My poor horizon concerning the systematic treatment of \(p\)-groups
did not extend until Aissa Derhem (in December \(2001\)) drew my attention 
to Nebelung's doctoral thesis
\cite{Ne1,Ne2}.
Suddenly I was able to remember numerous finite metabelian \(3\)-groups
according to their characteristic positions in \textit{descendant trees},
more precisely in \textit{coclass trees}
with stable difference \(\mathrm{cc}=\mathrm{lo}-\mathrm{cl}\)
between logarithmic order and nilpotency class,
and I started a big revival of the capitulation problem
for \textit{imaginary} quadratic number fields \(K\)
with \(3\)-class group \(\mathrm{Cl}_3(K)\simeq (3,3)\)
in \(2003\),
now always with the higher goal to identify the corresponding metabelian group
\(M=\mathrm{Gal}(\mathrm{F}_3^2(K)/K)\).
In \(2006\), I extended my investigations to \textit{real} quadratic fields \(K\),
for which the capitulation problem was nearly unknown
(except for \(5\) cases in the paper of Heider and Schmithals
\cite{HeSm}).
This scientific phase got its coronation in a series of papers
\cite{Ma2012,Ma2010,Ma2014,Ma2013}
where the second \(3\)-class groups \(M\) were determined
for all \(4596\) quadratic fields \(K\) with discriminants \(-10^6<d_K<10^7\)
and \(\mathrm{Cl}_3(K)\simeq (3,3)\).
It was the final realization of a dream
I had twenty years earlier
\cite[p. 77]{Ma1991}.

However, in spite of my success in determining
the second stage \(M\) of \(3\)-class field towers in the year \(2010\),
I began to realize that neither Nebelung's theory
\cite{Ne1}
nor the original work of Scholz and Taussky
\cite{SoTa}
provided the required techniques for answering the problem
concerning the length \(\ell_3(K)\) of the tower
for the mysterious imaginary quadratic fields \(K\)
with capitulation types in Section E,
which had been raised by Brink and Gold
\cite{Br,BrGo}
without any conclusive decision.
Scholz and Taussky had proved \(\ell_3(K)=2\) for imaginary quadratic fields \(K\)
with capitulation types in Section D
\cite[p. 39]{SoTa}.
In some hasty remarks they also claimed a two-stage tower
for other types, in particular those in Section E
\cite[pp. 20, 41]{SoTa}.
In \(1992\), Franz Lemmermeyer drew my attention to
the doubts of Brink and Gold about these remarks by Scholz and Taussky.

In \(2012\), a chain of lucky coincidences enabled the rigorous solution
of this annoying problem.
The dissertation of Tobias Bembom acquainted me with \textit{Schur \(\sigma\)-groups}
and unpublished drafts by Boston, Bush and Hajir
(published much later in
\cite{BBH})
about the role of these groups
as pro-\(3\) Galois groups \(G=\mathrm{Gal}(\mathrm{F}_3^\infty(K)/K)\)
of \(3\)-class field towers of imaginary quadratic fields \(K\).
In August \(2012\), I had opportunity to meet these three authors at a conference in Vienna.
In a discussion of only one or two hours,
Boston, Bush and myself succeeded in strictly proving \textit{exact length} \(\ell_3(K)=3\)
for imaginary quadratic fields \(K\) with capitulation types in Section E.
The reason why Brink and Gold did not succeed in the definite exclusion of length two
in \(1987\) was their lack of knowledge about Schur \(\sigma\)-groups, although
\cite{Sh}
appeared in \(1964\) and
\cite{KoVe}
in \(1975\).

I shall not become tired of pointing out again and again
the fundamental importance of the Shafarevich cohomology criterion
for the relation rank \(d_2(G)\) of the group \(G=\mathrm{Gal}(\mathrm{F}_3^\infty(K)/K)\).
In discussions at mathematical conferences,
I repeatedly realized that even experts in number theory
only know the famous results of Golod and Shafarevich on infinite class field towers,
but are not aware of the Shafarevich criterion
\cite{Sh}.

\newpage



\end{document}